\documentclass{alggeom}

\usepackage{mathrsfs}
\usepackage{amscd}
\usepackage{amsmath}
\usepackage{amssymb}
\usepackage{latexsym}
\usepackage{comment}
\usepackage{amscd}
\usepackage{wasysym}  
\usepackage{tikz}
\usetikzlibrary{matrix,arrows}
\usepackage{tikz-cd}
\usepackage{appendix}
\usepackage[nice]{nicefrac}
\usepackage{fix-cm} 
\usepackage{comment}
\usepackage{amscd}
\usepackage{wasysym}  
\usepackage{tikz}
\usetikzlibrary{matrix,arrows}
\usepackage{tikz-cd}
\usepackage{appendix}
\usepackage[nice]{nicefrac}

\usepackage{marvosym}
\usepackage{tikz}
\usepackage[all,cmtip]{xy}
\usetikzlibrary{matrix,arrows}

\usepackage[pagebackref,hyperindex,linktocpage=true]{hyperref}
\hypersetup{
    colorlinks,
    linkcolor={blue!40!black},
    citecolor={blue!40!black},
    urlcolor={blue!40!black}
}

\usepackage{color}

\usetikzlibrary{decorations.markings}

\makeatletter
\tikzcdset{
  open/.code     = {\tikzcdset{hook, circled};},
  closed/.code   = {\tikzcdset{hook, slashed};},
  open'/.code    = {\tikzcdset{hook', circled};},
  closed'/.code  = {\tikzcdset{hook', slashed};},
  circled/.code  = {\tikzcdset{markwith = {\draw (0,0) circle (.375ex);}};},
  slashed/.code  = {\tikzcdset{markwith = {\draw[-] (-.4ex,-.4ex) -- (.4ex,.4ex);}};},
  markwith/.code ={
    \pgfutil@ifundefined%
    {tikz@library@decorations.markings@loaded}%
    {\pgfutil@packageerror{tikz-cd}{You need to say %
      \string\usetikzlibrary{decorations.markings} to use arrows with markings}{}}{}%
    \pgfkeysalso{/tikz/postaction = {
      /tikz/decorate,
      /tikz/decoration={markings, mark = at position 0.5 with {#1}}}
    }
  },
}
\makeatother

\renewcommand{\geq}{\geqslant}
\renewcommand{\leq}{\leqslant}

\newtheorem{thm}{Theorem}[section]

\newtheorem{propo}[thm]{Proposition}
\newtheorem{lem}[thm]{Lemma}
\newtheorem{sublem}[thm]{Sublemma}
\newtheorem{lem-def}[thm]{Lemma-Definition}
\newtheorem{cor}[thm]{Corollary}
\newtheorem{conject}[thm]{Conjecture}
\newtheorem{propert}[thm]{Properties}
\newtheorem{observ}[thm]{Observation}

\newtheorem{fac}[thm]{Fact}

\theoremstyle{definition}
\newtheorem*{ack}{Acknowledgement}
\newtheorem*{con}{Conventions}

\newtheorem{dfn}[thm]{Definition}
\newtheorem{const}[thm]{Construction}

\newtheorem{quest}[thm]{Question}
\newtheorem{expec}[thm]{Expectation}

\theoremstyle{remark}
\newtheorem{rmk}[thm]{Remark}
\newtheorem{ex}[thm]{Example}

\numberwithin{equation}{section}

\newcommand{\nc}{\newcommand}

\nc{\theo}{\begin{thm}} \nc{\xtheo}{\end{thm}}
\nc{\prop}{\begin{propo}} \nc{\xprop}{\end{propo}}
\nc{\lemm}{\begin{lem}} \nc{\xlemm}{\end{lem}}
\nc{\sublemm}{\begin{sublem}} \nc{\xsublemm}{\end{sublem}}
\nc{\lemmdefi}{\begin{lem-def}} \nc{\xlemmdefi}{\end{lem-def}}
\nc{\coro}{\begin{cor}} \nc{\xcoro}{\end{cor}}
\nc{\conj}{\begin{conject}} \nc{\xconj}{\end{conject}}
\nc{\proper}{\begin{propert}} \nc{\xproper}{\end{propert}}
\nc{\obse}{\begin{observ}} \nc{\xobse}{\end{observ}}
\nc{\ques}{\begin{quest}} \nc{\xques}{\end{quest}}

\nc{\fact}{\begin{fac}} \nc{\xfact}{\end{fac}}
\nc{\expe}{\begin{expec}} \nc{\xexpe}{\end{expec}}

\nc{\ackn}{\begin{ack}} \nc{\xackn}{\end{ack}}
\nc{\exam}{\begin{ex}} \nc{\xexam}{\end{ex}}
\nc{\rema}{\begin{rmk}} \nc{\xrema}{\end{rmk}}
\nc{\defi}{\begin{dfn}} \nc{\xdefi}{\end{dfn}}
\nc{\abst}{\begin{abs}} \nc{\xabst}{\end{abs}}

\nc{\pf}{\begin{proof}} \nc{\xpf}{\end{proof}}

\nc{\on}{\operatorname}
\nc{\fraka}{{\mathfrak a}} \nc{\bba}{{\mathbf a}}
\nc{\frakb}{{\mathfrak b}}
\nc{\frakc}{{\mathfrak c}}
\nc{\frakd}{{\mathfrak d}}
\nc{\frake}{{\mathfrak e}}
\nc{\frakf}{{\mathfrak f}}
\nc{\frakg}{{\mathfrak g}}
\nc{\frakh}{{\mathfrak h}}
\nc{\fraki}{{\mathfrak i}}
\nc{\frakj}{{\mathfrak j}}
\nc{\frakk}{{\mathfrak k}}
\nc{\frakl}{{\mathfrak l}}
\nc{\frakm}{{\mathfrak m}}
\nc{\frakn}{{\mathfrak n}}
\nc{\frako}{{\mathfrak o}}
\nc{\frakp}{{\mathfrak p}}
\nc{\frakq}{{\mathfrak q}}
\nc{\frakr}{{\mathfrak r}}
\nc{\fraks}{{\mathfrak s}}
\nc{\frakt}{{\mathfrak t}}
\nc{\fraku}{{\mathfrak u}}
\nc{\frakv}{{\mathfrak v}}
\nc{\frakw}{{\mathfrak w}}
\nc{\frakx}{{\mathfrak x}}
\nc{\fraky}{{\mathfrak y}}
\nc{\frakz}{{\mathfrak z}}
\nc{\frakA}{{\mathfrak A}}
\nc{\frakB}{{\mathfrak B}}
\nc{\frakC}{{\mathfrak C}}
\nc{\frakD}{{\mathfrak D}}
\nc{\frakE}{{\mathfrak E}}
\nc{\frakF}{{\mathfrak F}}
\nc{\frakG}{{\mathfrak G}}
\nc{\frakH}{{\mathfrak H}}
\nc{\frakI}{{\mathfrak I}}
\nc{\frakJ}{{\mathfrak J}}
\nc{\frakK}{{\mathfrak K}}
\nc{\frakL}{{\mathfrak L}}
\nc{\frakM}{{\mathfrak M}}
\nc{\frakN}{{\mathfrak N}}
\nc{\frakO}{{\mathfrak O}}
\nc{\frakP}{{\mathfrak P}}
\nc{\frakQ}{{\mathfrak Q}}
\nc{\frakR}{{\mathfrak R}}
\nc{\frakS}{{\mathfrak S}}
\nc{\frakT}{{\mathfrak T}}
\nc{\frakU}{{\mathfrak U}}
\nc{\frakV}{{\mathfrak V}}
\nc{\frakW}{{\mathfrak W}}
\nc{\frakX}{{\mathfrak X}}
\nc{\frakY}{{\mathfrak Y}}
\nc{\frakZ}{{\mathfrak Z}}
\nc{\bbA}{{\mathbb A}}
\nc{\bbB}{{\mathbb B}}
\nc{\bbC}{{\mathbb C}}
\nc{\bbD}{{\mathbb D}}
\nc{\bbE}{{\mathbb E}}
\nc{\bbF}{{\mathbb F}} \nc{\bbf}{{\mathbf f}}
\nc{\bbG}{{\mathbb G}}
\nc{\bbH}{{\mathbb H}}
\nc{\bbI}{{\mathbb I}}
\nc{\bbJ}{{\mathbb J}}
\nc{\bbK}{{\mathbb K}}
\nc{\bbL}{{\mathbb L}}
\nc{\bbM}{{\mathbb M}}
\nc{\bbN}{{\mathbb N}}
\nc{\bbO}{{\mathbb O}}
\nc{\bbP}{{\mathbb P}}
\nc{\bbQ}{{\mathbb Q}}
\nc{\bbR}{{\mathbb R}}
\nc{\bbS}{{\mathbb S}}
\nc{\bbT}{{\mathbb T}}
\nc{\bbU}{{\mathbb U}}
\nc{\bbV}{{\mathbb V}}
\nc{\bbW}{{\mathbb W}}
\nc{\bbY}{{\mathbb Y}}
\nc{\bbZ}{{\mathbb Z}}
\nc{\calA}{{\mathcal A}}
\nc{\calB}{{\mathcal B}}
\nc{\calC}{{\mathcal C}}
\nc{\calD}{{\mathcal D}}
\nc{\calE}{{\mathcal E}}
\nc{\calF}{{\mathcal F}}
\nc{\calG}{{\mathcal G}}
\nc{\calH}{{\mathcal H}}
\nc{\calI}{{\mathcal I}}
\nc{\calJ}{{\mathcal J}}
\nc{\calK}{{\mathcal K}}
\nc{\calL}{{\mathcal L}}
\nc{\calM}{{\mathcal M}}
\nc{\calN}{{\mathcal N}}
\nc{\calO}{{\mathcal O}}
\nc{\calP}{{\mathcal P}}
\nc{\calQ}{{\mathcal Q}}
\nc{\calR}{{\mathcal R}}
\nc{\calS}{{\mathcal S}}
\nc{\calT}{{\mathcal T}}
\nc{\calU}{{\mathcal U}}
\nc{\calV}{{\mathcal V}}
\nc{\calW}{{\mathcal W}}
\nc{\calX}{{\mathcal X}}
\nc{\calY}{{\mathcal Y}}
\nc{\calZ}{{\mathcal Z}}

\nc{\scrA}{{\mathscr A}}
\nc{\scrE}{{\mathscr E}}
\nc{\scrR}{{\mathscr R}}

\nc{\Bmu}{\mbox{$\raisebox{-0.59ex}{$l$}\hspace{-0.18em}\mu\hspace{-0.88em}\raisebox{-0.98ex}{\scalebox{2}{$\color{white}.$}}\hspace{-0.416em}\raisebox{+0.88ex}{$\color{white}.$}\hspace{0.46em}$}{}}

\nc{\bnu}{{\bar{ \nu}}}

\nc{\olO}{\bar{\calO}}

\nc{\al}{{\alpha}} 
\nc{\be}{{\beta}}
\nc{\ga}{{\gamma}} \nc{\Ga}{{\Gamma}}
 \nc{\hGa}{\hat{\Gamma}}
\nc{\ve}{{\varepsilon}} 
\nc{\la}{{\lambda}} \nc{\La}{{\Lambda}}
\nc{\om}{\omega} \nc{\Om}{\Omega} 
\nc{\sig}{{\sigma}} \nc{\Sig}{{\Sigma}}

\nc{\tnb}{\psi_{\rm tame}}
\nc{\oM}{\overline{{M}}}
\nc{\op}{{\on{op}}}
\nc{\ad}{{\on{ad}}}
\nc{\alg}{{\on{alg}}}
\nc{\Ad}{{\on{Ad}}}
\nc{\Adm}{{\on{Adm}}} 
\nc{\Aut}{{\on{Aut}}}
\nc{\Bun}{{\on{Bun}}}
\nc{\cha}{{\on{char}}}
\nc{\der}{{\on{der}}}
\nc{\Der}{{\on{Der}}}
\nc{\Rad}{{\on{Rad}}}
\nc{\diag}{{\on{diag}}}
\nc{\End}{{\on{End}}}
\nc{\Fl}{{\calF\!\ell}}
\nc{\Tr}{{\on{Transp}}}
\nc{\TR}{{\calT\!\calR}}
\nc{\Gal}{{\on{Gal}}}
\nc{\Gr}{{\on{Gr}}}
\nc{\rH}{{\on{H}}}
\nc{\Hom}{{\on{Hom}}}
\nc{\IC}{{\on{IC}}}
\nc{\id}{{\on{id}}}
\nc{\Id}{{\on{Id}}}
\nc{\ind}{{\on{ind}}}
\nc{\Ind}{{\on{Ind}}}
\nc{\Lie}{{\on{Lie}}}
\nc{\Pic}{{\on{Pic}}}
\nc{\pr}{{\on{pr}}}
\nc{\Res}{{\on{Res}}}
\nc{\res}{{\on{res}}} \nc{\Sat}{{\on{Sat}}}
\nc{\s}{{\on{sc}}}
\nc{\drv}{{\on{der}}}
\nc{\sgn}{{\on{sgn}}}
\nc{\Spec}{{\on{Spec}}}
\nc{\Quot}{{\on{Quot}}}
\nc{\Spf}{\on{Spf}} 
\nc{\Sph}{\on{Sph}}
\nc{\St}{{\on{St}}}
\nc{\tr}{{\on{tr}}}
\nc{\Mod}{{\mathsf{Mod}}}
\nc{\Hilb}{{\on{Hilb}}} 
\nc{\Ext}{{\on{Ext}}} 
\nc{\vs}{{\on{Vec}}}
\nc{\ev}{{\on{ev}}}
\nc{\nO}{{\breve{\calO}}}
\nc{\tS}{{\tilde{S}}}
\nc{\spe}{{\on{sp}}}
\nc{\loc}{{\on{loc}}}
\nc{\Sym}{{\on{Sym}}}
\nc{\Cone}{{\on{C}}}
\nc{\syn}{{\on{syn}}}
\nc{\reg}{{\on{reg}}}
\nc{\colim}{{\on{colim}}}
\nc{\Norm}{{\on{N}}}

\nc{\nscrR}{{\mathscr{R}^{\on{nr}}}}

\nc{\GL}{{\on{GL}}}
\nc{\U}{{\on{U}}}
\nc{\Gl}{\on{Gl}} 
\nc{\GSp}{{\on{GSp}}}
\nc{\gl}{{\frakg\frakl}}
\nc{\SL}{{\on{SL}}} 
\nc{\SU}{{\on{SU}}} 
\nc{\SO}{{\on{SO}}}
\nc{\PGL}{{\on{PGL}}}

\nc{\Conv}{{\on{Conv}}}
\nc{\Rep}{{\on{Rep}}}
\nc{\Dom}{{\on{Dom}}}
\nc{\red}{{\on{red}}}
\nc{\act}{{\on{act}}}
\nc{\nr}{{\on{nr}}}
\nc{\ctf}{{\on{ctf}}}

\nc{\str}{{\on{-}}} 
\nc{\os}{{\bar{s}}}
\nc{\oeta}{{\bar{\eta}}}

\nc{\hookto}{\hookrightarrow}
\nc{\longto}{\longrightarrow}
\nc{\leftto}{\leftarrow}
\nc{\onto}{\twoheadrightarrow}
\nc{\lonto}{\twoheadleftarrow}

\nc{\uG}{{\underline{G}}}
\nc{\uA}{{\underline{A}}}
\nc{\uS}{{\underline{S}}}
\nc{\uT}{{\underline{T}}}
\nc{\uM}{{\underline{M}}}
\nc{\uP}{{\underline{P}}}
\nc{\uB}{{\underline{B}}}
\nc{\uN}{{\underline{N}}}

\nc{\ucG}{{\underline{\calG}}}
\nc{\ucA}{{\underline{\calA}}}
\nc{\ucS}{{\underline{\calS}}}
\nc{\ucT}{{\underline{\calT}}}
\nc{\ucalM}{{\underline{\calM}}}
\nc{\ucP}{{\underline{\calP}}}
\nc{\ucalN}{{\underline{\calN}}}

\nc{\bF}{{\breve{F}}}

\nc{\oFl}{{\overline{\Fl}}} 
\nc{\tGr}{{\tilde{\Gr}}}
\nc{\cGr}{\calG\! r}
\nc{\oGr}{\overline{\on{Gr}}} 
\nc{\ocGr}{\overline{\calG\! r}}
\nc{\co}{{\colon}}
\nc{\sch}[1]{(Sch/{#1})}
\nc{\HypLoc}[1]{HypLoc({#1})}

\nc{\ohtimes}{\stackrel{!}{\otimes}}
\nc{\boxtilde}{\widetilde{\boxtimes}}
\nc{\vstar}{{\varhexstar}}

\nc{\Div}{\on{Div}}
\nc{\Sht}{\on{Sht}}
\nc{\Frob}{\on{Frob}}

\nc{\x}{\times}
\nc{\bsl}{\backslash}
\nc{\algQl}{{\bar{\bbQ}_\ell}}
\nc{\sF}{{\bar{F}}}
\nc{\nF}{{\breve{F}}}
\nc{\nW}{{W^{\on{nr}}}}
\nc{\sk}{{\bar{k}}}
\nc{\cont}{\on{c}}
\nc{\Supp}{\on{Supp}}
\nc{\blt}{\bullet}  
\nc{\dom}{\on{dom}}
\nc{\scon}{{\on{sc}}} 
\nc{\Affine}{\on{Aff}} 
\nc{\nscrA}{\mathscr{A}^{\on{nr}}} 
\nc{\nfraka}{{\bbf^{\on{nr}}}}
\nc{\ran}{{\rangle}}
\nc{\lan}{{\langle}}
\nc{\tF}{{\tilde{F}}}
\nc{\sS}{{\bar{S}}}
\nc{\LG}{{^\text{L}\hspace{-0.04cm}G}}
\nc{\LL}{{^\text{L}\hspace{-0.07cm}L}}
\nc{\et}{{\text{\rm \'et}}}
\nc{\inv}{{\on{inv}}}
\nc{\Hecke}{{\on{Hecke}}}
\nc{\Isom}{{\on{Isom}}}
\nc{\oSht}{{\overline{\on{Sht}}}}
\nc{\umu}{{\underline \mu}}
\nc{\AIJ}{{\calO_X[{\scriptstyle{\calI\over \calJ}}]}}
\nc{\Proj}{{\on{Proj}}}
\nc{\Bl}{{\on{Bl}}}
\nc{\Stab}{{\on{Stab}}}
\nc{\cl}{{\on{cl}}}

\nc{\Pos}{{\on{Pos}}}
\nc{\Sets}{{\on{Sets}}}
\nc{\AffSch}{{\on{AffSch}}}
\nc{\Groups}{{\on{Groups}}}
\nc{\Gpds}{{\on{Groupoids}}}
\nc{\Sch}{{\on{Sch}}}
\nc{\fl}{{\on{flat}}}

\nc{\pot}[1]{ [\hspace{-0,5mm}[ {#1} ]\hspace{-0,5mm}] }
\nc{\rpot}[1]{ (\hspace{-0,7mm}( {#1} )\hspace{-0,7mm}) }

\nc{\defined}{\hspace{0.1cm}\stackrel{\text{\tiny \rm def}}{=}\hspace{0.1cm}}

\usepackage{moresize}
\usepackage{graphicx}

\newcommand{\Z}{\mathbb{Z}}
\newcommand{\N}{\mathbb{N}}
\newcommand{\bk}{\Bbbk}
\newcommand{\bG}{\mathbf{G}}
\newcommand{\bB}{\mathbf{B}}
\newcommand{\bT}{\mathbf{T}}
\newcommand{\bU}{\mathbf{U}}
\newcommand{\bbX}{\mathbb{X}}
\newcommand{\fR}{\mathfrak{R}}
\newcommand{\simto}{\xrightarrow{\sim}}

\newcommand{\bu}{\mathbf{u}}
\newcommand{\bg}{\mathbf{g}}
\newcommand{\Spr}{\widetilde{\mathcal{N}}}
\newcommand{\QCoh}{\mathsf{QCoh}}
\newcommand{\Coh}{\mathsf{Coh}}
\newcommand{\aff}{\mathrm{aff}}

\begin{document}

\title
%[Proj schemes] % ’optional short form; for the running head’
{On multi-graded Proj schemes}

%\shortauthors{}

\author{Arnaud Mayeux}
\email{arnaud.mayeux@mail.huji.ac.il}
\address{Einstein Institute of Mathematics,
Edmond J. Safra Campus,
The Hebrew University of Jerusalem,
Givat Ram. Jerusalem, 9190401, Israel}

\author{Simon Riche}
\email{simon.riche@uca.fr}
\address{Universit\'e Clermont Auvergne, CNRS, LMBP, F-63000 Clermont-Ferrand, France}

\classification{14A05,14A15}
\keywords{Proj construction, Brenner--Schröer, graded algebra, graded rings, potions of graded rings, twisting sheaves, algebraic dilatations, flag variety, Springer resolution.}

\begin{abstract}
We review the construction (due to Brenner--Schr\"oer) of the Proj scheme associated with a ring graded by a finitely generated abelian group. This construction generalizes the well-known Grothendieck Proj construction for $\bbN$-graded rings; we extend some classical results (in particular, regarding quasi-coherent sheaves on such schemes) from the $\bbN$-graded setting to this general setting, and prove new results that make sense only in the general setting of Brenner--Schröer. Finally, we show that flag varieties of reductive groups, as well as some vector bundles over such varieties attached to representations of a Borel subgroup, can be naturally interpreted in this formalism. 
\end{abstract}

\maketitle 

\begin{center}
\today
\end{center}

\tableofcontents

%\ackn 
%A.M. is supported by ISF grant 1577/23.
%This project has received funding from the European Research
%Council (ERC) under the European Union’s Horizon 2020 research and innovation programme
%(grant agreement No 101002592).
%\xackn

%%%%%%%%%%%%%%%%%%%%%%%%%%%%%%%%%%%%%%%%
\section{Introduction}
%%%%%%%%%%%%%%%%%%%%%%%%%%%%%%%%%%%%%%%%

%---------------------------------------
\subsection{Proj schemes}
%---------------------------------------

As part of his refoundation of Algebraic Geometry during the second part of the twentieth century, A.~Grothendieck introduced the Proj scheme of an $\bbN$-graded ring. This generalization of the notion of projective space is natural in the framework of his theory of schemes. Grothendieck's Proj construction has since been studied, taught, treated and used many thousands of times. 
During the first years of the third millennium, H.~Brenner and S.~Schr\"oer generalized Grothendieck's Proj construction to rings graded by arbitrary finitely generated abelian groups. Surprisingly, this very nice generalization still seems to be relatively confidential,
% (e.g. \cite{ab} does not cite \cite{BS03}), 
although it has been recently used and studied in some references, including \cite{Za19, KSU21, MR24}. In this paper, we provide a slightly different perspective on the Brenner--Schr\"oer Proj construction, and prove a couple of new results on this construction. In particular, we explain a generalization in this setting of Serre's classical description of coherent sheaves on a projective space (or, more generally, a Proj scheme) as a Serre quotient of a category of graded modules over a graded ring, see e.g.~\cite[\href{https://stacks.math.columbia.edu/tag/01YR}{Tag 01YR}]{stacks-project}.

To illustrate this theory,
we also explain how to describe flag varieties of reductive algebraic groups over algebraically closed fields, and vector bundles over such flag varieties associated with some representations of a Borel subgroup (including the Springer resolution) as the Proj scheme of a ring graded by characters of a maximal torus. The description of flag varieties as Proj schemes of $\N$-graded rings is classical, see e.g.~\cite{wang}, but it requires a choice of a strictly dominant weight. In the multi-graded setting one does \emph{not} make any such choice, hence gets a more canonical construction, which is sometimes useful. In fact, versions of this construction appear implicitly or explicitly (but without proper references) in various constructions in Geometric Representation Theory, including~\cite{ab,abg,arider}, and it was one of our motivations for this work to make this construction completely explicit and rigorous. Note that \cite[Exercise 14.16]{MS05} indicates a special case of this construction.

%We use this general Proj construction to interpret Flag Varieties and Springer Resolutions as Proj of some graded rings. We expect this will shed light and more attention on Brenner-Schröer construction. 

\begin{rmk}
 Other Proj-like constructions appear in the literature, often in relation with various forms of ``toric geometry,'' and some versions of the results explained below can be found in such frameworks (see e.g.~\cite{rohrer}). What we find particularly satisfying in the Brenner--Schr\"oer theory is its elementary, while very general, nature. 
\end{rmk}

%-------------------------------------------------
\subsection{The Brenner--Schr\"oer construction}
%-------------------------------------------------

Let us now introduce the ideas of our work more precisely. Recall that given a commutative ring $A$, 
%one associates 
its (affine) spectrum $\Spec(A)$ is defined as the set of prime ideals of $A$, endowed with a certain topology and a sheaf of rings. Most treatments of Grothendieck's Proj construction proceed by analogy with this construction,
% of the spectrum of a ring, 
namely 
%they start 
by associating to a graded ring the set of its graded prime ideals, and endowing it with a topology and a sheaf of rings after (somewhat tedious) identifications of various sets. However, 
%one can easily observe that 
Brenner--Schr\"oer proceed differently in their generalization in~\cite[\S 2]{BS03}: they start from 
%Demazure-Grothendieck's result saying 
the fact that, if $M$ is an abelian group, the datum of a grading on a commutative ring $A$ is equivalent to the datum of an action of the associated diagonalizable group scheme $\mathrm{D}_{\Spec(\Z)}(M)$ on $\Spec(A)$ (see e.g.~\cite[Exp.~I, \S 4.7.3]{sga3}), and then consider a certain open subset in a quotient ringed space constructed in terms of this action. 
%Then Brenner-Schröer uses this action to define their Proj. 
Brenner--Schr\"oer's Proj scheme is covered by some special affine open subschemes, which generalize similar special affine subschemes appearing in Grothendieck's Proj construction. In case $M=\Z$ and the negative-degree components in $A$ are $0$, this recovers Grothendieck's construction. (We will refer to this case as the ``$\N$-graded case.'')

In the present treatment of the Proj construction, we put these ``special affine subschemes'' at the heart of the Proj construction. In fact, we give them a name: ``potions.'' Potions glue together to give birth to the Proj scheme thanks to a purely algebraic statement we call ``the magic of potions.''
%  that we prove in this paper. 
To use this statement and glue potion schemes, we have to assume 
%some finiteness conditions (e.g. 
that $M$ is finitely generated. (In the $\N$-graded case, this result is part, as an auxiliary lemma, of the justification of the tedious identifications mentioned before in most treatments of the Proj construction.)
% in the $\bbN$-graded case. 
%The terminology "potion" is a novelty of the present paper. 
Note that, in general,
%outside the $\bbN$-graded case, 
%(Grothendieck's Proj scheme), 
the underlying set of the Proj scheme is \emph{not} the set of graded prime ideals of the given graded ring;
%one has to take more general type of ideals, cf. 
it can be described as a set of graded ideals (see~\cite[Remark 2.3]{BS03}), but we believe that this point of view is not really useful, and in any case it does not play any role in the present paper.
%pointless.
Note that \cite{MS05}, already mentioned, provides some specific multi-graded Proj schemes (under the name spector) by glueing, cf. \cite[Definition 10.25]{MS05}. It seems that \cite{BS03} and \cite{MS05} are independent (each does not cite the other). Some other works related to multi-graded Proj schemes are cited in \cite{BS03, Za19, KSU21, MR24}, see e.g. \cite{Ro98, Pe07}.

%-----------------------------------------------
\subsection{Contents}
%-----------------------------------------------

We now present the results contained in this document. Let as above $M$ be an abelian group and $A$ be a commutative ring endowed with a grading by $M$. In Sections~\ref{sec:potions} and~\ref{sec:projdef} we define the notions of $M$-relevant families of $A$ (following Brenner--Schröer in the case of singletons, see Definition \ref{def:rele}) and of potions associated with these families, and explain the definition of the Proj scheme $\Proj^M(A)$ of $A$ (in case $M$ is finitely generated). 
%Generalizing a similar concept of \cite{BS03}, $M$-relevant families of $A$ are certains homogeneous multiplicative subsets $S$ of $A$ (cf. Definition \ref{def:rele}). Potions of $A$, 
The potion associated with a relevant family $(f_i : i \in I)$ is the degree-$0$ part in the localization of $A$ with respect to the multiplicative subset generated by $(f_i : i \in I)$;
%$S$ cf. 
see Definition~\ref{potionsdef}. The spectra $\Spec (A_{(S)})$ are called potion schemes and are denoted $D_{\dagger} (S)$ (cf. Construction \ref{defproj}), and
the scheme $\Proj^M(A)$ is defined by glueing the potions schemes $D_{\dagger } (S)$ over all relevant families (cf. Construction \ref{defproj}). This glueing is possible by the magic of potions (cf. Proposition \ref{prop:magical}). In other words, Section \ref{sec:potions} focuses on commutative algebra around potions while Section \ref{sec:projdef} applies this material to explain the construction of Proj schemes. 

In Sections~\ref{sec:potions} and~\ref{sec:projdef}, we also prove some basic results that may be known to experts. In particular, using the concept of quasi-relevant element (Definition~\ref{defquasirel}), we prove that the Proj scheme of a tensor product of graded rings (with the product grading) is the fiber product of the Proj schemes of the rings (Proposition \ref{produitvstensoriel}). We prove some compatibility results regarding the radical of the ideal generated by relevant elements (Propositions \ref{radical} and \ref{prop:image-relevant}). In Example~\ref{exmulti}, we explain the relation between potions and dilatations of rings (as studied in~\cite{stacks-project, MRR20, Ma24, Ma23}).  We also show that the Proj construction is functorial (Proposition \ref{foncto}). In \S\ref{ss:dilatations} we discuss multi-centered blowups and its relations with dilatations. 

In Section~\ref{s:examplereductive} we explain how to realize flag varieties of reductive groups, and some vector bundles attached to representations of a Borel subgroup, as the Proj scheme of a natural graded ring (without any choice of strictly dominant character). Finally, in Section~\ref{s:qcohproj} we study, given an $M$-graded commutative ring $A$, the natural functors relating the categories of $M$-graded $A$-modules and quasi-coherent sheaves on $\Proj^M(A)$. In particular, we explain how, under a certain technical assumption satisfied in many examples of interest, one can describe $\QCoh(\Proj^M(A))$ as a Serre quotient of the category of $M$-graded $A$-modules, and therefore obtain a version of Serre's celebrated theorem.

Many of these results (in particular, those concerned with quasi-coherent sheaves) are counterparts of well-known results regarding Grothendieck's Proj construction, which can be found e.g.~in~\cite{stacks-project}.

%We generalize several other results that are classical for $\bbN$-Proj but do not appear in \cite{BS03}, for example we study quasi-coherent sheaves on Proj and twisting sheaves (Section \ref{s:qcohproj}).

%A motivation for this paper is \cite{ab}: details, etc.

%-----------------------------------------------------------------
\subsection{Relations between various constructions for graded rings}
%-----------------------------------------------------------------

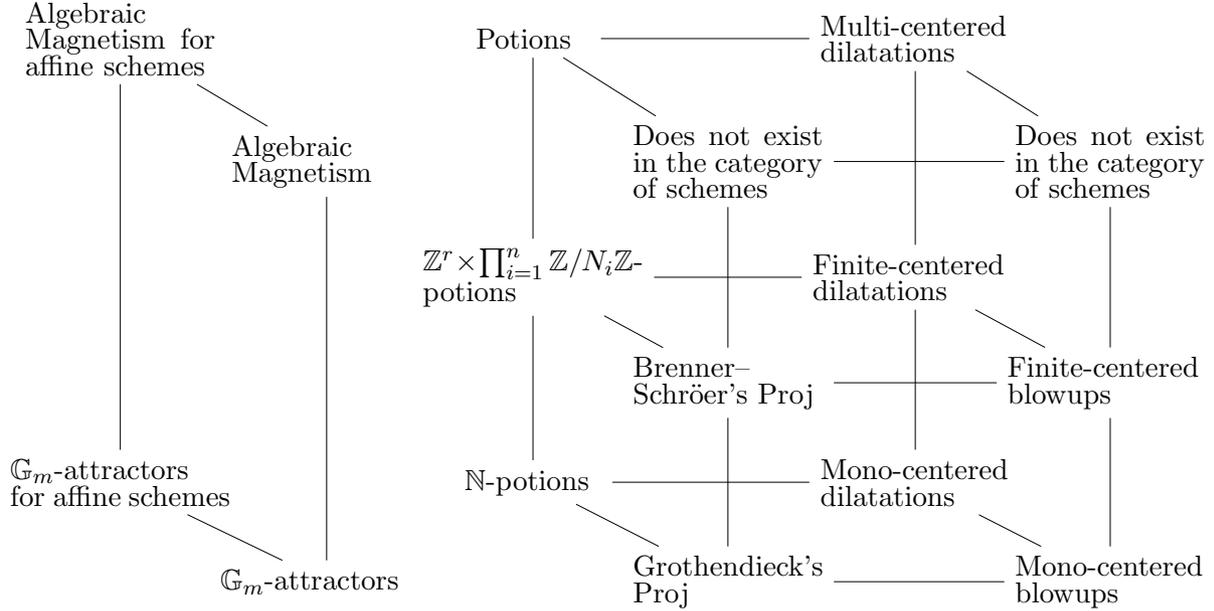
\begin{figure}[t]
\centering
{\footnotesize
\begin{tikzcd}[row sep=1.6ex, column sep=-3ex]  \begin{minipage}{2.5cm}{\normalsize Algebraic Magnetism for affine schemes}\end{minipage} \ar[dddddddd, dash]\arrow[rdd, dash] & & 
\begin{minipage}{1.5 cm}{{ {\normalsize  Potions}}} \end{minipage}\arrow[rr,dash] \arrow[dddd, dash] \arrow[rdd, dash]& &\begin{minipage}{2.5cm}{\normalsize Multi-centered dilatations}\end{minipage}  \arrow[dddd, dash] \arrow[rdd, dash]\\ 
& & & \\&
  \begin{minipage}{2.5cm}{\normalsize Algebraic Magnetism}\end{minipage} \ar[dddddddd, dash] & & \begin{minipage}{2.5cm}{\normalsize Does not exist in the category of schemes}\end{minipage} \arrow[rr,dash] \arrow[dddd, dash]& & \begin{minipage}{2.5cm}{\normalsize Does not exist in the category of schemes}\end{minipage}\arrow[dddd, dash]   \\
 & & &\\
 & &\begin{minipage}{2.9cm}{\normalsize $\bbZ^r \times \prod_{i=1}^{n} \bbZ / N_i \bbZ $-potions} \end{minipage} \arrow[rr,dash] \arrow[dddd, dash] \arrow[rdd, dash]& & \begin{minipage}{2.7cm}{\normalsize Finite-centered  dilatations}\end{minipage}  \arrow[dddd, dash] \arrow[rdd, dash]& \\
  &&&\\
& & &\begin{minipage}{2.5cm}{\normalsize Brenner--Schröer's Proj} \end{minipage} \arrow[rr,dash] \arrow[dddd, dash] & & \begin{minipage}{2.7cm}{\normalsize Finite-centered blowups} \end{minipage}\arrow[dddd, dash] \\
 &&&\\
 \begin{minipage}{2.9cm}{\normalsize $\bbG _m$-attractors for affine schemes} \end{minipage}\arrow[rdd, dash] & &
\begin{minipage}{1.8cm}{\normalsize $\bbN $-potions }\end{minipage}\arrow[rr,dash]  \arrow[rdd, dash] & & \begin{minipage}{2.5cm}{\normalsize Mono-centered dilatations}\end{minipage}   \arrow[rdd, dash]  \\
&&&\\
  &\begin{minipage}{2.8cm}{\normalsize $\bbG_m$-attractors} \end{minipage} & & \begin{minipage}{2.5cm}{\normalsize Grothendieck's Proj }\arrow[rr, dash]\end{minipage} & & \begin{minipage}{2.5cm}{\normalsize Mono-centered blowups} \end{minipage}
\end{tikzcd}
}
\caption{Some constructions in algebraic geometry associated with gradings}
\label{fig:graded-constructions}
\end{figure}

Recently, several other basic constructions for $\bbN$-graded rings were extended to more general gradings. In particular, in~\cite{Ma23a} the theory of $\bbG_{\mathrm{m}}$-attractors was refined and extended in the setting of actions of diagonalizable group schemes, and in~\cite{Ma24} mono-centered dilatations (involving $\bbN$-gradings) were extended to multi-centered dilatations. The algebraic point of view adopted in the present paper was partly inspired by these parallel generalizations. 
The diagram in Figure~\ref{fig:graded-constructions} summarizes some relations between these graded constructions. 
Lines materialize immediate connections. If two lines are parallel (not necessarily connected), then the associated connections are also parallel in some appropriate sense. 

\begin{enumerate}
\item Formalisms at the bottom of the diagram involve mono-generated gradings ($\bbN$ or $\bbZ$), while formalisms on higher floors involve more general gradings (finitely generated at the intermediate floor and arbitrary at the top floor).
\item Formalisms in the background are affine in nature (involve rings or affine morphisms of schemes), while formalisms in the foreground involve non-affine morphisms of schemes and use as prerequisites corresponding formalisms in the background.
\item Formalisms on the left study some subschemes of a given scheme endowed with a diagonalizable action (locally, a grading on the background). Formalisms in the middle of the foreground study  certain (non-affine) quotients of some given affine schemes with a diagonalizable action, and on the background we have affine local constituents to build them. Formalisms on the right are special cases of formalisms on the middle, but they behave very specifically.
\item 
\label{it:comments-cube-4}
On the right and middle, formalisms at the top of the foreground do not exist: in the arbitrary graded case, it is not possible to glue potions as at lower floors. The reason is that if $A$ is a ring and $S$ is an arbitrary multiplicative subset of $A$, the canonical morphism $\Spec (S^{-1} A) \to \Spec (A)$ is not an open immersion of schemes in general. 

\item 
Less recent and already widely used, the reference \cite{HS04} provides another multi-graded generalization of a classical construction.
It is possible that even more $\bbN$-graded constructions in schemes theory can be extended significantly to more general gradings.
\end{enumerate}
\rema
\label{rema:topolintro} 
In relation with Comment~\eqref{it:comments-cube-4} above, in Remark \ref{rema:topol} we explain that one cannot refine the theory of schemes, in the framework of locally ringed spaces, so that $\Spec (S^{-1} A) \to \Spec (A)$ is an ``open immersion'' in this refined topology. Therefore, it seems really difficult geometrically to give meaning to $\Proj (A)$ if $A$ is graded by a non-finitely generated group in general. Note however that, algebraically, Proposition~\ref{prop:magical}\eqref{it:magic-2} holds for arbitrary graded rings. 
\xrema

%--------------------------------------------
\subsection{Acknowledgements}
%--------------------------------------------

We thank A.-M.~Aubert, R.~Bezrukavnikov, S.~Brochard, J.~P.~dos Santos, A.~Ducros, S.~Evra, P.~Scholze, T.~Wedhorn and J.~Zhang for useful discussions and/or comments related to the subject of this paper. We also thank the referee for very useful comments.

Foundational parts of this paper are formalized in Lean4 by A. Mayeux and J. Zhang \cite{MZ25}, in particular potions and Proj are implemented.
Some of the conventions (e.g. the beginning of §\ref{ss:notconv} and Definition \ref{def:deg}) in the present paper were naturally adopted during the formalization in dependent type theory \cite{MZ25}.

A.~M. is supported by ISF grant 1577/23.
This project has received funding from the European Research
Council (ERC) under the European Union’s Horizon 2020 research and innovation programme
(grant agreement No 101002592).

%--------------------------------------
\subsection{Some notations and conventions} 
\label{ss:notconv}
%--------------------------------------

All the rings considered in this paper will be tacitly assumed to be associative and unital. (In practice, all the rings we want to consider will be commutative, but we will mention this assumption when it is necessary.)

Let us make explicit our conventions on graded rings.
For this we fix an arbitrary abelian group $M$.
\defi
An $M$-graded ring is a ring $A$ endowed with a direct sum decomposition $A = \bigoplus_{m \in M } A_m$ such that $A_m \cdot A_{m'} \subset  A_{m+m'}$  for all $m , m' \in M$.
\xdefi 

\defi  
Let $A$ be an $M$-graded ring. A \emph{homogeneous} pair in $A$ is a pair $(a,m)$ such that $a \in A$, $m \in M $ and $a \in A_m$. An element $a$ in $A$ is called \emph{homogeneous} if there exists $m \in M$ such that $(a,m)$ is a homogeneous pair. The degree of a homogeneous pair $x=(a,m)$ is the element of $M$ defined as $\deg(x)=m$.
\xdefi
\rema
Note that if $a$ is homogeneous and nonzero, there exists a unique $m$ such that $(a,m)$ is a  homogeneous pair (because $A_m \cap A_{m'}= \{0\}$ for $m \ne m'$). For all $m \in M$, $(0,m)$ is a homogeneous pair. 
\xrema
\begin{con}[(Degrees of homogeneous elements)] As a convention, given a homogeneous element $a$, the notation $\deg(a)$ means that we have fixed, implicitly or explicitly, an element $\deg (a) \in M$ such that $(a,\deg(a))$ is a homogeneous pair. The sentence ``let $a$ be a homogeneous element of degree $m$'' means ``let $a $ be in $ A_m$.'' The element $0$, regarded as an element of $A$, has no degree. But, as element of $A_m$ for some $m \in M$, $0$ is an element of degree $m$. 
\end{con}

\rema   
One can easily check that the unit $1 \in A$ always belongs to $A_0$, in particular it is homogeneous.  
\xrema

 If $R$ is a ring,
an $M$-graded $R$-algebra is an $M$-graded ring $A$ endowed with a ring homomorphism $R \to A_0$. (In particular, $A$ is then an $R$-algebra.)
%an $R$-algebra $A$ such that $A$ is an $M$-graded ring and such that $A_m$ is an $R$-module for all $m \in M$.

If $I$ is a set and $\bbM$ is a commutative unital semiring we put $\bbM_I = \bigoplus_{i \in I} \bbM$. Then $\bbM_I$ is an $\bbM$-semimodule. In this setting we will denote by $(e_i : i \in I)$ the canonical basis of $\bbM_I$ (i.e.~$e_i$ is the $I$-tuple with $1$ in place $i$ and $0$ in other places).

If $A$ is a ring and $a \in A$ is an element, we set $a^\bbN := \{ a^n : n \in \bbN \} \subset A$. (Here $\bbN$ denotes the set of \emph{non negative} integers; this does not follow the conventions of~\cite{stacks-project}.) For any $a \in A$, by convention we have $a^0=1$.

If $X$ is a scheme, we will denote by $\calO_X$ its structure sheaf, and set $\calO(X)=\Gamma(X, \calO_X)$. We will also set $X_\aff = \Spec(\calO(X))$; then we have a canonical morphism of schemes $X \to X_\aff$.

If $A$ is a commutative ring and $I \subset A$ is an ideal, we will denote by $V(I) \subset \Spec(A)$ the closed subset defined by $I$. If $f \in A$ we will denote by $D(f) \subset \Spec(A)$ the open subscheme defined by $f$, i.e.~the complement of $V(A \cdot f)$.

If $N$ is any commutative monoid, there is a universal homomorphism from
$N$ to a group $N^{\mathrm{gp}}$ (cf. \cite[\S 1.2]{Og}). In case $N$ is cancellative\footnote{The terminology used in~\cite{Og} is \emph{integral}.} (i.e.~if the equality $x+y=x'+y$ implies that $x=x'$ for any $x,x',y \in N$) this map is injective. If $N$ is a submonoid of a given abelian group $M$, then $N$ is cancellative and the group $N^{\mathrm{gp}}$ naturally identifies with the subgroup $N -N := \{ n - n ' : n , n ' \in N \}$ of $M$.

Let $(N, +)$ be a commutative monoid. A submonoid $F$ is called a \emph{face} of $N$ if, for all $a,b \in N, $ whenever $a+b\in  F $ then both $a$
and $b$ belong to $F$.

%%%%%%%%%%%%%%%%%%%%%%%%%%%%%%%%%%%%%%%%
\section{Potions of graded rings}
\label{sec:potions}
%%%%%%%%%%%%%%%%%%%%%%%%%%%%%%%%%%%%%%%%

%--------------------------------------
\subsection{Graded rings and localizations} 
%--------------------------------------

Recall the notion of multiplicative subset in a commutative ring, see~\cite[\href{https://stacks.math.columbia.edu/tag/00CN}{Tag 00CN}]{stacks-project}. (In particular, a multiplicative subset always contains $1$, and it might contain $0$.) If a given commutative ring $A$ is $M$-graded for some abelian group $M$, we will say that a multiplicative subset is homogeneous if it consists of homogeneous elements. 

If $A$ is a commutative $M$-graded ring and $S \subset A$ is a homogeneous multiplicative subset, the localization $A_S$ of $A$ with respect to $S$ (which is denoted $S^{-1} A$ in many references, including~\cite[\href{https://stacks.math.columbia.edu/tag/00CM}{Tag 00CM}]{stacks-project}) is canonically $M$-graded, explicitly, for $m \in M$ we have
\[
( A_S)_m =  \left\{ \frac{a}{s} : \exists m', m'' \in M \text{ such that }  a \in A_{m'}, \, s \in (S \cap A_{m''}) \text{ and
$m'-m'' =m$ }\right\}.
%\subset A_S.
\]
Given a graded $A$-module $Q$, one can also consider the associated localization $Q_S$ (denoted $S^{-1}Q$ in~\cite[\href{https://stacks.math.columbia.edu/tag/07JZ}{Tag 07JZ}]{stacks-project}), which has a natural structure of a graded $A_S$-module.

Given a homogeneous multiplicative subset $S \subset A$, we will denote by $\underline{S}$ the homogeneous multiplicative subset consisting of homogeneous divisors of elements in $S$. Note that we have a canonical isomorphism of graded rings
\begin{equation}
\label{eqn:isom-localizations}
 A_{ \underline{S}} \cong A_{S}.
\end{equation}
If $Q$ is a graded $A$-module we also have a canonical identification of graded abelian groups 
\begin{equation}
\label{eqn:isom-localizations-modules}
Q_{\underline{S}} \cong Q_S
\end{equation}
compatible with the actions of $A_{ \underline{S}} = A_{S}$.

%--------------------------------------
\subsection{Relevant families} 
\label{ss:relevant-families}
%--------------------------------------

Consider an abelian group $M$, and a commutative $M$-graded ring $A$. By a homogeneous subset of $A$ we mean a subset consisting of homogeneous elements.

 \defi \label{def:deg}
 Let $S$ be a homogeneous subset of $A$. We denote by $\deg(S)$ the subset of $M$ defined as 
\[\deg(S):= \{ m \in M : \exists s \in S , s  \in A_m \}. \]
 \xdefi
Note that $\deg(\{0\})=M$.
More generally, if $S$ is a homogeneous multiplicative subset of $A$, then $\deg(S)$ is a submonoid of $M$, also denoted $M[S\rangle$.

\defi \label{MN}
Let $S$ be a homogeneous multiplicative subset of $A$. 
\begin{enumerate} 
\item We put $M[S]=M[S \rangle^{\mathrm{gp}}= M[S \rangle- M[S \rangle$, the subgroup of $M$ generated by $\deg(S)$.
\item If $M$ is finitely generated, we 
denote by $M[S \rangle_{\mathbb{R}_{\geq 0}}= M[S \rangle \otimes_{\N} \mathbb{R}_{\geq 0}$
the closed convex cone of $M \otimes_{\Z} \bbR$ generated by $\deg(S)$.
\end{enumerate}
\xdefi

The following definition gives a slight generalization of a terminology used in~\cite{BS03}.

\begin{dfn}
\phantomsection
%[(Relevant families)] 
\label{def:rele}
%Let $S$ be a multiplicative family of homogeneous elements of $A$ (that is $S$ is a submonoid of $(A,\times)$ and any element in $S$ is homogeneous). The subgroup of $M$ generated by the degrees of elements of  $S$ is denoted $M _{S}$ and is called the \textit{group of degrees generated} by $S$. Let $\underline{S}$ be the \textit{family of homogeneous divisors} of elements in $S$. 
% The family  $S$
\begin{enumerate}
 \item 
 A homogeneous multiplicative subset $S$ of $A$ is called \emph{$M$-relevant} (or just \emph{relevant} if $M$ is clear from the context) if for any $m$ in $M$ there exists $n \in \bbZ_{> 0}$ such that $nm$ belongs to $M[\underline{S}]$, i.e. if $M/(M[\underline{S}])$ is a torsion abelian group.
 \item
 A family $\{a_i : i \in I\}$ of homogeneous elements in $A$ is called $M$-relevant if the multiplicative subset generated by the $a_i$'s is relevant.
 \item
 \label{it:relevant-element}
 A homogeneous element $a \in A$ is called $M$-relevant if the family $\{a\}$ is $M$-relevant.
 \item The ideal of $A$ generated by all $M$-relevant elements of $A$ is denoted $A_{\dagger}$. 
\end{enumerate}
\end{dfn}

\begin{rmk}
\phantomsection
\label{rmk:relevant}
\begin{enumerate}
 \item 
 \label{it:existence-relevant}
 Note that the homogeneous multiplicative subset $0^\bbN = \{0,1\}$ is always relevant.
 \item
 \label{it:relevant-0}
 In case $M$ is finite, all homogeneous multiplicative subsets of $A$ are relevant, and all homogeneous elements are relevant.
 \end{enumerate}
\end{rmk}

\prop 
\label{prop:char-relevant}
Let $S$ be a homogeneous multiplicative subset of $A$.
\begin{enumerate}
\item 
\label{it:char-relevant}
We have $M[\underline{S}] = \{ m \in M :  (A_{\underline{S}})^\times \cap (A_{\underline{S}})_m \ne \varnothing \}$, where $(A_{\underline{S}})^\times \subset A_{\underline{S}}$ is the group of invertible elements. In particular $S$ is $M$-relevant if and only if 
\[
M/\{ m \in M :  (A_{\underline{S}})^\times \cap (A_{\underline{S}})_m \ne \varnothing \}
\]
is a torsion abelian group.
\item 
\label{it:fact-relevant}
Assume that $M$ is finitely generated. A homogeneous element $x$ in $A$ is relevant if and only if there exist $k \in \mathbb{N}_{>0} $ and a factorization $x^k=x_1 \cdots x_l$ into homogeneous factors $x_i$ of degree $d_i$ such that the degrees $d_i$ generate a subgroup of finite index in $M$. 
\end{enumerate}
\xprop 

\pf
\eqref{it:char-relevant}
The inclusion $M[\underline{S}] \subset \{ m \in M :  (A_{\underline{S}})^\times \cap (A_{\underline{S}})_m \ne \varnothing \}$ is obvious. Reciprocally let $m \in M$, and assume that there exists a fraction $\frac{a}{s}$ (with $a \in A$, $s \in \underline{S}$ homogeneous) which is invertible in $A_{\underline{S}}$ and of degree $m$. Since this element is invertible, there exist homogeneous elements $a' \in A$ and $s' \in \underline{S}$ such that $\frac{a}{s}\frac{a'}{s'}=1$. So there exists $s'' \in \underline{S}$ such that $s''aa'=ss's''$. This implies $a$ belongs to $\underline{S}$ and hence that $m=\deg(a)-\deg(s)$ belongs to $M[\underline{S}] $.

\eqref{it:fact-relevant}
Let $x$ be a relevant homogeneous element in $A$. Consider all homogeneous divisors of all positive powers of $x$, and the degrees of all these elements. By assumption, these degrees generate a subgroup $M'$ of $M$ of finite index. Since $M$ is finitely generated, so is $M'$ as a group. So we can find a finite number of homogeneous elements $x_1, \dots, x_m$ such that their degrees generate $M'$ and such that $x_i$ divides a power of $x$, say $x^{k_i}.$ Set $k = \sum_{k=1}^m k_i.$ Then the product $x_1\cdots x_m$ divides $x^k.$ So there exists $y$ such that $x^k=x_1 \cdots x_my.$ To finish the proof of the direct implication, it is enough to prove that $y$ can be chosen homogeneous, which follows from the following observation: if $x$ and $ z$ are homogeneous in a graded ring and $z$ divides $x$, then there exists a homogeneous $y$ in $A$ such that $x=zy.$ The reverse implication is immediate. 
\xpf

The following Proposition generalizes \cite[Lemma 2.7]{BS07}.\footnote{We insist that this statement only appears in~\cite{BS07}, and not in~\cite{BS03}.} Here, we denote by $\Rad(I)$ the radical of an ideal $I$. 
%(note that we have to use the 2007 v3 arXiv version of \cite{BS03} here, as Lemma 2.7 does not appear in the printed version of this work of Brenner-Schröer).  

\prop
\label{radical}
Let $M$ and $M'$ be two finitely generated abelian groups. Let $R$ be a commutative ring and let $A$ (resp. $A'$) be a commutative $M$-graded (resp. $M'$-graded) $R$-algebra. Then, considering $A \otimes_R A'$ with its natural structure of a $M \times M'$-graded $R$-algebra:
\begin{enumerate} 
\item
\label{it:radic-1}
for any relevant element $x \in A \otimes_R A'$, there exists a positive integer $k$, a finite set $J$, and for any $j \in J$ an $M$-relevant element $s_j$ in $A$ and an $M'$-relevant element $s'_j$ in $A'$ such that $x^k = \sum_{j \in J} s_j \otimes s'_j$;
\item 
\label{it:radic-2}
Denote by $\langle A_\dag \otimes_R A'_\dag \rangle$ the image of $A_\dag \otimes_R A'_\dag$ in $A \otimes_R A'$. Then $\langle A_\dag \otimes_R A'_\dag \rangle$ is an ideal in $A \otimes_R A'$, and
we have $\Rad \big( (A \otimes_R A')_{\dagger} \big) = \Rad \big( \langle A_{\dagger} \otimes_R A'_{\dagger} \rangle \big)$.
\end{enumerate}
\xprop 

\pf
\eqref{it:radic-1}
Let $x \in A \otimes_R A'$ be relevant. By Proposition~\ref{prop:char-relevant}\eqref{it:fact-relevant}, there exist $k \in \Z_{>0}$ and a factorization $x^k = x_1 \cdots x_l$ into homogeneous factors $x_i$ of degree $d_i \in M \times M'$, such that the degrees $d_i$ generate a subgroup of $M \times M'$ of finite index. For $i \in \{1 , \dots , l\}$, write $d_i = (m_i , m'_i)$ with $m_i \in M$ and $ m_i' \in M'$. Note that the degrees $m_i$ (resp. $m'_i$) generate a subgroup of $M$ (resp. $M'$) of finite index. For $i \in \{1 , \ldots , l\}$, write 
\[
x_i = \sum_{ j \in J_i } s_{ij} \otimes s'_{ij}
\]
for an index set $J_i$, and homogeneous elements $s_{ij} \in A $ and $s'_{ij} \in A'$ of respective degrees $m_i$ and $m'_i$. Then
 \[ 
 x^k= \prod_{i =1}^l \big( \sum_{ j_i \in J_i } s_{ij_i} \otimes s'_{ij_i}  \big) = \sum_{(j_i) \in \prod_{i=1}^l J_i} \big( \prod_{i=1}^l s_{ij_i}\big) \otimes \big( \prod_{i=1}^l s'_{ij_i}\big).
 \]
For every choice $(j_i) \in \prod_{i=1}^l J_i$, the element $\prod_{i=1}^l s_{ij_i}$, resp.~$\prod_{i=1}^l s'_{ij_i}$, is a product of elements of degrees $m_1 , \ldots , m_l$, resp.~$m'_1 , \ldots, m'_l$, and so it is $M$-relevant in $A$, resp.~$M'$-relevant in $A'$. This finishes the proof, taking $J =\prod_{i=1}^l J_i$ and the elements $\prod_{i=1}^l s_{ij_i}$ and $\prod_{i=1}^l s'_{ij_i}$.

\eqref{it:radic-2}
It is clear that $\langle A_{\dagger} \otimes_R A'_{\dagger} \rangle$ is an ideal contained in $(A \otimes_R A')_{\dagger}$, hence $\Rad \big( \langle A_{\dagger} \otimes_R A'_{\dagger} \rangle \big) \subset \Rad \big( (A \otimes_R A')_{\dagger} \big)$. On the other hand, by~\eqref{it:radic-1} we have $(A \otimes_R A')_{\dagger} \subset \Rad \big( \langle A_{\dagger} \otimes_R A'_{\dagger} \rangle \big)$, hence $\Rad \big( (A \otimes_R A')_{\dagger} \big) \subset \Rad \big( \langle A_{\dagger} \otimes_R A'_{\dagger} \rangle \big)$, which finishes the proof.
%The inclusion $\supset $ is immediate. The inclusion $\subset$ follows from \ref{radical}(\ref{it:radic-1}).
\xpf 

\begin{rmk}
\label{rmk:relevant-projection}
 Assume that $M'=\{0\}$, so that $A'$ is simply an $R$-algebra, which we will denote $R'$ for clarity, and we consider the $M$-graded algebra $A \otimes_R R'$. In this case we have $R'_\dag=R'$ (see Remark~\ref{rmk:relevant}\eqref{it:relevant-0}), hence the proposition says that $\mathrm{Rad}((A \otimes_R R')_\dag) = \mathrm{Rad}( \langle A_\dag \otimes_R R' \rangle)$, which is exactly~\cite[Lemma 2.7]{BS07}. This statement implies that the preimage of the closed subset $V(A_\dag) \subset \Spec(A)$ under the projection
 \[
  \Spec(A \otimes_R R') = \Spec(A) \times_{\Spec(R)} \Spec(R') \to \Spec(A)
 \]
is
\[
 V(\langle A_\dag \otimes_R R' \rangle) = V(\mathrm{Rad}( \langle A_\dag \otimes_R R' \rangle)) = V(\mathrm{Rad}((A \otimes_R R')_\dag)) = V((A \otimes_R R')_\dag).
\]
\end{rmk}

\begin{propo}
\label{prop:image-relevant}
 Let $M$ be a finitely generated abelian group, let $A,B$ be $M$-graded rings, and let $\Psi : A \to B$ be a homomorphism of $M$-graded rings.
 
 \begin{enumerate}
  \item 
  \label{it:image-relevant}
  We have $\Psi(A_\dag) \subset B_\dag$.
  \item 
  \label{it:image-relevant-surj}
  If $\Psi$ is surjective, then $\mathrm{Rad}(\Psi(A_\dag)) = \Rad(B_\dag)$.
 \end{enumerate}
\end{propo}

\begin{proof}
 It is easily seen that the image under $\Psi$ of any relevant element of $A$ is a relevant element of $B$, which proves~\eqref{it:image-relevant}.
 
 Now, assume that $\Psi$ is surjective. Then $\Psi(A_\dag)$ is an ideal, and by~\eqref{it:image-relevant} we have $\mathrm{Rad}(\Psi(A_\dag)) \subset \Rad(B_\dag)$. Let $f \in B$ be a relevant element. Then there exists $N \in \Z_{\geq 1}$, elements $m_1, \dots, m_r \in M$ which generate a subgroup of finite index, and for any $i$ an element $f_i \in B_{m_i}$, such that 
% and homogeneous elements $f_1, \dots, f_r \in B$ of respective degrees $m_1, \dots, m_r$, and a homogeneous element $g \in B$ such that 
 $f^N = f_1 \cdots f_r$.
 % and the subgroup of $M$ generated by $m_1, \dots, m_r$ has finite index. Choose homogeneous elements $\tilde{f}_1, \dots, \tilde{f}_r, \tilde{g} \in A$ whose respective images are $f_1, \dots, f_r, g$. 
Choose for any $i$ a preimage $\tilde{f}_i \in A_{m_i}$ for $f_i$. Then $\tilde{f}=\tilde{f}_1 \cdots \tilde{f}_r$ is relevant and satisfies $\Psi(\tilde{f})=f^N$, hence $f \in \mathrm{Rad}(\Psi(A_\dag))$. This implies that $B_\dag \subset \mathrm{Rad}(\Psi(A_\dag))$, hence $\mathrm{Rad}(B_\dag) \subset \mathrm{Rad}(\Psi(A_\dag))$, which finishes the proof of~\eqref{it:image-relevant-surj}.
\end{proof}

\begin{rmk}
\label{rmk:preimage-Adag-surj}
 In the setting of Proposition~\ref{prop:image-relevant}\eqref{it:image-relevant-surj}, the preimage of the closed subset $V(A_\dag)$ under the closed immersion $\Spec(B) \subset \Spec(A)$ is
 \[
 V(\Psi(A_\dag)) = V(\mathrm{Rad}(\Psi(A_\dag))) = V(\Rad(B_\dag))=V(B_\dag).
 \]
\end{rmk}

%FIN RELECTURE

%--------------------------------------
\subsection{Potions} 
%--------------------------------------

Let $M$ be an abelian group, and let $A$ be a commutative $M$-graded ring.

\defi[(Potions)]
\label{potionsdef} 
Let $S$ be a homogeneous multiplicative subset of $A$. 
\begin{enumerate}
\item
The degree-$0$ part $(A_S)_0$ of the localization $A_{S}$ is denoted $A_{(S)}$ and is called the \emph{potion} of $A$ with respect to $S$. 
\item
If $Q$ is a graded $A$-module, we will also denote by $Q_{(S)}$ the degree-$0$ part $(Q_S)_0$ of $Q_S$; it admits a canonical structure of an $A_{(S)}$-module.
\end{enumerate}
\xdefi

Note that in the setting of Definition~\ref{potionsdef}, by~\eqref{eqn:isom-localizations}--\eqref{eqn:isom-localizations-modules} we have canonical identifications
$A_{(\underline{S})} \cong A_{(S)}$ and $Q_{(\underline{S})} \cong Q_{(S)}$.
%If $a \in A$ is homogeneous, we set $A_{(a)} := A_{(a^\bbN)}$. 
If  $\{a_i : i \in I \}$ is a family of homogeneous elements of $A$, we will denote by
$A_{(\{a_i : i \in I \})}$
the potion associated with the multiplicative subset of $A$ generated by $\{a_i : i \in I \}$; in case $\#I=1$ we will write $A_{(a)}$ for $A_{(\{a\})}$.
% Note that $A_{( \underline{S})}=A_{( S)}$.

%If $a,b \in A$ are homogeneous and $\deg (a)= \deg (b)$, then we write $a \smile b$.
If $S$ and $T$ are multiplicative subsets of $A$, we will denote by $ST$ the multiplicative subset of $A$ generated by $S \cup T$, i.e. $ST = \{st : s \in S, \, t \in T \}$. Of course, $ST$ is homogeneous if $S$ and $T$ are.
The following result generalizes \cite[Lemma 2.2.2]{Gr61} and \cite[Lemma 3.7]{Za19}, and is the key result that makes the Proj construction work.

\prop[(Magic of potions)] 
\label{prop:magical} 
Let $S$ and $T$ be homogeneous multiplicative subsets of $A$. 
%We denote by $ST$ the multiplicative family of homogeneous elements of $A$ generated by $S \cup T$, i.e. $ST = \{st |s \in S , t \in T \}$. 
\begin{enumerate}
\item 
\label{it:magic-1}
We have a canonical homomorphism of potion rings 
%\[\chi : 
$A_{(S )} \to  A_{(ST )}$.
%\]
\item 
\label{it:magic-2}
Assume that $S$ is relevant. Fix a subset $T' \subset T$ which generates $T$ as a submonoid of $(A,\times)$ and, for any $t $ in $T'$, fix $n_t \in \mathbb{N}_{>0}$ and $s_t, s_t' \in \underline{S}$ such that $\deg(t^{n_t}) = \deg(s_t)-\deg(s_t')$.
%$t^{n_t} \smile s_t$. 
Then $\frac{t^{n_t} s_t'} { s_t}$ belongs to the potion $A_{(\underline{S} )}=A_{(S)}$. 
 Moreover we have a canonical isomorphism of $A_{(S )}$-algebras between $ A_{(ST )}$ and the localization of $A_{(S )}$ with respect to the multiplicative subset of $A_{(S)}$ generated by the elements $\{\frac{t^{n_t} s_t'} { s_t} : t \in T' \}$.
\item 
\label{it:magic-3}
Assume that $S$ is relevant and that $T$ is finitely generated as submonoid of $(A,\times)$. The morphism of schemes 
%dual to $\chi$ 
\[
\mathrm{Spec} (A_{(ST )}) \longto \mathrm{Spec} (A_{(S)})
\] 
induced by the ring homomorphism in~\eqref{it:magic-1} is an open immersion of schemes.
\item 
\label{it:magic-4}
Let $f_1, \ldots, f_n \in A$ be nonzero relevant homogeneous elements of the same degree.
Then we have a canonical open immersion
\[
\Spec(A_{ (f_1 + \cdots + f_n) }) \longto \Spec (A_{(f_1)}) \cup \cdots \cup \Spec (A_{(f_n)})
\]
where the right-hand side is defined as the glueing (in the sense of~\cite[\href{https://stacks.math.columbia.edu/tag/01JA}{Tag 01JA}]{stacks-project}) of the affine schemes $\Spec (A_{(f_i)})$ along the open subschemes $\Spec (A_{(f_i \cdot f_j)}) \subset \Spec (A_{(f_i)})$ (see~\eqref{it:magic-3}).
%pieces using \eqref{it:magic-3} (see also Construction~\ref{defproj} below for details).
\end{enumerate}
\xprop

\begin{proof}
\eqref{it:magic-1}
The desired morphism 
%$\chi$ 
is obtained as the degree-zero part of the canonical homomorphism of graded rings $A_{S} \to A_{ST}$.

\eqref{it:magic-2}
Without loss of generality, we can (and will) assume that $S= \underline{S}$. Consider the localization $\big( A_{(S)}  \big)_{\{\frac{t^{n_t} s_t'} { s_t} : t \in T' \}}$ of $A_{(S )}$ with respect to the homogeneous multiplicative subset generated by $\{\frac{t^{n_t} s_t'} { s_t} : t \in T' \}$. Since the image of $\frac{t^{n_t} s_t'} { s_t}$ in $A_{(ST)}$ is invertible for any $t \in T' $, the universal property of localizations gives us a canonical homomorphism of $A_{(S)}$-algebras
\[
\phi : \big( A_{(S)}  \big)_{\{\frac{t^{n_t} s_t'} { s_t} : t \in T' \}} \to A_{(ST)}
\]
sending $\frac{\frac{a}{s}}{\prod_{t} (\frac{t ^{n_t} s_t'}{ s_t}) ^{k_t} } $ to  $\frac{a}{s} \cdot \prod_{t} (\frac{ s_t}{t ^{n_t} s_t'}) ^{k_t} $. (Here and below we consider products indexed by $T'$; we tacitly assume that only finitely many of the exponents are nonzero, i.e.~that $(k_t : t \in T')$ belongs to $\N_{T'}$.)

Let us prove that $\phi$ is an isomorphism.  If an element $\frac{\frac{a}{s}}{\prod_{t} (\frac{t ^{n_t} s_t'}{ s_t}) ^{k_t} } $ belongs to $\ker (\phi )$, then we have $\frac{a (\prod_t s_t^{k_t})}{s (\prod_t  {t}^{n_t k_t} (s_t')^{k_t})}=0$ in  $A_{(ST )}$. So there exist $\mathfrak{s} \in S$ and $\mathfrak{t} \in T $ such that $a \cdot (\prod_t s_t^{k_t}) \cdot \mathfrak{s} \mathfrak{t} =0$ in $A$. Now the element $\frac{\frac{a}{s}}{\prod_{t} (\frac{t ^{n_t} s_t'}{ s_t}) ^{k_t} }  $ equals zero in $\big( A_{(S)}  \big)_{\{\frac{t^{n_t} s_t'} { s_t} : t \in T' \}} $ if and only if there exists $(d_t ) \in \mathbb{N}_{T'}$ such that $\frac{a}{s} \prod_t (\frac{t ^{n_t} s_t'}{ s_t}) ^{d_t}$ equals zero in $ A_{(S )} $; that is if and only if there exists $\mathsf{s} \in S$ such that $ a \cdot \mathsf{s} \cdot \prod_t {t ^{n_t d_t} }=0$ in $A$. This shows that $\phi$ is injective. 
%We now prove that $\phi$ is surjective. So 
Let now $\frac{a}{\mathfrak{s} \mathfrak{t}}$ be an element in $A_{(ST)}$, with $\mathfrak{s} \in S$, $\mathfrak{t} \in T$, and $a \in A$ homogeneous of degree $\deg(\mathfrak{st})$. Write $\mathfrak{t} = \prod_{t} t^{d_t}$ for some $(d_t) \in \mathbb{N}_{T'}$.
%Recall that we have fixed $n_{\mathfrak{t}} \in \bbN _{>0}$ and $s_{\mathfrak{t}}\in S$ such that $\deg(\mathfrak{t} ^{n_{\mathfrak{t}}}) = \deg(s_{\mathfrak{t}})$.
 The equality 
\[
 \frac{a}{\mathfrak{st}} = \frac{a \cdot \prod_{t} t^{d_t (n_t-1)} (s_t')^{d_t}}{\mathfrak{s} \prod_{t} s_{t}^{d_t}} \cdot \prod_t \left( \frac{s_{t}}{t^{n_t} s_t'} \right)^{d_t}
 \]
holds in $A_{(ST)}$, and implies that $\frac{a}{\mathfrak{st}}$ belongs to the image of $\phi$.
%= \phi \Big( \frac{ \frac{a {\mathfrak{t}}^{n_{\mathfrak{t}}-1}}{\mathfrak{s} s_{\mathfrak{t}}}}{\frac{\mathfrak{t}^{n_{\mathfrak{t}}}}{s_{\mathfrak{t}}}} \Big)$. 
This proves that $\phi$ is surjective, and concludes the proof.

\eqref{it:magic-3}
In case $T$ contains $0$, the claim is clear. Otherwise, it follows immediately from~\eqref{it:magic-2}, since finite localizations of rings induce open immersions.

\eqref{it:magic-4}
For $i \in \{1 , \ldots, n \}$, put $f_i'= \frac{f_i}{f_1 + \cdots + f_n} \in A_{(f_1 + \cdots + f_n)}$ and consider the open subscheme
\[
D(f_i')= \Spec \left( \big(A_{(f_1 + \cdots + f_n)}\big)_{f_i'} \right) \subset \Spec (A_{(f_1 + \cdots + f_n)} ),
\]
where $(A_{(f_1 + \cdots + f_n)})_{f_i'}$ is the localization of $A_{(f_1 + \cdots + f_n)}$ with respect to the multiplicative subset generated by $f_i'$.
We have $f_1'+ \cdots + f_n' =1$ in
 $A_{(f_1 + \cdots + f_n)}$, so by (3) in~\cite[\href{https://stacks.math.columbia.edu/tag/01HS}{Tag 01HS}]{stacks-project} we obtain that
 \[
 D(f_1') \cup \cdots \cup D(f_n') = \Spec (A_{(f_1 + \cdots + f_n)} ).
 \]
  Now by~\eqref{it:magic-2} we have $(A_{(f_1 + \cdots + f_n)}\big)_{f_i'} = A_{\big( (f_1 + \cdots + f_n) f_i \big)} $ and therefore $D(f_i')$ identifies with an open subscheme in $\Spec (A_{(f_i)}) $, which finishes the proof.
\end{proof}

Proposition~\ref{prop:magical}\eqref{it:magic-4} leads us to introduce the following terminology. 

\defi
\label{defquasirel}
A nonzero homogeneous element $a \in A$ is called \emph{quasi-$M$-relevant} (or just quasi-relevant) if it is a sum of $M$-relevant elements of the same degree. 
\xdefi

\rema
Note that if $M \neq \Z$ there might exist quasi-relevant elements which are not relevant. For example, let $A= \bbZ [X,Y,Z,T] $ be $\bbZ^2$-graded with $\deg(X)=\deg(Y)= e_1$ and $\deg(Z)=\deg(T)= e_2$. Then $XZ $ and $YT$ are relevant. So $XZ+YT$ is quasi-relevant, but it is not relevant. 
\xrema 

%-----------------------------------
\subsection{Examples}
\label{ss:examples}
%-----------------------------------

\exam[(Multi-centered dilatations)] \label{exmulti}
Let $A$ be a commutative ring and let $\{[M_i , a_i ] : i \in I \}$ be a multi-center in $A$ in the sense of \cite{Ma24}, i.e. $M_i$ is an ideal of $A$ and $a_i $ is an element of $A$ for any $i \in I $. For $i \in I $, let  $L_i $ be the ideal $ M_i +(a_i)$ of $A$.
Let $\Bl_{\{L_i : i \in I\}} A= \bigoplus_{ \nu \in \bbN_I } L^{\nu}$ be the multi-Rees $A$-algebra associated with $A$ and $\{L_i : i \in I\}$. (Here, if $\nu=(\nu_i : i \in I)$, then $L^\nu$ is the product of the $(L_i)^{\nu_i}$'s.) The ring $\Bl_{\{L_i : i \in I\}} A$ is naturally $\bbZ_I $-graded, and each $a_i$ defines an element in its degree-$e_i$ part, also denoted $a_i$. Then the family $\{a_i : i \in I\}$ is relevant, 
and~\cite[Fact 2.35]{Ma24} shows that the potion ring $(\Bl_{\{L_i : i \in I \}} A)_{(\{a_i : i \in I\})}$ identifies with the dilatation $A [\{ \frac{M_i}{a_i} : i \in I \}]$.
\xexam 

\exam 
\label{examMo}
Any commutative ring $A$ can be considered as graded by the abelian group $M=0$. In this case any element $a$ in $A$ is relevant and $A_{(a)}= A_a$ is the localization of $A$ with respect to $a$ as in~\cite[\href{https://stacks.math.columbia.edu/tag/02C5}{Tag 02C5}]{stacks-project}.
\xexam 

\section{Brenner--Schr\"oer Proj} \label{sec:projdef}
%%%%%%%%%%%%%%%%%%%%%%%%%%%%%%%%%%%%%%%%

%--------------------------------------
\subsection{Proj scheme of a graded ring}
\label{subsectprojring}
%--------------------------------------

%We proceed with the notation of Section~\ref{sec:potions}, but we now assume that 
From now on we assume that $M$ is a \emph{finitely generated} abelian group, and fix a commutative $M$-graded ring $A$. 
Consider the affine scheme $\Spec(A_0)$, and the diagonalizable group scheme
\[
\mathrm{D}_{\Spec(A_0)}(M) = \Spec(A_0[M])
\]
over $\Spec(A_0)$
associated with $M$. The $M$-grading on $A$ defines an action of $\mathrm{D}_{\Spec(A_0)}(M)$ on $\Spec(A)$; we will denote by
\[
a,p : \mathrm{D}_{\Spec(A_0)}(M) \times_{\Spec(A_0)} \Spec(A) \to \Spec(A)
\]
the action and projection morphisms, respectively.

In~\cite[Definition~2.2]{BS03}, Brenner and Schr\"oer define the Proj scheme\footnote{Brenner--Schr\"oer use the terminology ``homogeneous spectrum'' of $A$. We prefer to use the term ``Proj scheme'' to emphasize the parallel with Grothendieck's construction.} associated with $A$ as a certain open subset in the ringed space $\mathrm{Quot}(A)$ obtained as the cokernel (in the sense of~\cite[Exp.~V, \S 1.b]{sga3}) of the maps $a$ and $p$. In Construction~\ref{defproj} below, we explain an equivalent description of this scheme, obtained by gluing spectra of certain potion rings. 
(See~\S\ref{ss:first-prop}
%the comments following the construction 
for a justification that our description is indeed equivalent to that in~\cite{BS03}.)

We will denote by $\mathcal{F}_A$ the set of all relevant homogeneous multiplicative subsets of $A$ which are finitely generated as submonoids of $(A,\times)$.

\begin{const}[(Proj as glueing potions)]
\label{defproj}
Let
$\mathcal{F} \subset \mathcal{F}_A$ be a subset.
% of relevant homogeneous multiplicative subsets of $A$ such that any element in $\mathcal{F}$ is finitely generated as a submonoid of $(A, \times)$. 
 For each $S \in \mathcal{F}$, let $D_{\dagger}(S)$ be the spectrum of the potion $A_{(S)}$. 
 By Proposition~\ref{prop:magical}\eqref{it:magic-3}, if $S,T \in \mathcal{F}$, the affine scheme $D_{\dagger}(ST)$ identifies canonically with an open subscheme of $D_{\dagger}(S)$. For each $S,T \in \mathcal{F}$, we have equalities
 \[
 D_{\dagger}({S S}) = D_{\dagger}(S) \quad \text{and} \quad D_{\dagger}({ST})=D_{\dagger}({TS}).
 \]
 Moreover, for each triple $S,T,U \in \mathcal{F}$, we have
 \[
 D_{\dagger}({ST} )\cap D_{\dagger}({SU}) = D_{\dagger}({TS}) \cap D_{\dagger}({TU})
 \]
 (intersections in $D_{\dagger}(S)$ and $D_{\dagger}(T)$ respectively; equalities in $D_{\dagger}(S)$, $D_{\dagger}(T)$ or $D_{\dagger}(U)$). Indeed, using Proposition~\ref{prop:magical}\eqref{it:magic-2}, these intersections identify with $\mathrm{Spec} (A_{(STU)} )$.
Now, by glueing~\cite[\href{https://stacks.math.columbia.edu/tag/01JA}{Tag 01JA}]{stacks-project}, from these data we obtain a scheme $\mathrm{Proj}^M_{\mathcal{F} } (A)$ and, for each $S \in \mathcal{F}$, an open immersion $\varphi_S : D_{\dagger}(S) \to \mathrm{Proj}^M_{\mathcal{F}} (A)$, such that
\[
\mathrm{Proj}^M_{\mathcal{F} } (A) = \bigcup_{S \in \mathcal{F} } \varphi_S( D_{\dagger}(S))
\]
and that for $S,T \in \mathcal{F}$ we have
\[
\varphi_S ( D_{\dagger}({ST}) )=\varphi_T ( D_{\dagger}({TS}) ) = \varphi_S (D_{\dagger}(S) ) \cap \varphi_{T} (D_{\dagger}(T)).
\] 
In practice, we will often identify $D_{\dagger}(S) $ and $\varphi_S (D_{\dagger}(S))$. 
%We now put $X=:\mathrm{Proj}_{\mathcal{F} } (A)$.

In the case when $\mathcal{F} = \mathcal{F}_A$,
% is the set of all relevant homogeneous multiplicative subsets of $A$ which are finitely generated as submonoids of $(A,\times)$, 
 the scheme $\mathrm{Proj}^M_{\mathcal{F}_A} (A)$ will be denoted $\mathrm{Proj}^M(A)$, or just $\mathrm{Proj}( A)$ when $M$ is clear from the context.
 \end{const}

\begin{rmk}
\phantomsection
\label{rmk:Proj-opens}
\begin{enumerate}
 \item 
 As explained in Remark~\ref{rmk:relevant}, $\{0\} \in \calF_A$. However, the localization $A_{\{0\}}$ is the zero ring, so that $D_\dag(\{0\})$ is empty.
\item
\label{rmk:open-contain-0-empty}
More generally,
if $0 \in S$ then $A_S=\{0\}$, so that $D_\dag(S)=\varnothing$. The scheme $\Proj^M(A)$ is therefore covered by the open subschemes $D_\dag(S)$ where $S$ does not contain $0$.
\item
\label{it:Proj-opens-3}
Let $S,T \in \calF_A$, and assume that $S \subset T \subset \underline{S}$. Then $T=ST$, so that we have an embedding $D_\dag(T) \subset D_\dag(S)$. Since $A_S=A_T$, this embedding is an equality.
\item Assume that $\{1\}$ is relevant. Then $D_\dag (\{1\})= \Spec (A_0)$ by definition. Moreover, for any relevant family $S$, we have $S=S\{1\}$ and $D_\dag (S ) = \left( D_\dag (S )\cap D_\dag (\{1\}) \right) \subset D_\dag (\{1\})$. So $\Proj^M(A)= \Spec (A_0)$.
\end{enumerate}
 \end{rmk}
 
 \begin{ex}
\label{ex:Proj-Spec}
In view of Example~\ref{examMo}, when $M=\{0\}$, for any ring $A$ we have $\Proj^{\{0\}}(A)=\Spec(A)$.
\end{ex}

 \begin{rmk} \label{rema:topol}
 If $A$ is any commutative ring, let $\mathrm{Prim}(A)$ be the set of prime ideals of $A$, without topology. For any subset $T \subset A$, let $D(T) \subset \mathrm{Prim}(A)$ be the set $\{ \mathfrak{p} \in \mathrm{Prim}(A) : \forall t \in T, \, t \not \in \mathfrak{p} \}.$ Let $\mathrm{Alex}(A)$ be the set $\mathrm{Prim}(A)$, endowed with the topology generated by all the $D(T)$ for $T \subset A$. Then $\mathrm{Alex}(A)$ is an Alexandroff topological space and in fact the Alexandroffication of $\Spec(A)$. We call the topology of $\mathrm{Alex}(A)$ the Alexandroff--Zariski topology. In general, this is not the discrete topology, nor the Zariski topology (e.g. take $A= \bbZ$). 
 
 Recall that the association $D(f) \mapsto A_f$ (for $f \in A$) extends to a sheaf of rings on $\Spec(A)$, namely the structure sheaf $\calO _{\Spec(A)}$ \cite[\href{https://stacks.math.columbia.edu/tag/01HU}{Tag 01HU}]{stacks-project}. It is a natural and elementary question to ask whether the association $D(T) \mapsto A_T$ (for $T \subset A$) extends to a sheaf of rings on $\mathrm{Alex}(A)$, refining $\calO _{\Spec(A)}$. In other words, considering the continuous Alexandroffication map $ f: \mathrm{Alex}(A) \to \Spec(A)$, does the pullback presheaf $f_p \calO_{\Spec(A)} $ coincide with the pullback sheaf $ f^{-1} \calO_{\Spec(A)}$? (See e.g.~\cite[\href{https://stacks.math.columbia.edu/tag/008C}{Tag 008C}]{stacks-project} for the notation $f_p$ and $f^{-1}$.) Even if it is easy to see that the answer is Yes in many cases (e.g. if $A$ is integral), the answer is negative in general. For example take $A= (\bbZ /2\bbZ)^\bbN$ (Boolean ring). In this case $\mathrm{Alex}(A)$ is the discrete topology, so that while the global sections of $f_p \calO_{\Spec(A)}$ give $A$, the global sections of $f^{-1} \calO_{\Spec(A)} $ are products of $\bbZ/2\bbZ$ over $\mathrm{Prim}(A)$. It is well-known that $\mathrm{Prim}(A)$ is in bijection with the Stone--Cech compactification of $\bbN$, which has ordinality strictly bigger than $\bbN$; as a consequence, $f_p \calO _{\Spec(A)} \ne f^{-1} \calO_{\Spec(A)}$. 
 
 Therefore, one cannot refine the structure sheaf of affine schemes (and the theory of schemes) in this way, using arbitrary localizations. So it is not possible to create a theory, similar to the theory of schemes in the framework of locally ringed space, allowing glueing along open subspaces in the Alexandroff--Zariski topology, which justifies Remark~\ref{rema:topolintro}. Note that Exercise~11.1 of~\cite{Sch19} says that one \emph{can} produce a theory inspired by the theory of schemes using arbitrary localizations, but it is not a theory based on locally ringed spaces, and will not behave formally as scheme theory regarding glueings. Note that the theory of Analytic Spaces of \cite{Sch19} uses sorts of quasi-compact coverings, so that one cannot define the Proj of an arbirary graded ring in this theory either. In the above geometric paradigms, it thus seems that the hypothesis that $A$ is graded by a \emph{finitely generated} abelian group cannot be removed without serious difficulties to define $\Proj$ schemes.
 \end{rmk}

If $\{a_i : i \in I\}$ is a finite relevant family of elements of $A$, we will write $D_\dag(\{a_i : i \in I\})$ for the open subscheme defined by the multiplicative subset of $A$ generated by the $a_i$'s. In case $\#I=1$ we will write $D_\dag(a)$ for $D_\dag(\{a\})$. 
 
 %--------------------------------------
\subsection{First properties}
\label{ss:first-prop}
%--------------------------------------
 
 By construction, $\Proj^M(A)$ is covered by the affine open subschemes $D_\dag(S)$, which have the property that the intersection of any two such subschemes is affine. This implies that the diagonal morphism of $\Proj^M(A)$ is affine (as in~\cite[Proposition~3.1]{BS03}), hence (since affine morphisms are quasi-compact) we get the following statement.

\begin{lem} 
\label{lem:Proj-qsep}
The scheme $\Proj^M(A)$ is quasi-separated.
\end{lem}

% \begin{proof}
% This follows from~\cite[\href{https://stacks.math.columbia.edu/tag/01KO}{Tag 01KO}]{stacks-project}, since by construction the intersections in $\Proj^M(A)$ of affine open subschemes of the form $D_\dag(S)$ are affine.
% \end{proof}

For any $S \in \mathcal{F}_A$ we have a canonical homomorphism $A_0 \to A_{(S)}$; these homomorphisms glue to a morphism of schemes
\begin{equation}
\label{eqn:morph-Spec-A0}
\Proj^M(A) \longto \Spec(A_0).
\end{equation}

%Using the notation introduced in Definition~\ref{def:rele}\eqref{it:relevant-element}, in case $S=a^{\N}$ where $a \in A$ is $M$-relevant, we will write $D_\dag(a)$ for $D_\dag(a^{\N})$. 
For the next statement, note that if $f \in A$ is relevant, for any homogeneous $g \in A$ the product $fg$ is relevant. 
%In the next lemma we describe the associated open subscheme in two special cases.

\begin{lem}
\label{lem:prod-relevant}
Let $f \in A$ be relevant, and let $g \in A$ be homogeneous.
\begin{enumerate}
\item
\label{it:prod-relevant-0}
We have $D_\dag(fg) \subset D_\dag(f)$.
\item
\label{it:prod-relevant-1}
If $g$ is also relevant, then we have
$D_\dag(fg) =
D_\dag(f) \cap D_\dag(g)$.
\item
\label{it:prod-relevant-2}
If $g \in A_0$, then $D_\dag(fg)$ is the preimage of $D(g) \subset \Spec(A_0)$ under the composition
\[
D_\dag(f) \hookrightarrow \Proj^M(A) \xrightarrow{\eqref{eqn:morph-Spec-A0}} \Spec(A_0).
\]
\end{enumerate}
\end{lem}

\begin{proof}
\eqref{it:prod-relevant-0}
We have $(fg)^\N \subset f^\N \cdot (fg)^\N \subset \underline{(fg)^\N}$, hence by Remark~\ref{rmk:Proj-opens}\eqref{it:Proj-opens-3} we have $D_\dag(f) \cap D_\dag(fg) = D_\dag(f^\N \cdot (fg)^\N) = D_\dag(fg)$, so that $D_\dag(fg) \subset D_\dag(f)$.

\eqref{it:prod-relevant-1}
%By definition we have $D_\dag(f) \cap D_\dag(g) = D_\dag(f^\N \cdot g^\N)$. Then the claim 
This follows similarly from Remark~\ref{rmk:Proj-opens}\eqref{it:Proj-opens-3} since $(fg)^\N \subset f^\N \cdot g^\N \subset \underline{(fg)^\N}$.

\eqref{it:prod-relevant-2}
By~\eqref{it:prod-relevant-0} we have
%Remark~\ref{rmk:Proj-opens}\eqref{it:Proj-opens-3} we have $D_\dag(fg) = D_\dag(f^\N \cdot (fg)^\N)$, hence 
$D_\dag(fg) \subset D_\dag(f)$. Now it is clear that $A_{(fg)}$ is the localization of $A_{(f)}$ with respect to $g$, which gives the claim.
\end{proof}

%By Fact~\ref{powerpotionequal}, 
In particular, Lemma~\ref{lem:prod-relevant}\eqref{it:prod-relevant-1} shows that if $a \in A$ is relevant,
for any $k \in \Z_{\geq 1}$ we have
\begin{equation}
\label{eqn:open-power}
D_\dag(a)=D_\dag(a^k).
\end{equation}

Note that if $S \in \mathcal{F}_A$, and if $f_1 , \ldots, f_n$ are multiplicative generators of $S$, then $f:=f_1 \cdots f_n$ is a relevant element of $A$, and by Remark~\ref{rmk:Proj-opens}\eqref{it:Proj-opens-3} we have $D_\dag(S)=D_\dag(f)$.
% and $A_{(f)}=A_{(S)}$. 
%Hence the open immersion
%\[
% D_\dag(S) \to D_\dag(f)
%\]
%(where we note that $S=S \cdot f^\N$) is an isomorphism.
It follows that
\begin{equation}
\label{eqn:singletonproj}
 \Proj^M(A) = \bigcup_{f \in A \text{ relevant}} D_\dag(f).
\end{equation}
In fact, in view of Remark~\ref{rmk:Proj-opens}\eqref{rmk:open-contain-0-empty} we can restrict this union to \emph{non nilpotent} relevant elements.

%in Construction~\ref{defproj} it is enough to glue potions associated with relevant homogeneous elements in $A$ to obtain the scheme $\Proj^M (A)$. 
Comparing Lemma~\ref{lem:prod-relevant}\eqref{it:prod-relevant-1} with the proof of~\cite[Proposition~3.1]{BS03} we see that the intersections of these open subschemes are as in the scheme constructed in~\cite{BS03},
which justifies that our scheme $\Proj^M(A)$ coincides with that defined in~\cite[Definition~2.2]{BS03}. 
%Note that if $f,g \in A$ are relevant then so is $fg$, and by Remark~\ref{rmk:Proj-opens}\eqref{it:Proj-opens-3} we have
%%similar considerations as above show that
%\[
% D_\dag(f) \cap D_\dag(g) = D_\dag(f^\N \cdot g^\N) = D_\dag((fg)^\N).
%\]
In particular the scheme $\Proj^M(A)$ can be constructed by gluing open subschemes associated with relevant elements along their intersections, which 
%by Lemma~\ref{lem:prod-relevant}\eqref{it:prod-relevant-1} 
are again open subschemes associated with relevant elements. We however feel that the possibility of defining open subschemes associated with relevant \emph{families} adds some flexibility which might be useful.

\begin{rmk}
\phantomsection
\label{rmk:def-Proj}
\begin{enumerate}
%  \item 
%  \label{it:singletonproj}
\item
In case $M=\Z$, any nonzero homogeneous element of nonzero degree is relevant. If we furthermore have $A_n = 0$ for any $n \in \Z_{<0}$ one recovers the usual Proj scheme associated with a (nonnegatively) graded ring as in~\cite[\href{https://stacks.math.columbia.edu/tag/01M3}{Tag 01M3}]{stacks-project}; we will refer to this setting as the $\N$-graded setting. 
\item
In the $\N$-graded setting, the scheme $\Proj^{\Z}(A)$ is always separated; see~\cite[\href{https://stacks.math.columbia.edu/tag/01MC}{Tag 01MC}]{stacks-project}. For general $M$ and $A$ this is not true, even when $M=\bbZ$. For examples and separatedness criteria, see~\cite[\S 3]{BS03} and Proposition \ref{prop:separated-family} below.  
\end{enumerate}
\end{rmk}

The following result is an immediate corollary of \cite[Proposition 3.3]{BS03}.
Given a finite family $T=\{a_i : i \in I\}$ of homogeneous elements of $A$, generating a multiplicative subset $S$, we will denote by $ C_T \subset M \otimes_{\Z} \mathbb{R}$ the closed convex cone $M[\underline{S} \rangle_{\mathbb{R}_{\geq 0}}$ (cf.~Definition~\ref{MN}). Note that $T$ is relevant if and only if $C_T$ has nonempty interior.
%nonempty interior.

\prop
\label{prop:separated-family}
Let $\calF$ be a collection of relevant finite homogeneous families such that $C_{T} \cap C_{T'} \subset M \otimes_\Z \mathbb{R}$ has nonempty interior for all $T,T' \in \calF$. Then the scheme $\Proj^M_{\calF} (A)$ is separated.
\xprop

\pf
This follows from~\cite[Proposition 3.3]{BS03} and the comment preceding~\eqref{eqn:singletonproj}.
%Definitions \ref{def:rele} \& \ref{defproj}.
\xpf 

The natural morphism of ringed spaces $\Spec(A) \to \Quot(A)$ restricts to a canonical morphism
\begin{equation}
\label{eqn:geometric-quotient}
 \Spec(A) \smallsetminus V(A_\dag) \longto \Proj^M(A).
\end{equation}
In the description given in~\eqref{eqn:singletonproj}, this morphism can be understood as follows. For any $f \in A$ relevant, we have a canonical embedding $A_{(f)} \subset A_f$, which provides a morphism $\Spec(A_f) \to \Spec(A_{(f)})$. The latter morphisms glue to the desired morphism from
\[
 \Spec(A) \smallsetminus V(A_\dag) = \bigcup_{f \in A \text{ relevant}} \Spec(A_f)
\]
to
\[
 \Proj^M(A) = \bigcup_{f \in A \text{ relevant}} \Spec(A_{(f)}).
\]
This morphism is affine since the preimage of $\Spec(A_{(f)})$ is $\Spec(A_{f})$, for any relevant $f$.

The morphism~\eqref{eqn:geometric-quotient} is discussed further in~\cite[Comments after Definition~2.2]{BS03}. The closed subscheme $V(A_\dagger) \subset \Spec(A)$ associated with $A_\dagger$ is stable under the action of $\mathrm{D}_{\Spec(A_0)}(M)$;
there is hence also an action on the open complement $\Spec(A) \smallsetminus V(A_\dagger)$. This action stabilizes each open subscheme $\Spec(A_f)$, and by~\cite[Lemma~2.1]{BS03} the morphism $\Spec(A_f) \to \Spec(A_{(f)})$ is a geometric quotient for the action of $\mathrm{D}_{\Spec(A_0)}(M)$ in the sense of~\cite[Definition~0.6]{MFK93}. It follows that~\eqref{eqn:geometric-quotient} is also a geometric quotient. 
%As explained in~\cite{BS03}, this morphism is also affine.

\begin{rmk}
\label{rmk:action-gp-scheme}
 Assume we are given a flat affine group scheme $H$ over $\Spec(A_0)$ and an action of $H$ on $A$ which preserves degrees. This induces an action on $\Spec(A)$, which stabilizes $V(A_\dag)$, hence also an action on $\Spec(A) \smallsetminus V(A_\dag)$. Since a geometric quotient is a categorical quotient (see~\cite[Proposition~0.1]{MFK93}), the composition
 \[
  H \times_{\Spec(A_0)} \bigl( \Spec(A) \smallsetminus V(A_\dag) \bigr) \longto \Spec(A) \smallsetminus V(A_\dag) \longto \Proj^M(A)
 \]
factors through a morphism
$H \times_{\Spec(A_0)} \Proj^M(A) \longto \Proj^M(A)$.
It is easily seen that this morphism defines an action of $H$ on $\Proj^M(A)$.
\end{rmk}

%--------------------------------------
\subsection{Open subschemes defined by quasi-relevant elements}
\label{ss:open-quasi-rel}
%--------------------------------------

Recall the notion of quasi-$M$-relevant element of $A$ from Definition~\ref{defquasirel}. If $a \in A$ is a nonzero $M$-quasi-relevant element, then we can write $a = \sum_{i=1}^n f_i$ where each $f_i$ is a nonzero $M$-relevant element of degree $\deg(a)$, and by Proposition~\ref{prop:magical}\eqref{it:magic-4} we have an open immersion
\[
 \Spec(A_{(a)}) \longto \bigcup_{i=1}^n D_\dag(f_i).
\]
By construction the right-hand side is an open subscheme in $\Proj^M(A)$, so we deduce an open immersion
\begin{equation}
\label{eqn:immersion-qrelevant}
 \Spec(A_{(a)}) \longto \Proj^M(A).
\end{equation}

\begin{lem}
\label{lem:open-qrelevant}
 The morphism~\eqref{eqn:immersion-qrelevant} is canonical, i.e.~it does not depend on the way $a$ is written as a sum of relevant elements.
\end{lem}

\begin{proof}
 Assume that $a=\sum_{i=1}^n f_i = \sum_{j=1}^m g_j$ where the $f_i$'s and the $g_j$'s are nonzero and relevant, all of the same degree. What we have to prove is that the compositions
 \[
  \Spec(A_{(a)}) \to \bigcup_{i=1}^n \Spec(A_{(f_i)}) \subset \Proj^M(A) \quad \text{and} \quad \Spec(A_{(a)}) \to \bigcup_{j=1}^m \Spec(A_{(g_j)}) \subset \Proj^M(A)
 \]
coincide. Here the first, resp.~second, immersion is obtained using the open covering
\[
\Spec(A_{(a)}) = \bigcup_{i=1}^n \Spec(A_{(a f_i)}), \quad \text{resp.} \quad
\Spec(A_{(a)}) = \bigcup_{j=1}^m \Spec(A_{(a g_j)})
\]
and the natural open immersions $\Spec(A_{(a f_i)}) \to \Spec(A_{(f_i)})$, resp.~$\Spec(A_{(a g_i)}) \to \Spec(A_{(g_i)})$.
However, in $A_{(a)}$ we have $\sum_{i,j} \frac{f_i g_j}{a^2}=1$, so that using similar considerations we can refine both of these decompositions into a decomposition
\[
\Spec(A_{(a)}) = \bigcup_{\substack{1 \leq i \leq n \\ 1 \leq j \leq m}} \Spec(A_{(a f_i g_j)}),
\]
and it is clear that both immersions can be obtained from this decomposition using the canonical open immersions $\Spec(A_{(a f_i g_j)}) \to \Proj^M(A)$.
%
%\[
%\Spec(A_{(a)}) = \bigcup_{\substack{1 \leq i \leq n \\ 1 \leq j \leq m}} D \left( \frac{f_i g_j}{a^2} \right)
%\]
%by~\cite[\href{https://stacks.math.columbia.edu/tag/01HS}{Tag 01HS}]{stacks-project}. Here, as in Proposition~\ref{prop:magical}\eqref{it:magic-2} , for any $i,j$ the open subset $D \left( \frac{f_i g_j}{a^2} \right)$ identifies with $\Spec(A_{(a f_i g_j)})$.
%
%both of these morphisms identify with the composition
%\[
%\Spec(A_{(a)}) \to \bigcup_{\substack{1 \leq i \leq n \\ 1 \leq j \leq m}} \Spec(A_{(f_i)}) \cap \Spec(A_{(g_j)}) \subset \Proj^M(A)
%\]
%obtained from the decomposition
%\[
%\Spec(A_{(a)}) = \bigcup_{\substack{1 \leq i \leq n \\ 1 \leq j \leq m}} \Spec(A_{(a f_i g_j)})
%\]
%(which follows from the fact that $\sum_{i,j} f_i g_j$ is invertible in $A_a$) and the identification of the scheme $\Spec(A_{(a f_i g_j)})$ with an open subscheme in $\Spec(A_{(f_i)}) \cap \Spec(A_{(g_j)}) = \Spec(A_{(f_i g_j)})$.
%images identify with the open subscheme of
%\[
% \left( \bigcup_{i=1}^n \Spec(A_{(f_i)}) \right) \cup 
% \left( \bigcup_{j=1}^m \Spec(A_{(g_j)}) \right)
%\]
%obtained as the union of the open subsets $\Spec(A_{(a \cdot f_i)})$ ($i \in \{1, \dots, n\}$) and $\Spec(A_{(a \cdot g_j)})$ ($j \in \{1, \dots, m\}$).
\end{proof}

Lemma~\ref{lem:open-qrelevant} justifies that the image of~\eqref{eqn:immersion-qrelevant} can be denoted $D_\dag(a)$; this open scheme identifies canonically with $\Spec(A_{(a)})$. In case $a$ is relevant, this coincides with the open subscheme denoted similarly above.

\begin{cor}
\label{cor:covering-Proj}
Let $(a_i : i \in I)$ be relevant elements in $A$.
%, and consider the ideal $\langle a_i : i \in I \rangle$ they generate. 
%Assume that any relevant element in $A$ belongs to 
If
\[
 A_\dag \subset \Rad \left( \sum_{i \in I} A \cdot a_i \right),
\]
%can be written as a linear combination (with coefficients in $A_0$) of products of the $a_i$'s. 
then we have
\[
\Proj^M(A) = \bigcup_{i \in I} D_\dag(a_i).
\]
\end{cor}

\begin{proof}
Of course we can assume that the $a_i$'s are nonzero.
By~\eqref{eqn:singletonproj}, $\Proj^M(A)$ is covered by the open subschemes $D_\dag(f)$ where $f$ is relevant and not nilpotent. By assumption, for such $f$ there exists $k \in \Z_{>0}$ such that $f^k = \sum_{j=1}^r b_j \cdot a_{i_j}$ for some nonzero elements $b_j \in A$ and some indices $i_j \in I$. Since $f$ is homogeneous, here we can assume that for any $j$ the element $b_j$ is homogeneous, with $\deg(b_j) + \deg(a_{i_j})=\deg(f)$. Then we have
\[
D_\dag(f) \overset{\eqref{eqn:open-power}}{=} D_\dag(f^k) \subset \bigcup_{j=1}^r D_\dag(b_j a_{i_j}),
\]
and by Lemma~\ref{lem:prod-relevant}\eqref{it:prod-relevant-0} we have $D_\dag(b_j a_{i_j}) \subset D_\dag(a_{i_j})$ for any $j$. The claim follows.
\end{proof}

\begin{rmk}
\begin{enumerate}
 \item
 Alternatively, the corollary can be proved by remarking that our assumptions imply that $V(A_\dag) = \bigcap_{i \in I} V(A \cdot a_i)$, so that $\Spec(A) \smallsetminus V(A_\dag) = \bigcup_{i \in I} D(a_i)$, and then using the map~\eqref{eqn:geometric-quotient}.
 \item
  In particular, in case $I$ can be chosen finite, this corollary implies that $\Proj^M(A)$ is quasi-compact.
\end{enumerate}
\end{rmk}

% \defi 
% \label{potionquasi}
% We will use the following terminology. 
% \begin{enumerate} 
% \item
% \label{it:potionquasiA}
% Let $a$ be a quasi-$M$-relevant element of $A$ (cf. Definition \ref{defquasirel}). Then by Proposition~\ref{prop:magical}\eqref{it:magic-4} we have a canonical open immersion $\Spec(A_{(a)}) \to \Proj^M(A)$, whose image will be denoted $D_{\dagger} (a)$. In other words, potions associated with quasi-relevant elements in $A$ also give rise to open subschemes of $\Proj^M(A)$.
% \item
% \label{it:potionquasiB} 
% If $\calF $ is a set of relevant elements of $A$, then we say that an element is \emph{$\calF$-quasi-$M$-relevant} (or just $\calF$-quasi-relevant)  if it is a sum of relevant elements in $\calF$ of the same degree. Then, by Proposition \ref{prop:magical}\eqref{it:magic-4}, $D_{\dagger} (a)$ is an open subscheme of $\Proj_{\calF} (A)$.
% \end{enumerate}
% \xdefi

%\pf 
%Because $A_{(f)}=A_{(f^\bbN)}= A_{(\underline{f^\bbN})}= A_{(\underline{S})}=A_{(S)}$.
%\xpf 

%Remark~\ref{remasingletonproj} shows that the scheme $\Proj^M(A)$ from Construction~\ref{defproj} coincides with the scheme defined by Brenner--Schr\"oer in~\cite[Definition~2.2]{BS03}. 

%--------------------------------------
\subsection{More basic properties}
\label{ss:prop-Proj}
%--------------------------------------

% The following result 
% % will be useful, it 
% generalizes a similar known result in the $\bbN$-graded case (see item~(3) in~\cite[\href{https://stacks.math.columbia.edu/tag/01HS}{Tag 01HS}]{stacks-project}, on which the proof is based). 
%However our argument is simpler and avoids descriptions of $\Proj$ in terms of homogeneous ideals of $A$.

We start with a functoriality property of the Proj scheme construction, which generalizes~\cite[\href{https://stacks.math.columbia.edu/tag/01MY}{Tag 01MY}]{stacks-project}. For this we note that if $\Psi: A \to B$ is a homomorphism of commutative $M$-graded rings and $S \subset A$ is a homogeneous multiplicative subset, then so is $\Psi(S)$. Clearly, $\Psi(S)$ is relevant if $S$ is, and finitely generated as a submonoid of $(B,\times)$ if $S$ is finitely generated as a submonoid of $(A,\times)$. Hence $\Psi$ induces a map $\mathcal{F}_A \to \mathcal{F}_B$, which will also be denoted $\Psi$.

\prop[(Functoriality of Proj)]
\label{foncto} 
Let $\Psi:A \to B $ be a homomorphism of $M$-graded rings. 
%\begin{enumerate}
%\item
For any $\mathcal{F} \subset \mathcal{F}_A$ 
%be a set of finite relevant families of $A$. Then $\Psi (\mathcal{F})$ is a set of relevant families of $B$ and 
we have a canonical morphism of schemes $\mathrm{Proj}^M_{\Psi (\mathcal{F})} (B) \to \mathrm{Proj}^M_{\mathcal{F} } (A)$. Moreover, for any $S \in \calF$ we have
\begin{equation}
\label{eqn:Proj-morph-affine}
\Proj^M_{\Psi(\calF)}(B) \times_{\Proj_{\calF}^M(A)} D_\dag(S) = D_\dag(\Psi(S));
\end{equation}
in particular, the morphism is affine.
%\item
%In case $\calF=\calF_A$, the morphism $\mathrm{Proj}^M_{\Psi (\mathcal{F}_A)} (B) \to \mathrm{Proj}^M(A)$ is affine.
%\end{enumerate}
\xprop

\pf 
The morphism is obtained by glueing from the morphisms induced by the canonical homomorphisms $A_{(S)} \to B_{(\Psi(S))}$ for $S \in \mathcal{F}$.
%If $f$ is a finite relevant family of $A$, then $\Psi (f) $ is finite relevant in $B$ and we have a morphism $A_{(f)} \to B _{(\Psi (f))}$. We obtain the desired morphism by glueing. $\qed$

Now fix $S \in \mathcal{F}$. From the definition of our morphism we see that we have an inclusion
\[
D_\dag(\Psi(S)) \subset \Proj^M_{\Psi(\calF)}(B) \times_{\Proj_{\calF}^M(A)} D_\dag(S)
\]
as open subschemes of $\Proj^M_{\Psi(\calF)}(B)$. On the other hand,
by definition $\Proj^M_{\Psi(\calF)}(B)$ is covered by the open subschemes $D_\dag(\Psi(T))$ where $T$ runs over $\mathcal{F}$, so that $\Proj^M_{\Psi(\calF)}(B) \times_{\Proj_{\calF}^M(A)} D_\dag(S)$ is covered by the open subschemes
$D_\dag(\Psi(T)) \times_{\Proj_{\calF}^M(A)} D_\dag(S)$.
Now the morphism $D_\dag(\Psi(T)) \to \Proj_{\calF}^M(A)$ factors through $D_\dag(T)$, so that
\[
D_\dag(\Psi(T)) \times_{\Proj_{\calF}^M(A)} D_\dag(S) = D_\dag(\Psi(T)) \times_{D_\dag(T)} (D_\dag(T) \cap D_\dag(S)).
\]
Using Proposition~\ref{prop:magical}\eqref{it:magic-2} one sees that the right-hand side coincides with $D_\dag(\Psi(ST))$; in particular it is contained in $D_\dag(\Psi(S))$, which finishes the proof of~\eqref{eqn:Proj-morph-affine}.

The equalities~\eqref{eqn:Proj-morph-affine} imply that our morphism is affine, see e.g.~\cite[\href{https://stacks.math.columbia.edu/tag/01S8}{Tag 01S8}]{stacks-project}.
\xpf

\begin{lem}
\label{lem:closed}
 In the setting of Proposition~\ref{foncto}, assume that $\Psi : A \to B$ is surjective. Then we have $\mathrm{Proj}^M_{\Psi (\mathcal{F}_A)} (B) = \mathrm{Proj}^M (B)$, and the canonical morphism
 $\Proj^M(B) \to \Proj^M(A)$
is a closed immersion.
\end{lem}

\begin{proof}
Let $S \in \calF_B$. By definition, there exist $s_1, \dots, s_r \in \underline{S}$ whose degrees generate a subgroup of $M$ of finite index. Let $S' \subset B$ be the homogenenous multiplicative subset generated by $s_1, \dots, s_r$, and let $S''=S \cdot S'$. Then $S'$ and $S''$ are finitely generated relevant homogeneous multiplicative subsets of $B$, and by Remark~\ref{rmk:Proj-opens}\eqref{it:Proj-opens-3} we have $D_\dag(S)=D_\dag(S'')$. We therefore have an embedding
% \[
%  D_\dag(S) \leftarrow D_\dag(S'') \to D_\dag(S').
% \]
% Moreover the morphism $A_{S} \to A_{S''}$ is an isomorphism, to that the left-hand arrow is an equality. We deduce an immersion 
$D_\dag(S) \to D_\dag(S')$. Now, choosing for each $i \in \{1, \dots, r\}$ a homogeneous preimage $t_i$ of $s_i$ in $A$ and denoting by $\tilde{S} \subset A$ the homogeneous multiplicative subset generated by $t_1, \dots, t_r$, it is clear that $\tilde{S}$ is relevant and that $\Psi(\tilde{S})=S'$. Hence $D_\dag(S') \subset \mathrm{Proj}^M_{\Psi (\mathcal{F}_A)} (B)$, which implies that $D_\dag(S)$ is also contained in $\mathrm{Proj}^M_{\Psi (\mathcal{F}_A)} (B)$. Since $S$ was arbitrary, this proves the equality $\mathrm{Proj}^M_{\Psi (\mathcal{F}_A)} (B) = \mathrm{Proj}^M (B)$.

Since the property of being a closed immersion is local on the target (see~\cite[\href{https://stacks.math.columbia.edu/tag/02L6}{Tag 02L6}]{stacks-project}), to prove the second assertion it suffices to prove that for any $S \in \calF_A$ the induced morphism
\[
 \Proj^M(B) \times_{\Proj^M(A)} D_\dag(S) \to D_\dag(S)
\]
is a closed immersion. This follows from~\eqref{eqn:Proj-morph-affine}, since the homomorphism $A_{(f)} \to B_{(\Psi(f))}$ is clearly surjective.
%
%, in view of~\eqref{eqn:singletonproj} it suffices to prove that for any $f \in A$ relevant the morphism
%\[
% \Proj^M(B) \times_{\Proj^M(A)} D_\dag(f) \to D_\dag(f)
%\]
%is a closed immersion. However, we claim that 
%\begin{equation}
%\label{eqn:Proj-quotient-open}
%\Proj^M(B) \times_{\Proj^M(A)} D_\dag(f) = D_\dag(\Psi(f)),
%\end{equation}
%which will finish the proof since the morphism $A_{(f)} \to B_{(\Psi(f))}$ is surjective. In fact, 
%%by Proposition~\ref{prop:image-relevant}\eqref{it:image-relevant-surj} 
%by Remark~\ref{rmk:preimage-Adag-surj}
%the preimage in $\Spec(B)$ of $A \smallsetminus V(A_\dag)$ is $B \smallsetminus V(B_\dag)$, and the diagram
%\[
% \xymatrix{
% B \smallsetminus V(B_\dag) \ar[r] \ar[d] & A \smallsetminus V(A_\dag) \ar[d] \\
% \Proj^M(B) \ar[r] & \Proj^M(A)
% }
%\]
%commutes, where the vertical arrows are the canonical morphisms considered in~\S\ref{subsectprojring} are the lower arrow is the morphism considered above. The preimage in $A \smallsetminus V(A_\dag)$ of $D_\dag(f)$ is $\Spec(A_f)$, whose preimage in $B \smallsetminus V(B_\dag)$ is $D_\dag(B_{\Psi(f)})$, which implies~\eqref{eqn:Proj-quotient-open}.
\end{proof}

\begin{rmk}
 The conclusion of Lemma~\ref{lem:closed} can often be obtained under weaker assumptions using more specific information on relevant elements in $A$. See~\cite[\href{https://stacks.math.columbia.edu/tag/01N0}{Tag 01N0}]{stacks-project} for the case of $\N$-graded rings, and Proposition~\ref{prop:vb-flag-Proj} below for an example with a more general $M$. 
\end{rmk}

 \prop 
 \label{produitvstensoriel}
 Let $M$ and $M'$ be two finitely generated abelian groups. Let $R$ be a commutative ring and let $A$ (resp.~$A'$) be a commutative $M$-graded (resp.~$M'$-graded) $R$-algebra. 
 Then for the natural $(M \times M')$-grading on $A \otimes_R A'$, we have a canonical isomorphism
 %Then $A \otimes _R A'$ is a commutative $M \times M'$-graded $R$-algebra and 
 \[
 \Proj ^{M \times M'} (A \otimes_R A') \cong \Proj ^M (A) \times_{\Spec(R)} \Proj ^{M'} (A'). 
 \]
\xprop 

\pf
 The identification is provided by the following equalities (to be explained below):
\begin{align*}
\Proj ^M ( A ) \times_{\Spec (R) } \Proj ^{M'} ( A')&= \left( \bigcup_{f\in A \text{ relevant} }\Spec (A_{(f)}) \right) \times_{\Spec (R)} \left( \bigcup_{f' \in A' \text{ relevant} } \Spec ( A_{(f')})\right) \\
&= \bigcup_{f \in A, \, f' \in A' \text{ relevant} } \Spec \left( A_{(f)} \otimes_R A'_{(f')} \right) \\
&= \bigcup_{f \in A, \, f' \in A' \text{ relevant} } \Spec \left( (A\otimes_R A')_{(f \otimes f')} \right) \\
& = \Proj ^{M \times M'}( A \otimes_{R} A'). 
\end{align*}
Here the first equality follows from~\eqref{eqn:singletonproj}.
%is immediate from the definition of $\Proj$ schemes and Remark~\ref{rmk:def-Proj}\eqref{it:singletonproj}. 
The second equality follows from a basic property of fiber products of schemes, see~\cite[\href{https://stacks.math.columbia.edu/tag/01JS}{Tag 01JS}]{stacks-project}. 
The third equality follows from the obvious isomorphism of $(M \times M')$-graded rings
$A_{f} \otimes_R A'_{f'} \cong (A\otimes_R A')_{f \otimes f'}$
(for $f \in A$ and $f' \in A'$ homogeneous)
by restriction to the components of degree $(0,0)$.

%To prove the third equality, note that the two $R$-morphims 
%$A_{(f)} \otimes _R A'_{(f')} \to (A \otimes _R A') _{(f \otimes f')}, \frac{a}{f^k} \otimes \frac{a'}{{f'}^{k'}} \mapsto \frac{af^{k'} \otimes a' {f'}^k }{(f\otimes f')^{k+k'}}$
% and $(A \otimes _R A')_{(f \otimes f')} \to A_{(f)} \otimes _R A'_{(f')}, \frac{a \otimes a'}{(f\otimes f') ^k} \mapsto \frac{a}{f^k} \otimes \frac{a'}{f'^k}$ are well-defined and mutual inverses. We now explain 
 
 To conclude we have to explain the fourth equality. Let $\calF_{A \otimes_R A'}^{\otimes} \subset \calF_{A \otimes_R A'}$ be the subset consisting of homogeneous multiplicative subsets of the form $(f \otimes f')^{\N}$ 
% family of $(M \times M')$-relevant elements of $A \otimes_R A'$ made of pure tensors $f \otimes f'$ 
 with $f \in A$ $M$-relevant and $f' \in A'$ $M'$-relevant. 
% (Here we use once again the notation introduced in Definition~\ref{def:rele}\eqref{it:relevant-element}.) 
 The desired equality will follow from the equality
 \[
 \Proj_{\calF_{A \otimes_R A'}^{\otimes}}^{M \times M'}( A \otimes_{R} A')=\Proj ^{M \times M'}( A \otimes_{R} A'),
 \]
 which can be justified as follows.
 %We now prove this last equality. The inclusion $\Proj_{\calF_{\otimes}}^{M \times M'}( A \otimes_{R} A')\subset \Proj ^{M \times M'}( A \otimes_{R} A')$ is trivial. 
 By definition, the left-hand side is an open subset of the right-hand side. To prove that this open subset is all of $\Proj ^{M \times M'}( A \otimes_{R} A')$, by~\eqref{eqn:singletonproj} it suffices to prove that it contains the open subscheme attached to any relevant element $x $ in $ A \otimes_R A'$.
% To prove the other inclusion, let us take a $M \times M'$-relevant element $x $ in $ A \otimes_R A'$. 
 Now if $x \in A \otimes_R A'$ is relevant,
 by Proposition \ref{radical}\eqref{it:radic-1} there exist a positive integer $k$, a finite set $J$ and for any $j \in J$ an $M$-relevant element $s_j$ in $A$ and an $M'$-relevant element $s'_j$ in $A'$ such that $x^k = \sum_{j \in J} s_j \otimes s'_j$. Then we have
 \[
 D_\dag(x) = D_\dag(x^k) \subset \bigcup_{j \in J} D_\dag(s_j \otimes s'_j)
 \]
 by~\eqref{eqn:open-power} and the constructions of~\S\ref{ss:open-quasi-rel},
 which finishes the proof.
 %So $x^k$ is $\calF_{A \otimes_R A'}^{\otimes}$-quasi-relevant. Now Definition \ref{potionquasi}\eqref{it:potionquasiB} shows that $D_\dagger (x^k) \subset \Proj_{\calF_{\otimes}}^{M \times M'}( A \otimes_{R} A').$ This finishes to explain the inclusion $\Proj_{\calF_{\otimes}}^{M \times M'}( A \otimes_{R} A')\supset \Proj ^{M \times M'}( A \otimes_{R} A')$ since $D_\dagger (x^k)= D_\dagger (x)$ by Fact \ref{powerpotionequal}.
% The proof is finished. Another proof of Proposition \ref{produitvstensoriel} might be obtained using Proposition \ref{radical}\eqref{it:radic-2} and Description \ref{geodescri}.
\xpf 

 \begin{ex}
 Let $A=\bbZ [X_1, \ldots, X_{n+1}]$ be $\bbZ$-graded by $\deg(X_1^{k_1}\cdots X_{n+1} ^{k_{n+1}} ) = \sum_{j=1}^{n+1} k_j $. Let $A'=\bbZ [X_1, \ldots , X_{n'+1}]$ be similarly $\bbZ$-graded. We have 
 $\Proj ^\bbZ (A) =\bbP ^n$ and $ \Proj ^\bbZ (A') = \bbP^{n'}$. Proposition \ref{produitvstensoriel} implies that $\bbP^n \times _{\Spec (\bbZ)} \bbP^{n'}= \Proj ^{\bbZ ^2}(A \otimes _{\bbZ} A' ) $ 
 \end{ex}

\begin{rmk}
Another proof of Proposition~\ref{produitvstensoriel} can be obtained using Proposition~\ref{radical}\eqref{it:radic-2} and the fact that~\eqref{eqn:geometric-quotient} is a geometric quotient.
%Description \ref{geodescri}.
\end{rmk}

The following corollary appears in the proof of \cite[Proposition 2.5]{BS07}.

\coro[(Base change)]
\label{propaffinebasechange}
Let $M$ be a finitely generated abelian group, let $R$ be a commutative ring, and let $A$ be a commutative $M$-graded $R$-algebra. If $R'$ is a commutative $R$-algebra, then for the natural $M$-grading on $A \otimes_R R'$ we have
\[
\Proj^M (A \otimes_R R' ) \cong \Proj^M (A) \times_{\mathrm{Spec}(R)} \mathrm{Spec} (R').
\]
\xcoro

\begin{proof}
By Example~\ref{ex:Proj-Spec}, the statement is the special case of Proposition~\ref{produitvstensoriel} when $M'=\{0\}$.
\end{proof}

\begin{rmk}
\label{rmk:open-basechange}
Let $S \in \calF_A$, and let $S_{R'}$ be its image in $A \otimes_R R'$. Then it is clear that under the isomorphism of Corollary~\ref{propaffinebasechange} the open subscheme $D_\dag(S_{R'}) \subset \Proj^M (A \otimes_R R' )$ identifies with $D_\dag(S) \times_{\mathrm{Spec}(R)} \mathrm{Spec} (R')$.
\end{rmk}

The following example generalizes~\cite[\href{https://stacks.math.columbia.edu/tag/01MI}{Tag 01MI}]{stacks-project}.

\exam
\label{triv}
Let $R$ be a commutative ring and $N$ be a commutative cancellative monoid such that the group $M= N^{\mathrm{gp}}$ is finitely generated. Let $R[N]:= \bigoplus_{m \in N} R \cdot X^m$ be the $R$-algebra of the monoid $N$, with its natural $M$-grading. The homogeneous elements in this ring are those of the form $r X^m$ with $r \in R$ and $m \in N$; such an element is relevant if and only if $N \cap (m-N)$ generates $M$ up to torsion. For each such element we have $(R[N])_{(r X^m)} = R_r$; therefore we have $\mathrm{Proj} (R[N])=\mathrm{Spec} (R) $. 
%(This example generalizes~\cite[\href{https://stacks.math.columbia.edu/tag/01MI}{Tag 01MI}]{stacks-project}.)
\xexam

%\fact \label{M0proj} Assume that $M=0$. Then $\Proj (A) = \Spec (A)$.
%\xfact
%\pf
%We have $\Proj (A) = \cup _{a \in A} \Spec (A_{(a)})= \cup _{a \in A} \Spec (A_a)$. So $\Proj (A) = \Spec (A)$.
%\xpf 

It is a standard fact that, even in the $\N$-graded setting, very different graded rings can have the same Proj scheme, see e.g.~the discussion of the Veronese embedding in~\cite{eg}. We generalize this observation in the following proposition.

\prop
\label{thin}
Let $M$ be a finitely generated abelian group, and let $A$ be a commutative $M$-graded ring.
Let $M' \subset M$ be a subgroup of finite index. Consider the subring $A'= \bigoplus_{m \in M'} A_{m}$; it can be considered naturally as an $M'$-graded ring, or as an $M$-graded ring.
We have canonical isomorphisms of schemes
\[
\Proj^M(A) \cong \Proj^{M} (A')\cong \Proj^{M'} (A').
\]
\xprop

\pf 
We use the description of these schemes using~\eqref{eqn:singletonproj}.
%arising from Remark~\ref{rmk:def-Proj}\eqref{it:singletonproj}.
If $f$ is an $M$-relevant element in $A$, then as in~\eqref{eqn:open-power}
%Fact~\ref{powerpotionequal} 
we have $A_{(f)}= A_{(f^d)}$ for all $d \in \bbZ_{>0}$. If $n$ is the index of $M'$ in $M$, 
%$M/M'$ is finite there exists an integer $n_f>0$ such that 
then $f^{n}$ belongs to $A'$ and is $M'$-relevant in $A'$; moreover we have $A_{(f)}= A_{(f^{n})}= A'_{(f^{n})}$. This shows that $\Proj^M(A)$ identifies with an open subscheme of $\Proj^{M'} (A')$. Reciprocally, if $f \in A'$ is $M'$-relevant, then $f$ is an $M$-relevant element of $A$; it follows that this open subscheme is the whole of $\Proj^{M'} (A')$. 
%This shows that $\Proj^{M'} (A') $ is an open subscheme of $\Proj ^{M} (A)$. 

The other isomorphism $\Proj^{M} (A')\cong \Proj^{M'} (A')$ is clear since
%, because of our assumption that $M/M'$ is finite, 
an element in $A'$ is $M'$-relevant if and only if it is $M$-relevant.
\xpf 

% \rema
% %\begin{enumerate}
% %\item
% Note that Proposition \ref{thin} is not true if $M/M'$ is infinite. For example, let $R$ be a ring, which we endow with the $\bbZ$-grading where $R$ is in degree $0$. One can also consider $R$ as a $\{0\}$-graded. There is no $\bbZ$-relevant element in $R$ so $\Proj ^\bbZ (R) = \varnothing$; however by Example~\ref{ex:Proj-Spec} we have $\Proj^{\{0\}}(R) = \Spec (R)$.
% %\item
% %\end{enumerate}
% \xrema

We now state a few immediate (but useful) consequences of Proposition~\ref{thin}. Here, intuitively, \eqref{it:cor-large} says that ``the Proj scheme only depends on components in any subgroup of finite index,'' and~\eqref{it:cor-free} allows to reduce the theory to the case $M$ is a free $\Z$-module. Finally,~\eqref{it:cor-M-finite} generalizes Example~\ref{ex:Proj-Spec}.

\coro 
%\phantomsection
\label{cor:isom-large-degrees}
Let $M$ be a finitely generated abelian group.
\begin{enumerate}
\item
\label{it:cor-large}
Let $A$ and $B$ be commutative $M$-graded rings, and assume that there exists a subgroup $M' \subset M$ such that $M / M'$ is finite and the subrings of $A$ and $B$ consisting of their homogeneous components associated with elements in $M'$ are isomorphic (as $M'$-graded rings). Then we have
$\Proj^M(A) \cong \Proj^M (B)$.
%This technique applies in particular when $M=\Z^r$ and $M'=d M$ for $d \in M$ with all of its components nonzero.
\item
\label{it:cor-free}
Let $A$ be a commutative $M$-graded ring.
Let $M' \subset M$ be a complement to the torsion subgroup $M_{\mathrm{tor}} \subset M$, and consider the $M'$-graded ring $A' = \bigoplus_{m \in M'} A_{m}$. Then we have $\Proj^M(A) \cong \Proj^{M'}(A')$.
%Write $M = \bbZ ^r \times D$ where $D$ is a finite abelian group of order $d$. Then $dM $ is torsion free and $\Proj (A) = \Proj (\bigoplus _{k \in dM} A_{k})$.
\item 
\label{it:cor-M-finite}
Assume that $M$ is finite. For any commutative $M$-graded ring, we have $\Proj (A) = \Spec (A_0)$. 
\end{enumerate}
\xcoro

%\pf
%\begin{enumerate}
%\item  This is an immediate corollary of Proposition \ref{thin}.
%\item Let $d$ be $\# M$. Then Proposition \ref{thin} implies that $\Proj (A) = \Proj ( \bigoplus _{k \in dM} A_{k} )$. Now $dM=0$, so $\Proj (A)= \Proj (A_0)$ where $A_0$ is $0$-graded. We now apply Fact \ref{M0proj} and get $\Proj (A_0) = \Spec (A_0)$.  
%\end{enumerate}
%\xpf 

%--------------------------------------
\subsection{Relative Proj schemes}
\label{ss:relative}
%--------------------------------------

Proposition~\ref{propaffinebasechange} allows us to extend the construction of Proj schemes to a relative setting, as follows. Let $X$ be a scheme, $M$ be a finitely generated abelian group, and $\mathcal{A}$ be an $M$-graded quasi-coherent $\mathcal{O}_X$-algebra. Note that for any quasi-compact open subscheme $U \subset X$ the algebra $\Gamma(U,\mathcal{A})$ admits a canonical $M$-grading by~\cite[\href{https://stacks.math.columbia.edu/tag/01AI}{Tag 01AI}]{stacks-project}.

\begin{const}[(Relative Proj)]
\label{relative}
For any affine open subscheme $U \subset X$ we consider the Proj-scheme $\Proj^M(\Gamma(U,\mathcal{A}))$ of the $M$-graded ring $\Gamma(U,\mathcal{A})$; it admits a canonical morphism to $\Spec(\Gamma(U,\mathcal{A}_0))$, hence to $U$. If $U' \subset U \subset X$ are open affine subschemes, by Corollary~\ref{propaffinebasechange} we have a canonical isomorphism of $U'$-schemes.
\[
U' \times_U \Proj^M(\Gamma(U,\mathcal{A})) \cong \Proj^M(\Gamma(U',\mathcal{A})).
\]
By relative glueing (see~\cite[\href{https://stacks.math.columbia.edu/tag/01LH}{Tag 01LH}]{stacks-project}), there exists an $X$-scheme $\mathrm{Proj}_X^M(\mathcal{A})$ such that
%There exists an $S$-scheme $\mathrm{Proj}_S^M(\mathcal{A})$ 
%%and a morphism of schemes $p: \mathrm{Proj}_S^M (\mathcal{A}) \to S$ 
%such that for every affine open subscheme $U \subset S$, there exists a canonical isomorphism
\[
\mathrm{Proj}_X^M(\mathcal{A}) \times_X U = \mathrm{Proj}^M (\mathcal{A} (U))
\]
for any affine open subscheme $U \subset X$.
%and
% for $U \subset U' \subset S $ affine open, the composition
%$\mathrm{Proj} (\mathcal{A} (U)) \cong p^{-1} (U) \to p^{-1}(U') \cong \mathrm{Proj} (\mathcal{A} (U'))  $ is an open immersion.
\end{const}

Most of the results of this section have obvious analogues in this setting; we leave it to the reader to formulate these variants and adapt the proofs.

% \prop 
% Let $S$ be a scheme. Let $M$ and $M'$ be two finitely generated abelian groups. Let $\calA$ (resp. $\calA '$) be a $M$-graded (resp. $M'$-graded) quasi-coherent $\calO_S$-algebra. Then for the natural $M \times M'$-grading on the $\calO_S$-algebra $\calA \otimes_{\calO_S} \calA'$, we have a canonical isomorphism 
% \[ 
% \mathrm{Proj}_S^{M} (\calA) \times_S \mathrm{Proj}_S^{M'} (\calA') \cong \mathrm{Proj}_S^{M \times M'} (\calA \otimes_{\calO_S} \calA').
% \]
% \xprop
% 
% \pf 
% This follows from Proposition \ref{produitvstensoriel}.
% \xpf 

%-------------------------------------------
\subsection{Proj schemes, blowups and dilatations} 
\label{ss:dilatations}
%-------------------------------------------

Let $A$ be a ring and $\{J_i: i \in I \}$ be a family of ideals of $A$ parametrized by a finite set $I$. Recall the ring $\Bl_{\{J_i : i \in I\}} A= \bigoplus_{\nu \in \mathbb{N}^I} J_i^\nu$ (see Example \ref{exmulti} and \cite[Definition 2.34]{Ma24}), which is $\mathbb{Z}^I$-graded. 
%We assume now that $I $ if finite.
 Any family $\{ a_i : i \in I \}$ of elements of $A$ such that $a_i \in J_i$ for each $i \in I$ can be considered as a relevant family in the $\mathbb{Z}^I$-graded ring $\Bl _{\{J_i : i \in I\}} A$, so that we have an associated affine open subscheme
\[
D_\dag(\{a_i : i \in I \}) \subset \Proj^{\Z^I} ( \Bl_{\{J_i : i \in I\}} A), 
\]
and this open subscheme identifies with the dilatation ring $\Spec(A[\{ \frac{J_i}{a_i} : i \in I \}])$ by Example \ref{exmulti}.

The goal of the present subsection is to study related connections and compatibilities between dilatations and these specific Proj schemes, 
and to show that this allows to easily provide a proof of the existence of (projective) blowups, possibly with multiple centers. This section was announced in \cite[\S~3.9]{Ma24}.
We start with a multi-centered version of~\cite[\href{https://stacks.math.columbia.edu/tag/0804}{Tag 0804}]{stacks-project}.

\begin{lem}
\label{lem:covering-blowup}
In the setting above, the scheme $\Proj^{\Z^I} ( \Bl_{\{J_i : i \in I\}} A)$ is covered by the open subschemes $D_\dag(\{a_i : i \in I \})$ where $\{a_i : i \in I\}$ runs over the families such that $a_i \in J_i$ for all $i \in I$.
\end{lem}

\begin{proof}
By Corollary~\ref{cor:covering-Proj}, the lemma will follow if we prove that any non nilpotent relevant element in the $\Z^I$-graded ring $\Bl_{\{J_i : i \in I\}} A$ 
%has a power which 
belongs to the sum of the ideals $(\Bl_{\{J_i : i \in I\}} A) \cdot \prod_i a_i$ where $\{ a_i : i \in I \}$ runs over the families as above. Now it is clear that this sum is the sum of the graded components in $\Bl_{\{J_i : i \in I\}} A$ associated with elements in $\Z^I$ all of whose components are positive. The claim therefore follows from the observation that the degree of any non nilpotent relevant element in $\Bl_{\{J_i : i \in I\}} A$ must have positive components, because if the $i$-component of the degree of a non nilpotent homogeneous element is $0$, then the same will hold for all divisors of all of its powers, so that this element cannot be relevant.
\end{proof}

We now work in a relative setting; we therefore fix a base scheme $X$, and let $\{Y_i : i \in I\}$ be a family of closed subschemes parametrized by a finite set $I$. We first recall the definition of blowups.
 
 \defi[(Multi-centered blowups)]
 \label{blowdef} 
 A (projective) blowup of $X$ with center $\{Y_i : i \in I\}$  is a pair consisting of a scheme $\widetilde{X}$ and a morphism $\pi : \widetilde{X} \to X$ such that $\pi^{-1} (Y_i)$ is an effective
Cartier divisor for all $i \in I$, and which is universal with respect to this property in the following sense. If $\pi ': \widetilde{X}'\to X$ is
any morphism such that $(\pi')^{-1} (Y_i)$ is an effective Cartier divisor for all $i $, then there exists a unique
morphism $ m : \widetilde{X}' \to \widetilde{X} $ such that $\pi' = \pi \circ m$.
 \xdefi 
 
 Clearly, a blowup of $X$ with center $\{Y_i : i \in I\}$  is unique up to unique isomorphism if it exists. 
 For each $i$, let $\mathcal{J}_i \subset \calO _X$ be the quasi-coherent ideal which defines $Y_i$, and
 %, so that $Y_i = V (\mathcal{J}_i )$.
let
\[
\mathrm{Bl} _{\{ \mathcal{J} _i : i \in I \} } \mathcal{O} _X  = \bigoplus _{ \nu \in \mathbb{N}^I } \mathcal{J} ^\nu
\]
be the associated Rees algebra (cf. e.g. \cite[\S 3.1]{Ma24}). 
%(Here, if $\nu = (\nu_i : i \in I)$ then $\mathcal{J}^\nu$ is the product of the $(\mathcal{J}_i)^{\nu_i}$'s.) 
This is canonically a $\mathbb{Z}^I$-graded quasi-coherent $\mathcal{O} _X$-algebra.
Using the construction of~\S\ref{ss:relative} we now define an $X$-scheme 
\[
{\mathrm{Bl}}_{\{Y_i : i \in I \}} X := \Proj^{\mathbb{Z}^I}_X ( \mathrm{Bl} _{\{ \mathcal{J} _i : i \in I \} } \mathcal{O}_X  ).
\] 
We will prove in Proposition~\ref{exiblo} below that this scheme is a blowup of $X$ with center $\{Y_i : i \in I\}$ in the sense of Definition~\ref{blowdef}. (The case $|I|=1$ of this statement is classical, see e.g.~\cite[\href{https://stacks.math.columbia.edu/tag/01OF}{Tag 01OF}]{stacks-project}, but our proof is different.)

\begin{lem}
\label{lem:Bl-separated}
The canonical morphism ${\mathrm{Bl}}_{\{Y_i : i \in I\}} X \longto X$ is separated.
\end{lem}

\begin{proof} 
%This is local on the target. 
Since the property of being separated is local on the target (see e.g.~\cite[\href{https://stacks.math.columbia.edu/tag/02KU}{Tag 02KU}]{stacks-project}), we can assume $X=\Spec(A)$ is affine, and by~\cite[\href{https://stacks.math.columbia.edu/tag/01KV}{Tag 01KV}]{stacks-project} the claim will follow if in this case we prove that $
{\mathrm{Bl}}_{\{Y_i : i \in I\}} X$ is a separated scheme.
We proceed with the usual local notations from the beginning of the subsection.
Consider the collection of families $\{a_i : i \in I\}$ with  $a_i \in J_i$ for all $i \in I$, and let $\mathcal{F} \subset \mathcal{F}_A$ be the subset consisting of the multiplicative subsets generated by such families. By Lemma \ref{lem:covering-blowup}, we have $\Proj^{\mathbb{Z}^I} ( \Bl_{\{J_i : i \in I\}} A) = \Proj^{\mathbb{Z}^I}_{\calF }( \Bl_{\{J_i : i \in I\}} A)$. Using the notation from Proposition~\ref{prop:separated-family}, for any $T \in \calF $ we have $C_T = (\mathbb{R}_{\geq 0})^I \subset \mathbb{R}^I$. So if $T,T' $ belong to $\calF $, $C_T \cap C_{T'}$ has nonempty interior. Now Proposition \ref{prop:separated-family} finishes the proof.
\end{proof}

\prop[(Definiteness)]
\label{defined} 
If $C$ is an effective Cartier divisor on $X$, then $\left( 
{\mathrm{Bl}}_{\{Y_i : i \in I\}} X \right) \times_X C$ is an effective Cartier divisor on $
{\mathrm{Bl}}_{\{Y_i : i \in I\}} X$. 
\xprop

\pf
This statement is local on $X$, so we can assume $X=\Spec(A)$ for some ring $A$, and write $J_i$ for the global sections of $\calJ_i$. 

In view of Lemma~\ref{lem:covering-blowup}, in this case it suffices to prove that for any family $\{ a_i : i \in I \}$ with $a_i \in J_i$ for any $i$, 
%Now that this claim is proved, to prove the proposition it suffices to prove that 
if $C$ is an effective Cartier divisor on $X$ then $\Spec(A[\{ \frac{J_i}{a_i} : i \in I \}]) \times_X C$ is an effective Cartier divisor on $\Spec(A[\{ \frac{J_i}{a_i} : i \in I \}])$. This statement is a special case of~\cite[Proposition 3.18]{Ma24}.
\xpf

In the following statement we use the notation of dilatations of schemes from~\cite{Ma24} (see also \cite{MRR20,DMdS23}). From now on, for simplicity of notation we set
\[
X' := 
{\mathrm{Bl}}_{\{Y_i : i \in I \}} X, \quad
% \]
% and for $i \in I$ we set
% \[
Y_i' := Y_i \times _X X'.
\]

 \prop[(Complementary)] 
 \label{complementary}
   For all $i \in I $, $ Y_i '$ is an effective Cartier divisor on $X'$. Moreover, the composition $ \mathrm{Bl}^{\{Y_i' : i \in I \}}_{\{\varnothing : i \in I \}} X' \to X' \to X$ is an open immersion\footnote{Recall that if $X$ is a scheme, $\varnothing$ denotes the empty closed subscheme corresponding to the quasi-coherent ideal $\calO _X$. In particular $\mathrm{Bl}^{\{Y_i' : i \in I \}}_{\{\varnothing : i \in I \}} X'$ denotes the dilatation of $X'$ relatively to the multi-center $\{ [\varnothing, Y_i'] : i \in I\} $ ($Y_i' $ is a Cartier divisor in $X'$ by assumption). }; in fact, this is just the inclusion $X \smallsetminus (\cup_{i \in I } Y_i ) \to X$.
 \xprop
 
 \pf 
 This statement is local on $X$, so we can assume $X=\Spec(A)$ for some ring $A$, and write $J_i$ for the global sections of $\calJ_i$. 
 In this case, 
 by Lemma~\ref{lem:covering-blowup} the scheme
 %we have explained in the course of the proof of Proposition~\ref{defined} that
 $X'={\mathrm{Bl}}_{\{Y_i : i \in I\}} X$ is covered by affine open subschemes
$\Spec \big(( \Bl_{\{J_i : i \in I\}} A )_{(\{a_i : i \in I \})}  \big)$ where $\{a_i : i \in I\}$ is a family of elements of $A$ such that $a_i \in J_i$ for any $i$. 
%Fact \cite[2.29]{Ma24} implies that 
As explained above, $( \Bl_{\{J_i : i \in I\}} A )_{(\{a_i : i \in I \})}$ is the dilatation ring $A[ \{ \frac{J_i}{a_i} : i \in I \}]$. Now the first assertion follows from~\cite[Facts~2.10 \& 2.28]{Ma24}. The second assertion follows from the equalities
\begin{align*}
\big(A[\big\{ \frac{J_i}{a_i} :i \in I \big\}] \big)\Big[ \Big\{ \frac{A[\{ \frac{J_i}{a_i} : i \in I \}]}{a_i} : i \in I \Big\} \Big] &= \big(A[\big\{ \frac{J_i}{a_i} :i \in I \big\}] \big)\Big[ \Big\{ \frac{A}{a_i} : i \in I \Big\} \Big]\\& = A[\big\{\frac{J_i}{a_i} : i \in I\big\}\sqcup \{ \frac{A}{a_i} : i \in I\big\} ]  \\
&= A[\big\{ \frac{A}{a_i} : i \in I\big\} ]\\
& =A_{\{a_i\}_{i \in I}}
\end{align*}
(use e.g. \cite[2.22, 2.28, 2.11]{Ma24})
and the fact that $A_{\{a_i\}_{i \in I}}$ is the localization of $A$ with respect to the $a_i$'s.
%The third assertion follows from universal properties. 
%
% To treat the general case, one takes an affine open covering of $X$ and uses the affine case.
 \xpf
 
\prop[(Existence of blowups)] \label{exiblo}
 Let $\pi: T \to X$ be a scheme such that, for all $i \in I $, $T \times _X Y_i$ is an effective Cartier divisor on $T$. Then there exists a unique morphism of $X$-schemes $ T \to  X'$. As a consequence, the pair consisting of $X'$ and the canonical morphism $X'\to X$ is a blowup of $X$ with center $\{ Y_i : i \in I\}$.
\xprop

\pf
First we show that for $T$ and $\pi$ as in the statement, there exists at most one morphism of $X$-schemes $T \to X'$. In fact, let $\phi,\phi' : T \to X'$ be two morphisms of $X$-schemes. Denote by $U$ the open complement of the sum of the effective Cartier divisors $T \times _X Y_i$; then $U$ is scheme theoretically dense in $T$ by~\cite[\href{https://stacks.math.columbia.edu/tag/07ZU}{Tag 07ZU}]{stacks-project}. By Lemma~\ref{lem:Bl-separated},
the morphism $X' \to X$ is separated; as a consequence, the difference kernel $\mathrm{ker}(\phi , \phi')$ (see~\cite[Definition~9.1]{GW20}) is a closed subscheme of $T$ by \cite[Definition/Proposition~9.7]{GW20}. 
%Let $U \subset T$ be the open subscheme $\mathrm{Spec} (R_\Gamma)$ where $\Gamma$ is the multiplicative subset of $R$ generated by $\{\gamma _i : i \in I\}$.  Then the scheme theoretic image of the open immersion $\iota_U : U \to T$ is $T$ by~\cite[\href{https://stacks.math.columbia.edu/tag/056A}{Tag 056A}]{stacks-project}. 
%Let $\iota_U$ be the open immersion $U \to T $. 
Since 
%$\iota _U ^{-1} \pi^{-1} (Y_i) = \emptyset $ 
$U \times_X Y_i = \varnothing$
for all $i$, the universal property of dilatations (see~\cite[Proposition~3.17]{Ma24}) implies that both $\phi_{|U}$ and $ \phi'_{|U}$ factor through $\mathrm{Bl}_{\{\varnothing : i \in I \}}^{\{Y_i' : i\in I \}}  X'$. By Proposition~\ref{complementary}, $\mathrm{Bl}_{\{\varnothing : i \in I \}}^{\{Y_i' : i\in I \}}  X'$ is an open subscheme of $X$, so that this implies that $\phi_{|U}= \phi'_{|U}$. We deduce that $U \subset \mathrm{ker} ( \phi , \phi' )$, and then that $\mathrm{ker} ( \phi , \phi' )= T$, so that $\phi  = \phi'$.

%First, we reduce to the case: 
%It suffices to prove the statement in 
To prove existence we first consider the following setting. We assume that $X= \Spec (A)$ and $ T=\Spec(R)$ for some rings $A$ and $R$, and denote by $J_i$ the global sections of $\calJ_i$, so that $ Y_i = \Spec (A / J_i)$, and by $f:A \to R$ the homomorphism of rings
%, and by assumption $T \times_X Y_i$ is principal for all $i\in I$. 
corresponding to $\pi : T \to X$. With this notation, we further assume that, for any $i \in I$, we have $f(J_i)R = \gamma _i R$ for some non-zero-divisor $\gamma_i \in R$. 

Let us first show the existence of a morphism of $X$-schemes $T \to X'$.
The morphism $f$ induces a canonical morphism of $\mathbb{Z}^I$-graded rings
\[
\bigoplus_{ \nu \in \mathbb{N}^I} J^\nu \to \bigoplus_{\nu \in \mathbb{N}^I }(f(J)R)^\nu
\]
where $(f(J)R)^\nu$ denotes the ideal $\prod_{i \in I }(f(J_i)R)^{\nu_i } \subset R$. Now since $f(J_i)R = \gamma_i R$ with $\gamma_i$ non-zero-divisor for any $i$, we have a canonical isomorphism of $R$-algebras 
$\bigoplus_{\nu \in \mathbb{N}^I }(f(J)R)^\nu \cong R[\mathbb{N}^I ]$.
So we get a canonical homomorphism of graded $A$-algebras
\[
\bigoplus _{ \nu \in \mathbb{N}^I} J^\nu \to R[\mathbb{N}^I ].
\]
By Example~\ref{triv} we have $\Proj^{\mathbb{Z}^I}( R[\mathbb{N}^I ]) = \Spec(R)$, so that
by Proposition~\ref{foncto} this morphism provides a morphism of $X$-schemes $\phi$ from an open subscheme of $\Spec(R)$ to $\Proj ( \bigoplus_{ \nu \in \mathbb{N}_I} J^\nu)$. We will now prove that this open subscheme is the whole of $\Spec(R)$, which will conclude the proof of existence. In fact, $\Proj^{\mathbb{Z}^I}( R[\mathbb{N}^I ])$ coincides with its open subscheme associated with the relevant element corresponding to the element in $\mathbb{N}^I$ all of whose components are $1$, which corresponds to $\prod_{i \in I} \gamma_i \in \bigoplus_{\nu \in \mathbb{N}^I }(f(J)R)^\nu$. Now for any $i$ there exist finite collections $(\delta_{i,j} : j \in A_i)$ of elements of $J_i$ and $(r_{i,j} : j \in A_i)$ of elements of $R$ such that $\gamma_i = \sum_{j \in A_i} r_{i,j} f(\delta_{i,j})$. Then we have
\[
 \prod_{i \in I} \gamma_i = \sum_{(j_i) \in \prod_i A_i} \prod_i r_{i,j_i} \cdot \prod_i f(\delta_{i,j_i}).
\]
In view of the construction of~\S\ref{ss:open-quasi-rel} and Lemma~\ref{lem:prod-relevant}\eqref{it:prod-relevant-0}, this implies that $\Proj^{\mathbb{Z}^I}( R[\mathbb{N}^I ])$ is covered by the open subschemes $D_\dag(\prod_i f(\delta_{i,j_i}))$. Now by construction $\phi$ is defined on each of these subsets, which finishes the justification.

Now we prove the existence of a morphism of $X$-schemes $T \to X'$ in general. For this we choose an affine open covering $X = \bigcup_{\alpha \in \mathscr{C}} X_\alpha$. Then for any $\alpha \in \mathscr{C}$ and $i \in I$, $(T \times_X X_\alpha) \times_{X_\alpha} (Y_i \times_X X_\alpha)$ is an effective Cartier divisor on $T \times_X X_\alpha$. Hence there exists an affine open covering $T \times_X X_\alpha = \bigcup_{\beta \in B_\alpha} T_{\alpha,\beta}$ such that each $T_{\alpha,\beta} \times_{X_\alpha} (Y_i \times_X X_\alpha)$ is the closed subscheme defined by a non-zero-divisor in $\mathcal{O}(T_{\alpha,\beta})$, see e.g.~\cite[\href{https://stacks.math.columbia.edu/tag/01WS}{Tag 01WS}]{stacks-project}. By the case treated above, for any $\alpha,\beta$ there exists a morphism of $X_\alpha$-schemes $T_{\alpha,\beta} \to X' \times_X X_\alpha$. By uniqueness, these morphisms coincide on the intersections of the $T_{\alpha,\beta}$'s, hence they glue to a morphism of $X$-schemes $T \to X'$, as desired.

What we have now proved and Proposition~\ref{complementary} imply that $X'$ satisfies the defining property of blowups (see Definition~\ref{blowdef}), so that the proof is complete.
\xpf

\rema 
In other words, Proposition \ref{exiblo} says that ${\Bl} _{\{Y_i : i \in I\}} X$ is a final object in the category $\Sch^{\{Y_i : i \in I \}\text{-reg}} _X$ defined in \cite[Definition~3.4]{Ma24}.
\xrema 

We finish this subsection with some applications of Proposition~\ref{exiblo}. For the first one, recall the ``$+$'' operation on closed subschemes of $X$, see~\cite[Notation~3.1]{Ma24}.

\begin{cor}
We have a canonical isomorphism of $X$-schemes ${\Bl} _{\{Y_i : i \in I\}} X \cong {\Bl} _{\sum _{i \in I} Y_i } X$.
\end{cor}

\pf
It is enough to prove that ${\Bl} _{\sum _{i \in I} Y_i } X$ is also a blowup of $X$ with center $\{Y_i : i \in I\}$. This is immediate using that ${\Bl}_{\sum _{i \in I} Y_i } X$ is a blowup of $X$ with center $\{ \sum_{i \in I} Y_i \} $ and the fact that effective Cartier divisors form a face of the monoid of closed subschemes (see~\cite[Notation~3.1]{Ma24}).
%(for the law $+$, cf. \cite[§3]{Ma24}).
\xpf

\begin{cor}
Let $K \subset I $ be a subset and set $J = I \smallsetminus K$. We have a canonical isomorphism of $X$-schemes
\[ {\Bl}_{\{Y_i : i \in I\}} X \cong {\Bl}_{\{Y_j \times _{X}  {\Bl}_{\{Y_k : k \in K\}} X : j \in J\}} {\Bl}_{\{Y_k : k \in K\}} X .
\]
\end{cor}

\pf
It is enough to prove that ${\Bl} _{\{Y_j \times _{X}  {\Bl} _{\{Y_k : k \in K\}} X : j \in J\}} {\Bl} _{\{Y_k : k \in K\}} X$ is also a blowup of $X$ with center $\{Y_i : i \in I\}$.  Using Proposition \ref{defined}, we see that  ${\Bl} _{\{Y_j \times _{X}  {\Bl} _{\{Y_k : k \in K\}} X : j \in J\}} {\Bl} _{\{Y_k : k \in K\}} X$ belongs to $\Sch _X ^{\{Y_i :i \in I \}\text{-reg}} $.
Let $T \to X$ be in $\Sch _X ^{\{Y_i :i \in I \}\text{-reg}} $. Applying twice Proposition \ref{exiblo}, we get a unique $X$-morphism $T \to {\Bl} _{\{Y_j \times _{X}  {\Bl} _{\{Y_k : k \in K\}} X : j \in J\}} {\Bl} _{\{Y_k : k \in K\}} X$. This finishes the proof.
\xpf 

%%%%%%%%%%%%%%%%%%%%%%%%%%%%%%%%%%%%%%%%
\section{Examples arising from reductive algebraic groups}
\label{s:examplereductive}
%%%%%%%%%%%%%%%%%%%%%%%%%%%%%%%%%%%%%%%%

In this section we show that a number of varieties of interest in Geometric Representation Theory can be described as Proj schemes associated with certain graded rings in a canonical way (which, in particular, does not require a specific choice of a dominant character). We see these examples as a motivation to study the general theory considered in the rest of the paper, but none of the considerations of the present sections will be used in proofs of results of other sections. (In Section \ref{s:qcohproj}, these results will only be used as illustrations of certain general statements.)

%-----------------------------------------------
\subsection{Flag varieties}
\label{ss:flag}
%-----------------------------------------------

In this section we fix an algebraically closed field $\bk$,\footnote{It is likely that most of the results in this section can be generalized to reductive group schemes over more general base schemes, but we will not go in this direction.} and a connected reductive algebraic group $\bG$ over $\bk$. We fix a Borel subgroup $\bB \subset \bG$ and a maximal torus $\bT \subset \bB$, and denote by $\bU$ the unipotent radical of $\bB$. We denote by $\bbX$ the character lattice of $\bT$ (a free $\Z$-module of finite rank), and by $\fR \subset \bbX$ the root system of $(\bG,\bT)$, i.e.~the subset of nonzero $\bT$-weights in the Lie algebra of $\bG$. We will denote by $\fR_+ \subset \fR$ the subset consisting of the opposites of the $\bT$-weights in the Lie algebra of $\bU$; this forms a positive system. We will also consider the cocharacter lattice $\bbX^\vee$ of $\bT$, which identifies with the $\Z$-dual of $\bbX$, and the system of coroots $\fR^\vee \subset \bbX^\vee$. There exists a canonical bijection $\fR \simto \fR^\vee$, which we will as usual denote $\alpha \mapsto \alpha^\vee$. Our choice of positive system $\fR_+$ lets us define a submonoid $\bbX_+ \subset \bbX$ of dominant weights as
\[
\bbX_+ = \{\lambda \in \bbX \mid \forall \alpha \in \fR_+, \, \langle \lambda, \alpha^\vee \rangle \geq 0\}.
\]
We will also consider the submonoid of strictly dominant weights, defined by
\[
\bbX_{++} = \{\lambda \in \bbX \mid \forall \alpha \in \fR_+, \, \langle \lambda, \alpha^\vee \rangle > 0\}.
\]

Consider the flag variety $\bG/\bB$, a smooth projective scheme over $\bk$.
%\footnote{Arnaud: On pourrait donner une référence précise utilisant "geometric quotient", [AR]? Wang? }. 
To any $\lambda \in \bbX$ is associated in a natural way an invertible $\calO_{\bG/\bB}(\lambda)$ on $\bG/\bB$; see e.g.~\cite[\S II.4.1]{jantzen}. (Our $\calO_{\bG/\bB}(\lambda)$ corresponds to $\mathcal{L}_{\bG/\bB}(\bk_{\bB}(\lambda))$ in the notation of~\cite{jantzen}, where $\bk_{\bB}(\lambda)$ is the $1$-dimensional $\bB$-module determined by $\lambda$.) It is a standard fact that the space
$\Gamma(\bG/\bB, \calO_{\bG/\bB}(\lambda))$
vanishes unless $\lambda \in \bbX^+$, and is finite-dimensional in this case; see e.g.~\cite[\S II.2.1 and Proposition~II.2.6]{jantzen}. For $\lambda, \mu \in \bbX$ there exists a canonical isomorphism
\begin{equation}
\label{eqn:tensor-line-bundles-flag}
\calO_{\bG/\bB}(\lambda) \otimes_{\calO_{\bG/\bB}} \calO_{\bG/\bB}(\mu) \cong \calO_{\bG/\bB}(\lambda+\mu),
\end{equation}
giving rise to a canonical morphism
\begin{equation}
\label{eqn:mult-line-bundles-flag}
\Gamma(\bG/\bB, \calO_{\bG/\bB}(\lambda)) \otimes_{\bk} \Gamma(\bG/\bB, \calO_{\bG/\bB}(\mu)) \to \Gamma(\bG/\bB, \calO_{\bG/\bB}(\lambda+\mu)),
\end{equation}
which is known to be surjective. (In case $\mathrm{char}(\bk)=0$ this fact easily follows from the simplicity of $\Gamma(\bG/\bB, \calO_{\bG/\bB}(\lambda+\mu))$ as a $\bG$-module, see e.g.~\cite[Corollary~II.5.6]{jantzen}; the case $\mathrm{char}(\bk)>0$ is more delicate, and treated e.g.~in~\cite[Theorem~3.1.2]{bk}.)

We consider the $\bbX$-graded ring
\[
A := \bigoplus_{\lambda \in \bbX_+} \Gamma(\bG/\bB, \calO_{\bG/\bB}(\lambda)),
\]
where multiplication is provided by the morphisms~\eqref{eqn:mult-line-bundles-flag}.
Its spectrum $\mathbf{X} := \Spec(A)$ is the ``affine completion'' of the space $\bG/\bU$; it is a classical object, whose main properties are summarized in~\cite[\S 6.2.1]{ar}\footnote{In this section of~\cite{ar} it is assumed that the base field has characteristic $0$. This assumption is however not used for the statements we use below.} (where it is denoted $\mathcal{X}$). 

\begin{lem}
\label{lem:A-domain}
The ring $A$ is a domain.
\end{lem}

\begin{proof}
Consider the quotient $\bG/\bU$.
Decomposing the space $\calO(\bG/\bU)$ according to the weight spaces of the natural action of $\bT$ and looking at the definitions we obtain an identification
\begin{equation}
\label{eqn:A-G/U}
 \calO(\bG/\bU) = \bigoplus_{\lambda \in \mathbb{X}} \Gamma(\bG/\bB, \calO_{\bG/\bB}(\lambda)) = A.
\end{equation}
This implies the claim, since $\bG/\bU$ is a smooth variety.
\end{proof}

The main result of this subsection is the following statement.

\begin{propo}
\label{prop:flag-var}
There exists a canonical isomorphism
\[
\bG/\bB \simto \Proj^{\bbX}(A).
\]
\end{propo}

\begin{rmk}
\begin{enumerate}
 \item 
  A similar description of $\bG/\bB$ as the Proj-scheme associated with an $\N$-graded ring, which requires a choice of a strictly dominant weight, is classical, although we were not able to find an explicit mention in the published literature. For an exposition of the proof, see the notes~\cite{wang}.
  \item
  The algebra $A$ has a canonical action of the group $\bG$. As explained in Remark~\ref{rmk:action-gp-scheme}, this induces a canonical action of $\bG$ on $\Proj^{\bbX}(A)$. It will be clear from the proofs below that this action corresponds to the obvious action on $\bG/\bB$ via the identification of Proposition~\ref{prop:flag-var}.
\end{enumerate}
\end{rmk}

We will give two independent proofs of Proposition~\ref{prop:flag-var}. The first one uses the description of the Proj-scheme as a geometric quotient under the action of a diagonalizable group scheme recalled in~\S\ref{ss:first-prop}. The second one uses results of~\cite{BS03} on ample families of invertible sheaves.

%the standard closed immersion of the flag variety in a product of projective spaces of fundamental Weyl modules.

The first proof will require a preliminary lemma.
Here we consider the ideal $A_\dag$ with respect to the $\bbX$-grading, as defined in Definition~\ref{def:rele}.

\begin{lem}
\label{lem:relevant-flag}
We have
\[
A_\dagger \subset \bigoplus_{\lambda \in \bbX_{++}} \Gamma(\bG/\bB, \calO_{\bG/\bB}(\lambda)).
\]
Moreover, if $\bG$ has a simply connected derived subgroup this inclusion is an equality.
\end{lem}

\begin{proof}
It is clear that $\bigoplus_{\lambda \in \bbX_{++}} \Gamma(\bG/\bB, \calO_{\bG/\bB}(\lambda))$ is an ideal in $A$; so to prove the inclusion it suffices to prove that the degree of any nonzero relevant element in $A$ belongs to $\bbX_{++}$.
%If $f \in A$ is a nonzero relevant element, 
In fact if $f$ is nonzero and homogeneous and if there exists $\alpha \in \fR_+$ such that $\langle \deg(f), \alpha^\vee \rangle = 0$, then any homogeneous divisor $g$ of a power of $f$ (which is necessarily nonzero by Lemma~\ref{lem:A-domain}) will also satisfy $\langle \deg(g), \alpha^\vee \rangle = 0$, so that $f$ is not relevant. 

% As a consequence, we have
% \[
% A_\dagger \subset \bigoplus_{\lambda \in \bbX_{++}} \Gamma(\bG/\bB, \calO_{\bG/\bB}(\lambda)).
% \]

Now, assume that $\bG$ has simply connected derived subgroup, and let $\fR_{\mathrm{s}}$ be the set of simple roots. Our assumption ensures that for any $\alpha \in \fR_{\mathrm{s}}$ there exists a weight $\varpi_\alpha \in \bbX$ such that
\[
\langle \varpi_\alpha, \beta^\vee \rangle = \delta_{\alpha,\beta}
\]
for $\beta \in \fR_{\mathrm{s}}$. (In case $\bG$ is semisimple these weights are uniquely determined, and called the fundamental weights. If $\bG$ is not semisimple, they are not unique.) Setting
\[
\bbX_0 = \{ \lambda \in \bbX \mid \forall \alpha \in \fR_{\mathrm{s}}, \, \langle \lambda, \alpha^\vee \rangle = 0 \},
\]
then we have
\begin{equation}
\label{eqn:X-sc}
\bbX = \bbX_0 \oplus \left( \bigoplus_{\alpha \in \fR_{\mathrm{s}}} \Z \varpi_\alpha \right), \quad \mathbb{X}_+ = \left( \bigoplus_{\alpha \in \fR_{\mathrm{s}}} \Z_{\geq 0} \cdot \varpi_{\alpha} \right) \oplus \mathbb{X}_0.
\end{equation}
Any $\lambda \in \bbX_0$ extends to a character of $\bG$, so that the space $\Gamma(\bG/\bB, \calO_{\bG/\bB}(\lambda))$ is $1$-dimensional (by the tensor identity~\cite[Proposition~I.3.6]{jantzen}), and its nonzero elements are invertible in $A$.

By the surjectivity of~\eqref{eqn:mult-line-bundles-flag} when $\lambda,\mu$ are dominant, the ideal
\[
\bigoplus_{\lambda \in \bbX_{++}} \Gamma(\bG/\bB, \calO_{\bG/\bB}(\lambda)) \subset A
\]
is generated by elements which are products $\prod_{\alpha \in \fR_{\mathrm{s}}} f_\alpha$ with $\deg(f_\alpha)=\varpi_\alpha$. Such products are relevant since they admit a divisor of degree $\varpi_\alpha$ for each $\alpha \in \fR_{\mathrm{s}}$, and also of degree $\lambda$ for any $\lambda \in \bbX_0$. This proves the inclusion
\[
\bigoplus_{\lambda \in \bbX_{++}} \Gamma(\bG/\bB, \calO_{\bG/\bB}(\lambda)) \subset A_\dagger
\]
and finishes the proof.
%of elements in the subspaces
%\[
%\Gamma(\bG/\bB, \calO_{\bG/\bB}(\varpi_\alpha)) \quad \text{for $\alpha \in \fR_{\mathrm{s}}$ and} 
%\]
\end{proof}

Now we can give the first proof of Proposition~\ref{prop:flag-var}.

\begin{proof}[First proof of Proposition~\ref{prop:flag-var}]
By~\cite[Proposition~4.2]{BS03}, there exists a canonical rational map from $\bG/\bB$ to $\Proj^{\bbX}(A)$; what we will show is that this rational morphism is defined everywhere, and an isomorphism.

It is a standard fact that there exists a connected reductive algebraic group $\widetilde{\bG}$ over $\bk$ with simply connected derived subgroup and a finite central isogeny $\widetilde{\bG} \to \bG$. (This follows e.g.~from the considerations in~\cite[\S I.1.18]{jantzen}.) The preimages $\widetilde{\bB}$ and $\widetilde{\bT}$ of $\bB$ and $\bT$ are a Borel subgroup and a maximal torus in $\widetilde{\bG}$ respectively. Let also $\widetilde{\bbX}$ and $\widetilde{A}$ be the counterparts of $\bbX$ and $A$ for $\widetilde{\bG}$, $\widetilde{\bB}$ and $\widetilde{\bT}$. We have an isomorphism
\begin{equation}
\label{eqn:flag-variety-identification}
\widetilde{\bG} / \widetilde{\bB} \simto \bG/\bB,
\end{equation}
and composition with the surjection $\widetilde{\bT} \to \bT$ induces an embedding $\bbX \hookrightarrow \widetilde{\bbX}$ whose cokernel is finite. For $\lambda \in \bbX$, under the identification~\eqref{eqn:flag-variety-identification} the invertible sheaf on $\bG/\bB$ associated with $\lambda$ identifies with the invertible sheaf on $\widetilde{\bG} / \widetilde{\bB}$ associated with its image in $\widetilde{\bbX}$; hence $A$ identifies with the subring of $\widetilde{A}$ given by the sum of the components whose label is in $\bbX$. In view of Proposition~\ref{thin}, we deduce an isomorphism
\[
\Proj^{\bbX}(A) \cong \Proj^{\widetilde{\bbX}}(\widetilde{A}).
\]
Comparing with~\eqref{eqn:flag-variety-identification}, this reduces the proof to the case $\bG$ has simply connected derived subgroup.

In this case, we have described the ideal $A_\dagger \subset A$ in Lemma~\ref{lem:relevant-flag}. 
This ideal coincides with that considered in~\cite[Remark~6.2.3]{ar}; 
in view of this statement we deduce an identification
\begin{equation}
\label{eqn:G/U}
\Spec(A) \smallsetminus V(A_\dagger) = \bG/\bU.
\end{equation}
Using the comments preceding Remark~\ref{rmk:action-gp-scheme}, we see that $\Proj^{\bbX}(A)$ identifies with the geometric quotient of the left-hand side by the action of $\bT$. (This quotient is unique by~\cite[Proposition~0.1]{MFK93}.) It is clear that the projection $\bG/\bU \to \bG/\bB$ identifies $\bG/\bB$ with the geometric quotient of $\bG/\bU$ by the action of $\bT$, which implies our identification.
%Passing to quotients under the actions of $\bT$ we deduce the desired isomorphism.\footnote{Arnaud: On pourrait faire référence à la Description \ref{geodescri}.}
\end{proof}

\begin{rmk}
 %The pushforward of the structure sheaf of $\bG/\bU$ under the canonical $\bT$-torsor $\bG/\bU \to \bG/\bB$ is $\bigoplus_{}
% Decomposing the space $\calO(\bG/\bU)$ according to the weight spaces of the natural action of $\bT$ and looking at the definitions we obtain an identification
% \[
%  \calO(\bG/\bU) = \bigoplus_{\lambda \in \mathbb{X}} \Gamma(\bG/\bB, \calO_{\bG/\bB}(\lambda)) = A.
% \]
% Hence we have $\Spec(A)=(\bG/\bU)_\aff$.
In view of the isomorphism~\eqref{eqn:A-G/U}, the above proof amounts to noticing that, in case $\bG$ has simply connected derived subgroup, the natural morphism
$\bG/\bU \to (\bG/\bU)_\aff$
is an open immersion, with image $\Spec(A) \smallsetminus V(A_\dag)$.
 \end{rmk}

\begin{proof}[Second proof of Proposition~\ref{prop:flag-var}]
 As in the first proof, one can assume that $\bG$ has simply connected derived subgroup. 
 In fact, one can even reduce the proof to the case $\bG$ is \emph{semisimple} and simply connected. Indeed, $\bG$ and its derived subgroup have the same flag variety. On the other hand, let $\bbX_0$, $\fR_{\mathrm{s}}$ and $(\varpi_\alpha : \alpha \in \fR_{\mathrm{s}})$ are as in the proof of Lemma~\ref{lem:relevant-flag}, 
 and recall the descriptions~\eqref{eqn:X-sc} of $\mathbb{X}$ and $\mathbb{X}_+$.
%  we have
%   \[
%   \mathbb{X}_+ = \left( \bigoplus_{\alpha \in \fR_{\mathrm{s}}} \Z_{\geq 0} \cdot \varpi_{\alpha} \right) \oplus \mathbb{X}_0.
%  \]
%  Moreover, 
 Denoting by $A_1$, resp.~$A_2$, the sum of the components of $A$ whose degrees belong to $\bigoplus_{\alpha \in \fR_{\mathrm{s}}} \Z_{\geq 0} \cdot \varpi_{\alpha}$, resp.~to $\mathbb{X}_0$, multiplication induces an isomorphism
 \[
  A_1 \otimes A_2 \simto A.
 \]
By Proposition~\ref{produitvstensoriel} we deduce an isomorphism
\[
 \Proj^{\mathbb{X}}(A) \cong \Proj^{\bigoplus_{\alpha \in \fR_{\mathrm{s}}} \Z \cdot \varpi_{\alpha}}(A_1) \times \Proj^{\mathbb{X}_0}(A_2).
\]
Now $A_2$ identifies with $\bk[\mathbb{X}_0]$, so that by Example~\ref{triv} we have $\Proj^{\mathbb{X}_0}(A_2)=\Spec(\bk)$. Since $A_1$ identifies with the version of $A$ associated with the derived subgroup of $\bG$, this finishes the reduction of the proof to the semisimple case.
 
 From now on, we therefore assume that $\bG$ is semisimple and simply connected. In this case, we will show that the claim follows from~\cite[Corollary~4.6]{BS03}. First, let us prove that $\Proj^{\mathbb{X}}(A)$ is separated. 
%  Let $\bbX_0$ be as above, and choose a basis $(\lambda_i : i \in I)$ of this lattice. Then
%  $(\varpi_\alpha : \alpha \in \fR_{\mathrm{s}}) \cup (\lambda_i : i \in I)$ is a basis of $\mathbb{X}$, and we have
%  \[
%   \mathbb{X}_+ = \left( \bigoplus_{\alpha \in \fR_{\mathrm{s}}} \Z_{\geq 0} \cdot \varpi_{\alpha} \right) \oplus \left( \bigoplus_{i \in I} \Z \lambda_i \right).
%  \]
% %The space $A_{\varpi_\alpha} = \Gamma(\bG/\bB, \calO_{\bG/\bB}(\varpi_\alpha))$ is the fundamental representation associated with $\alpha$, and 
% For any $i \in I$, as noted in the proof of Lemma~\ref{lem:relevant-flag} the vector space $A_{\lambda_i}$ has dimension $1$, and
As explained above we have $\mathbb{X}_+ = \bigoplus_{\alpha \in \fR_{\mathrm{s}}} \Z_{\geq 0} \cdot \varpi_{\alpha}$, and by the surjectivity of~\eqref{eqn:mult-line-bundles-flag} the natural morphism of $\mathbb{X}$-graded rings
\begin{equation}
\label{eqn:surjection-Sym-A}
 \left( \bigotimes_{\alpha \in \fR_{\mathrm{s}}} \mathrm{S}(A_{\varpi_\alpha}) \right) 
 %\otimes \left( \bigotimes_{(i,n) \in I \times \Z} (A_{\lambda_i})^{\otimes n} \right) 
 \longto A
\end{equation}
is surjective,
where we denote by $\mathrm{S}(V)$ the symmetric algebra of a vector space $V$.
By Lemma~\ref{lem:closed} we deduce a closed immersion
\[
 \Proj^{\mathbb{X}}(A) \hookrightarrow \Proj^{\mathbb{X}} \left( 
 %\left( 
 \bigotimes_{\alpha \in \fR_{\mathrm{s}}} \mathrm{S}(A_{\varpi_\alpha}) \right), 
 %\otimes \left( \bigotimes_{(i,n) \in I \times \Z} (A_{\lambda_i})^{\otimes n} \right) \right),
\]
and by Proposition~\ref{produitvstensoriel} the right-hand side identifies with
\[
% \left( 
 \prod_{\alpha \in \fR_{\mathrm{s}}} \Proj^{\Z \cdot \varpi_\alpha} \bigl( \mathrm{S}(A_{\varpi_\alpha} ) \bigr). 
 %\right) \times \left( \prod_{i \in I} \Proj^{\Z \lambda_i} \left( \bigotimes_{n \in \Z} (A_{\lambda_i})^{\otimes n} \right) \right).
\]
By the theory of $\N$-graded Proj schemes, for $\alpha \in \fR_{\mathrm{s}}$ the scheme $\Proj^{\Z \cdot \varpi_\alpha} \bigl( \mathrm{S}(A_{\varpi_\alpha} ) \bigr)$ identifies with the projective space of the dual vector space $(A_{\varpi_\alpha})^*$ (which is the Weyl module of highest weight $-w_\circ(\varpi_\alpha)$, where $w_\circ$ is the longest element in the Weyl group). 
%On the other hand, for $i \in I$ we have a ring isomorphism $\bigotimes_{n \in \Z} (A_{\lambda_i})^{\otimes n} \cong \bk[x^{\pm 1}]$, so that by Example~\ref{triv} we have $\Proj^{\Z \lambda_i} \bigl( \bigotimes_{n \in \Z} (A_{\lambda_i})^{\otimes n} \bigr) = \Spec(\bk)$. 
Hence we obtain a closed immersion
\begin{equation}
\label{eqn:embedding-flag-projective-spaces}
 \Proj^{\mathbb{X}}(A) \hookrightarrow \prod_{\alpha \in \fR_{\mathrm{s}}} \mathbb{P} \bigl( (A_{\varpi_\alpha})^* \bigr).
\end{equation}
Since projective spaces (and, more generally, Proj schemes associated with $\N$-graded rings, see Remark~\ref{rmk:def-Proj}) are separated and closed immersions are separated (see~\cite[\href{https://stacks.math.columbia.edu/tag/01QU}{Tag 01QU}]{stacks-project}), this indeed shows that $\Proj^{\mathbb{X}}(A)$ is separated.

Let us note also that we have $(\bG/\bB)_\aff = \Spec(\bk)$, so that the morphism $\bG/\bB \to (\bG/\bB)_\aff$ is projective, hence proper.

Finally, choose a numbering $\alpha_1, \dots, \alpha_r$ of the elements in $\mathfrak{R}_{\mathrm{s}}$.
It is a standard fact (see e.g.~\cite[Proposition~II.4.4]{jantzen}) that $\calO_{\bG/\bB}(\lambda)$ is ample for any $\lambda \in \mathbb{X}_{++}$; in particular, $\calO_{\bG/\bB}(\varpi_{\alpha_1} + \cdots + \varpi_{\alpha_r})$ is ample. Comparing the characterization of ample invertible sheaves in~\cite[\href{https://stacks.math.columbia.edu/tag/01PS}{Tag 01PS}]{stacks-project} with~\cite[Proposition~1.1]{BS03} we deduce that the collection
  \[
  \calO_{\bG/\bB}(\varpi_{\alpha_1}), \dots, \calO_{\bG/\bB}(\varpi_{\alpha_r})
  \]
  is an ample family in the sense of~\cite{BS03}.
  
We have now checked that the conditions in~\cite[Corollary~4.6]{BS03} are satisfied, and this statement allows us to conclude.
\end{proof}

\begin{rmk}
\label{rmk:flag-sc}
\begin{enumerate}
 \item
 In the proof above, the fact that $\Proj^{\bbX}(A)$ is separated can also be deduced from Proposition~\ref{prop:separated-family}.
 \item
 In the case when $\bG$ is semisimple and simply connected, the closed immersion~\eqref{eqn:embedding-flag-projective-spaces} recovers the standard closed immersion of the flag variety in a product of projective spaces of fundamental Weyl modules.
\end{enumerate}
\end{rmk}

\begin{rmk}
\label{rmk:flag-covering}
One can construct an explicit open covering of $\Proj^{\bbX}(A)$ as follows. 
%
%First, assume that $\bG$ has simply connected derived subgroup and,
%%choose a basis $(\lambda_i : i \in I)$ of $\bbX_0$.
%%as in Remark~\ref{rmk:flag-sc}.
%%consider a family $(\varpi_\alpha : \alpha \in \fR_{\mathrm{s}})$ as above. 
%for any $\alpha \in \fR_{\mathrm{s}}$, fix a basis $(f^\alpha_i : i \in I_\alpha)$ of $A_{\varpi_\alpha}$. Let us denote by $E$ the set of sections of the obvious map
%\[
% \bigsqcup_{\alpha \in \fR_{\mathrm{s}}} I_\alpha \to \fR_{\mathrm{s}},
%\]
%and for $\sigma \in E$ set $f_\sigma = \prod_{\alpha \in \fR_{\mathrm{s}}} f^\alpha_{\sigma(\alpha)}$.
%Then each $f_\sigma$ is relevant, and by surjectivity of~\eqref{eqn:surjection-Sym-A} and Corollary~\ref{cor:covering-Proj} we have
%\[
% \Proj^{\bbX}(A) = \bigcup_{\sigma \in E} D_\dag ( f_\sigma ).
%\]
%%In fact, it is clear that we have such an open covering of the product in the right-hand side of~\eqref{eqn:embedding-flag-projective-spaces}, and then the claim follows by the proof of Lemma~\ref{lem:closed}. 
%Note that for any $\sigma \in E$ we have $\bbX[\underline{(f_\sigma)^\N}]=\bbX$; in other words each $f_\sigma$ is maximally relevant in the sense considered in~\S\ref{ss:max-relevant} below.
%
Choose a finite central isogeny $\widetilde{\bG} \to \bG$ where $\widetilde{\bG}$ has simply connected derived subgroup, see the first proof of Proposition~\ref{prop:flag-var}. Let $\widetilde{\bbX}$ and $\widetilde{A}$ be as in this proof, and recall that $A$ identifies with the sum of the graded components of $A$ labelled by elements in $\bbX \subset \widetilde{\bbX}$.
%, and that $\Proj^{\bbX}(A) = \Proj^{\widetilde{\bbX}}(\widetilde{A})$. 
Choose elements $(\varpi_\alpha : \alpha \in \fR_{\mathrm{s}})$ as above for the group $\widetilde{\bG}$ and, for any $\alpha \in \fR_{\mathrm{s}}$, choose a 
basis $(f^\alpha_i : i \in I_\alpha)$ of $\widetilde{A}_{\varpi_\alpha}$.
%$(\lambda_i : i \in I)$ 
%as above for the group $\widetilde{\bG}$, and 
Let $E$ be the set of sections of the obvious map
\[
 \bigsqcup_{\alpha \in \fR_{\mathrm{s}}} I_\alpha \longto \fR_{\mathrm{s}};
\]
for $\sigma \in E$ we set $f_\sigma = \prod_{\alpha \in \fR_{\mathrm{s}}} f^\alpha_{\sigma(\alpha)}$.
%$E$ be the set defined as above. 
Since $\bbX \subset \widetilde{\bbX}$ has finite index, for any $\sigma \in E$ there exists $n_\sigma \in \Z_{\geq 1}$ such that $g_\sigma := (f_\sigma)^{n_\sigma} \in A$. Then we claim that 
\begin{equation*}
%\label{eqn:Adag-flag}
A_\dag \subset \Rad \left( \sum_{ \sigma \in E} A \cdot g_\sigma \right),
\end{equation*}
so that by Corollary~\ref{cor:covering-Proj} we
have
\[
 \Proj^{\bbX}(A) = \bigcup_{\sigma \in E} D_\dag ( g_\sigma ).
\]
%Here again, for $\sigma \in E$ we have $\bbX[\underline{(g_\sigma)^\N}]=\bbX$, so that each $g_\sigma$ is maximally relevant.

In fact, by Lemma~\ref{lem:relevant-flag}, to prove this inclusion it suffices to prove that 
%for any $\lambda \in \bbX_{++}$ we have 
\begin{equation}
\label{eqn:Alambda-flag}
A_\lambda \subset \Rad \left( \sum_{ \sigma \in E} A \cdot g_\sigma \right) \quad \text{for any $\lambda \in \bbX_{++}$}.
\end{equation}
Fix such $\lambda$.
If $a \in A_\lambda$, there exist elements $(b_\sigma : \sigma \in E)$ in $\widetilde{A}_{\lambda-\sum_\alpha \varpi_\alpha}$ such that $a = \sum_\sigma b_\sigma \cdot f_\sigma$. Then if $N > |E| \cdot (\max_{\sigma} n_\sigma -1)$, we have $a^N = \sum_\sigma c_\sigma g_\sigma$ for some homogeneous elements $c_\sigma \in \widetilde{A}$, where $\deg(c_\sigma) + \deg(g_\sigma) = N\lambda$ for any $\sigma$ such that $c_\sigma g_\sigma \neq 0$. Here we have $\deg(c_\sigma) \in \bbX$, hence $c_\sigma \in A$, and the desired claim follows.
\end{rmk}

%-----------------------------------------------
\subsection{Some vector bundles over the flag variety}
\label{ss:vector-bundles}
%-----------------------------------------------

We continue with the setting of~\S\ref{ss:flag}, and consider a finite-dimensional 
%$\bB$-module $V$.
$\bG$-module $\widetilde{V}$ and a $\bB$-stable subspace $V \subset \widetilde{V}$. 
We can then consider the induced scheme
\[
 \bG \times^{\bB} V,
\]
i.e.~the quotient of the product $\bG \times V$ by the (free) action of $\bB$ defined by $b \cdot (g,x) = (gb^{-1}, b \cdot x)$. (This construction is a special case of that considered in~\cite[\S I.5.14]{jantzen}.) It is a vector bundle over $\bG/\bB$ (in particular, a smooth variety). For $\lambda \in \bbX$ we will denote by $\calO_{\bG \times^{\bB} V}(\lambda)$ the pullback to $\bG \times^{\bB} V$ of the invertible sheaf $\calO_{\bG/\bB}(\lambda)$. By~\eqref{eqn:tensor-line-bundles-flag}, for $\lambda,\mu \in \bbX$ we have a canonical isomorphism
\[
\calO_{\bG \times^{\bB} V}(\lambda) \otimes_{\calO_{\bG \times^{\bB} V}} \calO_{\bG \times^{\bB} V}(\mu) \cong \calO_{\bG \times^{\bB} V}(\lambda+\mu),
\]
which provides a structure of a $\bbX$-graded ring on
\[
A_V := \bigoplus_{\lambda \in \bbX_+} \Gamma(\bG \times^{\bB} V, \calO_{\bG \times^{\bB} V}(\lambda)).
\]

\begin{rmk}
 One can also consider the quotient $\bG \times^{\bU} V$ of the restriction of the action above to $\bU$, which is a (Zariski locally trivial) $\bT$-torsor over $\bG \times^{\bB} V$. As in the proof of Lemma~\ref{lem:A-domain} one can see that $A_V$ identifies with a subring of $\mathcal{O}(\bG \times^{\bU} V)$, hence is a domain.
\end{rmk}

%We will assume that there exists an affine $\bk$-scheme $\widetilde{V}$ endowed with an action of $\bG$ and a $\bB$-equivariant closed immersion $V \to \widetilde{V}$.
The $\bB$-equivariant embedding $V \subset \widetilde{V}$ 
%This immersion 
induces a closed immersion
\begin{equation}
\label{eqn:embedding-vb}
 \bG \times^{\bB} V \hookrightarrow \bG \times^{\bB} \widetilde{V},
\end{equation}
and the right-hand side identifies with $\bG/\bB \times \widetilde{V}$ (via the map $[g,x] \mapsto (g\bB, g \cdot x)$ for $g \in \bG$ and $x \in \widetilde{V}$). For any $\lambda \in \mathbb{X}_+$ we deduce a canonical morphism
\begin{equation}
\label{eqn:morph-restriction-vb}
 \calO(\widetilde{V}) \otimes A_\lambda \longto (A_V)_\lambda.
\end{equation}
Since $\bG/\bB$ is projective, these considerations also show that there exists a canonical projective morphism
\begin{equation}
\label{eqn:morph-vb}
 \bG \times^{\bB} V \longto \widetilde{V}.
\end{equation}

Our goal in this subsection is to prove the following statement.

\begin{propo}
\label{prop:vb-flag-Proj}
%Assume that the morphism $\calO(\widetilde{V}) \to \calO(\bG \times^{\bB} V)$ induced by~\eqref{eqn:morph-vb} is surjective. Then
There exists a canonical isomorphism
\[
\bG \times^{\bB} V \cong \Proj^{\bbX}(A_V).
\]
\end{propo}

\begin{rmk}
% \phantomsection
% \begin{enumerate}
% \item
In case $V=\widetilde{V}=\{0\}$, Proposition~\ref{prop:vb-flag-Proj} recovers Proposition~\ref{prop:flag-var} (but our proof will use the latter statement).
%\item
Another interesting case is when $\widetilde{V}$ is the Lie algebra $\bg$ of $\bG$ and $V$ is the Lie algebra $\bu$ of the unipotent radical of $\bB$. In this case, $\bG \times^{\bB} V$ is the so-called \emph{Springer resolution}; see~\S\ref{ss:Springer} below for more details about this case.
%our assumption is satisfied at least when $\bG$ has simply connected derived subgroup and the characteristic of $\bk$ is good for $\bG$. (In fact, this follows from a form of Zariski's Main Theorem, see e.g.~\cite[\href{https://stacks.math.columbia.edu/tag/0AY8}{Tag 0AY8}]{stacks-project}, combined with the results gathered in~\cite[\S 10.1]{jantzen-nilp}.)
%\end{enumerate}
\end{rmk}

%As for Proposition~\ref{prop:flag-var}, we will give two independent proofs of this statement. Both will use 
The proof of Proposition~\ref{prop:vb-flag-Proj} will rely on the following preliminary result. 

\begin{lem}
\label{lem:restriction-vb-surjective}
 There exists $N \in \Z_{> 0}$ such that for any $\lambda \in \mathbb{X}_+$ which satisfies $\langle \lambda, \alpha^\vee \rangle \geq N$ for any simple root $\alpha$, the morphism~\eqref{eqn:morph-restriction-vb} is surjective.
\end{lem}

\begin{proof}
 Consider the natural (affine) morphism $\pi : \bG \times^{\bB} V \to \bG/\bB$. If we denote by $\mathcal{V}$ the vector bundle on $\bG/\bB$ associated with the $\bB$-module $V$ (this vector bundle is denoted $\mathcal{L}_{\bG/\bB}(V)$ in~\cite[\S I.5.8]{jantzen}), and by $\mathcal{V}^\vee$ the dual vector bundle, then $\pi_* \calO_{\bG \times^{\bB} V}$ identifies with the symmetric algebra $\mathrm{S}_{\calO_{\bG/\bB}}(\mathcal{V}^\vee)$. Similarly, the pushforward to $\bG/\bB$ of $\calO_{\bG \times^{\bB} \widetilde{V}}$ identifies with $\mathcal{O}(\widetilde{V}^*) \otimes \calO_{\bG/\bB}$, and the closed immersion~\eqref{eqn:embedding-vb} 
% $\bG \times^{\bB} V \hookrightarrow \bG \times^{\bB} \widetilde{V}$ 
 corresponds to a surjection of sheaves
 \[
  \mathcal{O}(\widetilde{V}) \otimes \calO_{\bG/\bB} \twoheadrightarrow \mathrm{S}_{\calO_{\bG/\bB}}(\mathcal{V}^\vee).
 \]
Denoting by $\mathcal{V}^\bot \subset \widetilde{V}^* \otimes \calO_{\bG/\bB}$ the orthogonal of $\mathcal{V} \subset \widetilde{V} \otimes \calO_{\bG/\bB}$ (a sub-vector bundle), this surjection can be ``extended'' to the Koszul resolution
\[
 \bigwedge \hspace{-2pt}{}^{-\bullet}(\mathcal{V}^\bot) \otimes_{\calO_{\bG/\bB}} \left( \mathcal{O}(\widetilde{V}) \otimes \calO_{\bG/\bB} \right) \to \mathrm{S}_{\calO_{\bG/\bB}}(\mathcal{V}^\vee),
\]
a quasi-isomorphism of complexes of (quasi-coherent) $\calO_{\bG/\bB}$-modules. For any $\lambda \in \mathbb{X}$, tensoring with $\calO_{\bG/\bB}(\lambda)$ we deduce a quasi-isomorphism
\begin{equation}
\label{eqn:Koszul-resolution}
 \bigwedge \hspace{-2pt}{}^{-\bullet}(\mathcal{V}^\bot) \otimes_{\calO_{\bG/\bB}} \left( \mathcal{O}(\widetilde{V}) \otimes \calO_{\bG/\bB} \right) \otimes_{\calO_{\bG/\bB}} \calO_{\bG/\bB}(\lambda) \to \pi_* \calO_{\bG \times^{\bB} V}(\lambda).
\end{equation}

Now $\mathcal{V}^\bot$ 
%is an iterated extension of the line bundles $\calO_{\bG/\bB}(\mu)$ where $\mu$ runs over the (finitely many) $T$-weights in 
is the vector bundle on $\bG/\bB$ associated with
the orthogonal $V^\bot \subset \widetilde{V}^*$ of $V \subset \widetilde{V}$. Hence, if we denote by $\Lambda \subset \mathbb{X}$ the (finite) subset consisting of the $\bT$-weights appearing in the various exterior powers of $V^\bot$, then each $\wedge^i \mathcal{V}^\bot$ is an iterated extension of invertible sheaves of the form $\calO_{\bG/\bB}(\mu)$ with $\mu \in \Lambda$.
If we choose $N$ such that $\langle \mu, \alpha^\vee \rangle \geq -N$ for any $\mu \in \Lambda$ and any simple root $\alpha$, then if $\lambda \in \bbX_+$ satisfies $\langle \lambda, \alpha^\vee \rangle \geq N$ for any simple root $\alpha$, each weight $\lambda+\mu$ with $\mu \in \Lambda$ is dominant. Since invertible sheaves on $\bG/\bB$ associated with dominant weights have no higher cohomology (this is Kempf's vanishing theorem, see~\cite[Proposition~II.4.5]{jantzen}), for such $\lambda$ we have
\[
 \mathsf{H}^n(\wedge^i \mathcal{V}^\bot \otimes_{\calO_{\bG/\bB}} \calO_{\bG/\bB}(\lambda))=0 \quad \text{for any $i$ and any $n>0$.}
\]
Breaking the resolution~\eqref{eqn:Koszul-resolution} into short exact sequences and using the long exact sequences obtained by applying the global sections functor and its derived functors, we deduce in this case a quasi-isomorphism
\[
\mathcal{O}(\widetilde{V}) \otimes \Gamma \left( \bG/\bB, \bigwedge \hspace{-2pt}{}^{-\bullet}(\mathcal{V}^\bot) \otimes_{\calO_{\bG/\bB}} \calO_{\bG/\bB}(\lambda) \right) \to \Gamma \bigl( \bG/\bB, \pi_* \calO_{\bG \times^{\bB} V}(\lambda) \bigr) = (A_V)_\lambda,
\]
proving in particular that the morphism of the lemma is surjective.
\end{proof}

\begin{proof}[Proof of Proposition~\ref{prop:vb-flag-Proj}]
 As in the second proof of Proposition~\ref{prop:flag-var}, one can (and will) assume that $\bG$ is semisimple and simply connected. Here again, to prove the proposition we will show that the conditions in~\cite[Corollary~4.6]{BS03} are satisfied.
 
 Recall the notation of Remark~\ref{rmk:flag-covering} and, for any $\sigma$, denote by $\tilde{f}_\sigma$ the image of $f_\sigma$ in $A_V$. We claim that
 \begin{equation}
 \label{eqn:covering-Proj-AV}
 \Proj^{\bbX}(A_V) = \bigcup_{\sigma \in E} D_\dag(\tilde{f}_\sigma).
 \end{equation}
 In fact, as in the proof of Lemma~\ref{lem:relevant-flag} one sees that the degree of any nonzero relevant element $f$ in $A_V$ belongs to $\bbX_{++}$. If $N$ satisfies the condition in Lemma~\ref{lem:restriction-vb-surjective}, then $f^N$ belongs to the image of~\eqref{eqn:morph-restriction-vb} for $\lambda=N\deg(f)$, hence to the ideal generated by the elements $\tilde{f}_\sigma$. This implies the claim in view of Corollary~\ref{cor:covering-Proj}.
 
 Next, we will prove that the scheme $\Proj^{\mathbb{X}}(A_V)$ is separated. In fact, combining the morphisms~\eqref{eqn:morph-restriction-vb} we obtain a canonical morphism of $\mathbb{X}$-graded rings
\begin{equation}
\label{eqn:morph-A-AV}
\calO(\widetilde{V}) \otimes A \longto A_V.
\end{equation}
 (Here, the grading on the left-hand side is such that the degree-$\lambda$ component is $\calO(\widetilde{V}) \otimes A_\lambda$.) 
% We deduce a morphism of $\mathbb{X}$-graded rings
% \begin{equation}
% \label{eqn:morph-gr-rings-vb}
% \calO(\bG \times^{\bB} V) \otimes A = \calO(\bG \times^{\bB} V) \otimes_{\calO(\widetilde{V})} \bigl( \calO(\widetilde{V}) \otimes A \bigr) \to A_V,
% \end{equation}
% where the morphism $\calO(\widetilde{V}) \to \calO(\bG \times^{\bB} V)$ is induced by~\eqref{eqn:morph-vb}.
Consider the associated rational morphism
\begin{equation}
\label{eqn:morph-Proj-AV}
 \Proj^{\mathbb{X}}(A_V) \dashrightarrow \Proj^{\mathbb{X}} (\calO(\widetilde{V}) \otimes A) = \widetilde{V} \times \Proj^{\mathbb{X}}(A)
 %\left(  \calO(\bG \times^{\bB} V) \otimes A \right)
\end{equation}
provided by Proposition~\ref{foncto}, where the identification is provided by Corollary~\ref{propaffinebasechange}. By~\eqref{eqn:covering-Proj-AV} this morphism is defined everywhere. We claim that it is a closed immersion. In fact, 
by Remarks~\ref{rmk:open-basechange} and~\ref{rmk:flag-covering} we have
\[
%\Proj^{\mathbb{X}} \left(  \calO(\bG \times^{\bB} V) \otimes A \right) 
\Proj^{\mathbb{X}} (\calO(\widetilde{V}) \otimes A)
= \bigcup_{\sigma \in E} D_\dag(1 \otimes f_\sigma).
\]
Using the fact that the property of being a closed immersion is local on the target (see~\cite[\href{https://stacks.math.columbia.edu/tag/02L6}{Tag 02L6}]{stacks-project}), we deduce that to prove the claim it suffices to prove that for any $\sigma \in E$ the projection
\[
\Proj^{\mathbb{X}}(A_V) \times_{\Proj^{\mathbb{X}} (\calO(\widetilde{V}) \otimes A)}
D_\dag(1 \otimes f_\sigma) \to D_\dag(1 \otimes f_\sigma)
\]
is a closed immersion. Fix such a $\sigma$.
By~\eqref{eqn:Proj-morph-affine}, we have
\[
\Proj^{\mathbb{X}}(A_V) \times_{\Proj^{\mathbb{X}} (\calO(\widetilde{V}) \otimes A)}
%_{\Proj^{\mathbb{X}} (  \calO(\bG \times^{\bB} V) \otimes A)} 
D_\dag(1 \otimes f_\sigma) = D_\dag(\tilde{f}_\sigma).
\]
Therefore, to conclude it suffices to prove that the morphism $\calO(\widetilde{V}) \otimes A_{(f_\sigma)} \to (A_V)_{(\tilde{f}_\sigma)}$
%$\calO(\bG \times^{\bB} V) \otimes A_{(f_\sigma)} \to (A_V)_{(\tilde{f}_\sigma)}$ 
is surjective. But, if $N$ is as above any nonzero element in $(A_V)_{(\tilde{f}_\sigma)}$ can be written in the form $\frac{a}{(\tilde{f}_\sigma)^{kN}}$ for some $k>0$, where $a \in A$ is homogeneous of degree $kN\deg(\tilde{f}_\sigma)$. Then $a$ belongs to the image of $\calO(\widetilde{V}) \otimes A$, which implies the desired surjectivity. 

% This result implies that the projection
% \[
% \Proj^{\mathbb{X}}(A_V) \times_{\Proj^{\mathbb{X}} (\calO(\widetilde{V}) \otimes A)}
% D_\dag(1 \otimes f_\sigma) \to D_\dag(1 \otimes f_\sigma)
% \]
% is a closed immersion, which finishes the proof of our claim since the property of being a closed immersion is local on the target (see~\cite[\href{https://stacks.math.columbia.edu/tag/02L6}{Tag 02L6}]{stacks-project}).
%In fact, by Proposition~\ref{thin}, for any $N \in \Z_{>0}$ we can replace $\mathbb{X}$ by $N \cdot \mathbb{X}$ and $A, A_V$ by the sum of their components whose degrees belong to $N \cdot \mathbb{X}$. Lemma~\ref{lem:restriction-vb-surjective} ensures that after an appropriate such replacement our ring morphism~\eqref{eqn:morph-gr-rings-vb} becomes surjective, and then the claim follows from Lemma~\ref{lem:closed}.

%By Corollary~\ref{propaffinebasechange} we have an identification
%\[
% \Proj^{\mathbb{X}} \left( \calO(\bG \times^{\bB} V) \otimes A \right) = \Spec(\calO(\bG \times^{\bB} V)) \times \Proj^{\mathbb{X}}(A),
%\]
%and 
By Proposition~\ref{prop:flag-var} the right-hand side in~\eqref{eqn:morph-Proj-AV}
%identifies with $\Spec(\calO(\bG \times^{\bB} V)) \times \bG/\bB$. In particular this scheme 
is separated, hence as in the second proof of Proposition~\ref{prop:flag-var} we deduce that $\Proj^{\mathbb{X}}(A_V)$ is separated.

%Since $\bG/\bB$ is proper, this fact also implies 
Next, we show
that the natural morphism
\[
 \bG \times^{\bB} V \longto (\bG \times^{\bB} V)_\aff
\]
is proper.
In fact this map fits in the commutative diagram
\[
 \xymatrix{
 \bG \times^{\bB} V \ar[rd]_-{\eqref{eqn:morph-vb}} \ar[rr] && \Spec(\calO(\bG \times^{\bB} V)) \ar[ld] \\
 & \widetilde{V}. &
 }
\]
Here the left diagonal arrow 
%is the composition of the closed immersion~\eqref{eqn:embedding-vb} with the natural projection $\bG \times^{\bB} \widetilde{V} \cong \bG/\bB \times \widetilde{V} \to \widetilde{V}$, hence 
is proper, and the right diagonal arrow is induced by the corresponding morphism on function algebras. The latter morphism is separated (as is any morphism whose domain is affine, see~\cite[\href{https://stacks.math.columbia.edu/tag/01KN}{Tag 01KN}]{stacks-project}), hence the horizontal morphism is proper by~\cite[\href{https://stacks.math.columbia.edu/tag/01W6}{Tag 01W6}]{stacks-project}.

Finally, we observe that if $\lambda \in \mathbb{X}_{++}$ the invertible sheaf $\calO_{\bG \times^{\bB} V}(\lambda)$ is ample by~\cite[\href{https://stacks.math.columbia.edu/tag/0892}{Tag 0892}]{stacks-project} and the similar claim for $\bG/\bB$. As in the second proof of Proposition~\ref{prop:flag-var}, these arguments show that the conditions in~\cite[Corollary~4.6]{BS03} are satisfied, and this statement provides the desired identification of $\bG \times^{\bB} V$ with $\Proj^{\mathbb{X}}(A_V)$.
\end{proof}

\begin{rmk}
\label{rmk:Springer-covering}
Recall the notation from Remark~\ref{rmk:flag-covering} and, for $\sigma \in E$, denote by $\tilde{g}_\sigma$ the image of $g_\sigma$ in $A_V$. We claim that
\[
(A_V)_\dag \subset \Rad \left( \sum_{\sigma \in E} (A_V) \cdot \tilde{g}_\sigma \right),
\]
so that as in this remark we have
\[
\Proj^{\bbX}(A_V) = \bigcup_{\sigma \in E} D_\dag(\tilde{g}_\sigma).
\]
(This generalizes~\eqref{eqn:covering-Proj-AV}, proved under the assumption that $\bG$ is semisimple and simply connected.) In fact, as in the proof of Lemma~\ref{lem:relevant-flag} the degree of any relevant element $f$ in $A_V$ belongs to $\bbX_{++}$; to conclude, it therefore suffices to prove that we have
\[
(A_V)_\lambda \subset \Rad \left( \sum_{\sigma \in E} (A_V) \cdot \tilde{g}_\sigma \right) \quad \text{for any $\lambda \in \bbX_{++}$.}
\]
Fix such a $\lambda$, and let $a \in (A_V)_\lambda$. If $N$ is as in Lemma~\ref{lem:restriction-vb-surjective}, then $a^N$ belongs to the image of $\calO(\widetilde{V}) \otimes A_{N\lambda}$; i.e.~there exists a finite set $I$ and elements $(b_i : i \in I)$ in $\calO(\widetilde{V})$ and $(c_i : i \in I)$ in $A_{N\lambda}$ such that
$a^N$ is the image of $\sum_{i \in I} b_i \otimes c_i$.
By~\eqref{eqn:Alambda-flag}, for any $i \in I$ there exists $n_i \in \Z_{>0}$ such that $(c_i)^{n_i} \in \sum_\sigma A\cdot g_\sigma$. Then if $M \in \Z$ satisfies $M > |I| \cdot (\max_i n_i - 1)$ we have $a^{MN} \in \sum_{\sigma} (A_V) \cdot \tilde{g}_\sigma$, which concludes the proof.

%In the course of the proof of Proposition~\ref{prop:vb-flag-Proj} we have produced an open covering of the scheme $\Proj^{\bbX}(A_V)$ by affine open subschemes associated with elements in $A_V$ which are maximally relevant. (See~\eqref{eqn:covering-Proj-AV} for the case $\bG$ is semisimple and simply connected. The general case can be deduced as in Remark~\ref{rmk:flag-covering}.)
%
% Fix a family
% $(f_j : j \in J)$ of elements of $A$ as in Remark~\ref{rmk:flag-covering}. By Remark~\ref{rmk:open-basechange}, the scheme $\Proj^{\bbX}( \calO(\bG \times^{\bB} V) \otimes_\bk A)$ is covered by the open subschemes $(D_\dag(f_j) : j \in J)$. Using~\eqref{eqn:Proj-quotient-open} we deduce that $\Proj^{\bbX}(A_V)$ is covered by the open subschemes $(D_\dag(\tilde{f}_j) : j \in J)$, where $\tilde{f}_j$ is the image of $f_j$ in $A_V$. Note that here again, each $\tilde{f}_j$ is maximally relevant.
\end{rmk}

%-----------------------------------------------
\subsection{The Springer resolution}
\label{ss:Springer}
%-----------------------------------------------

We continue with the setting of~\S\ref{ss:vector-bundles}, in the special case when $\widetilde{V}$ is the Lie algebra $\bg$ of $\bG$ and $V$ is the Lie algebra $\bu$ of the unipotent radical $\bU$ of $\bB$. In this case the vector bundle $\bG \times^\bB \bu$ is called the \emph{Springer resolution}, and usually denoted $\Spr$. (The name is justified by the fact that, under appropriate technical assumptions, $\Spr$ is a resolution of singularities of the nilpotent cone $\mathcal{N} \subset \bg$ of $\bG$, first introduced by Springer.) This variety appears in numerous works in Geometric Representation Theory, among which~\cite{ab} (in the case when $\bk$ is an algebraic closure of the field $\mathbb{Q}_\ell$). In that reference the authors introduce an $\bbX$-graded ring, the spectrum of which is denoted $\hat{\tilde{\mathcal{N}}}_{af}$ in~\cite[\S 3.1]{ab}. This construction is reproduced in~\cite[\S 6.2.2]{ar}, where the spectrum is denoted $\widehat{\mathcal{N}}_{\mathcal{X}}$. The definition makes sense over an arbitrary algebraically closed field; here, we will denote this ring by $A_{\bu}'$. By construction, we have morphisms of $\bbX$-graded rings
\[
\calO(\bg) \otimes A \longto A_{\bu}' \to A_\bu,
\]
the first of which is surjective. (The main point of this construction is that $A_\bu'$ is a somewhat explicit quotient of $\calO(\bg) \otimes A$, whereas $A_\bu$ does not have a very explicit description.)

\begin{lem}
There exists a canonical isomorphism
\[
\Proj^{\bbX}(A_\bu) \simto \Proj^{\bbX}(A'_\bu).
\]
\end{lem}

\begin{proof}
The morphism $A_\bu' \to A_\bu$ provides a rational morphism $\Proj^{\bbX}(A_\bu) \dashrightarrow \Proj^{\bbX}(A_\bu')$, see Proposition~\ref{foncto}. As in the proof of Proposition~\ref{prop:vb-flag-Proj}, it follows from~\eqref{eqn:covering-Proj-AV} that this morphism is defined everywhere. To show that it is an isomorphism, we will construct an isomorphism $\Proj^{\bbX}(A_\bu') \cong \Spr$, under which this morphism is the identity.

Consider the morphisms
\[
 \Spec(A'_\bu) \hookrightarrow \Spec(\calO(\bg) \otimes_\bk A) = \bg \times_{\Spec(\bk)} \Spec(A) \to \Spec(A),
\]
where the left arrow is the closed immersion induced by the surjection $\calO(\bg) \otimes A \to A_{\bu}'$, and the right one is the obvious projection.
By Remark~\ref{rmk:relevant-projection}, $V((\calO(\bg) \otimes_\bk A)_\dag)$ is the preimage of $V(A_\dag)$ under the second morphism, and by Remark~\ref{rmk:preimage-Adag-surj} $V((A'_\bu)_\dag)$ is the preimage of $V((\calO(\bg) \otimes_\bk A)_\dag)$ under the first morphism.
Combining these informations, we obtain that
% by the same arguments as in Lemma~\ref{lem:relevant-flag} and surjectivity of the morphism $\calO(\bg) \otimes_\bk A \twoheadrightarrow A''$, one sees that
% \[
% A''_\dagger = \bigoplus_{\lambda \in \bbX_{++}} A''_\lambda = A_\dagger \cdot A''.
% \]
% Hence 
$\Spec(A'_\bu) \smallsetminus V((A'_\bu)_\dagger)$ is the preimage of $\Spec(A) \smallsetminus V(A_\dagger)$, which we have identified with $\bG/\bU$ in~\eqref{eqn:G/U}. By~\cite[Equation~(6.2.10)]{ar} this preimage is the variety denoted $\widehat{\mathcal{N}}$ in~\cite[\S 6.2]{ar}, whose quotient by the action of $\bT$ is $\Spr$ (see the discussion preceding~\cite[Equation~6.2.7]{ar}). This concludes the proof.
\end{proof}

\begin{rmk}
\phantomsection
\label{rmk:Springer-AB}
\begin{enumerate}
\item
The second part of the proof can alternatively be replaced by an argument based on the fact that the morphism $(A'_\bu)_\lambda \to (A_\bu)_\lambda$ is an isomorphism if $\lambda$ is sufficiently dominant, see~\cite[Lemma~6.2.4]{ar}.
\item
The same comments as in Remark~\ref{rmk:Springer-covering} apply in this case.
\end{enumerate}
\end{rmk}

\section{Quasi-coherent sheaves on Proj schemes} \label{s:qcohproj}
%%%%%%%%%%%%%%%%%%%%%%%%%%%%%%%%%%%%%%%%

In this section we study quasi-coherent shaves on Proj schemes, generalizing classical results for $\Proj$ schemes of $\bbN$-graded rings. Our presentation follows and extends \cite[\href{https://stacks.math.columbia.edu/tag/01MJ}{Tag 01MJ}]{stacks-project} and \cite[\href{https://stacks.math.columbia.edu/tag/01MM}{Tag 01MM}]{stacks-project}. Other pioneering works on quasi-coherent sheaves on Brenner--Schr\"oer Proj schemes include \cite{MR24}.

Most of the results of this section have obvious analogues in the relative setting of~\S\ref{ss:relative}. We leave it to the reader to formulate these variants and adapt the proofs.

%--------------------------------------
\subsection{Sheaves associated with graded modules}
\label{ss:sheaves-gr-modules}
%--------------------------------------

We proceed with the notation of~\S\ref{subsectprojring}; in particular we fix a finitely generated abelian group $M$ and a commutative $M$-graded ring $A$. Recall also the notation $\calF_A$. 

Let $Q$ be an $M$-graded $A$-module. For any homogeneous multiplicative subset $S \subset A$, we have considered in Definition~\ref{potionsdef} the $A_{(S)}$-module $Q_{(S)}$. The following fact is immediate by glueing of quasi-coherent sheaves.

%Facts \ref{factMtilde} and \ref{0part} extends part of \cite[\href{https://stacks.math.columbia.edu/tag/01M7}{Tag 01M7}]{stacks-project}.

\fact
\label{factMtilde}
There exists 
%on $X=\mathrm{Proj} (A)$ 
a unique quasi-coherent $\mathcal{O}_{\Proj^M(A)}$-module $\widetilde{Q}$ such that
$\Gamma \bigl( D_\dag(S)  , \widetilde{Q} \bigr) = Q_{(S)}$
for every $S \in \calF_A$.
\xfact

% \pf
% This is immediate by glueing.
% \xpf 

It is clear that the assignment $Q \mapsto \widetilde{Q}$ is functorial. More specifically,
denoting by $\Mod^M(A)$ the abelian category of $M$-graded $A$-modules, and by $\QCoh(\Proj^M(A))$ the abelian category of quasi-coherent sheaves on the scheme $\Proj^M(A)$, this assignment defines a functor
\[
\Mod^M(A) \to \QCoh(\Proj^M(A)),
\]
which is exact by exactness of localization. Note that this functor commutes with all colimits. (This follows from the facts that restriction to open subscheme and localization commute with colimits.)

An $M$-graded $A$-module $Q$ will be called \emph{negligible} if $\widetilde{Q}=0$. We will denote by $\Mod^M(A)_{\mathrm{neg}}$ the full subcategory of $\Mod^M(A)$ whose objects are the negligible modules. Since the functor $Q \mapsto \widetilde{Q}$ is exact, this is a Serre subcategory, see~\cite[\href{https://stacks.math.columbia.edu/tag/02MQ}{Tag 02MQ}]{stacks-project}, and our functor factors through an exact functor
\[
\mathsf{L} : \Mod^M(A) / \Mod^M(A)_{\mathrm{neg}} \longto \QCoh(\Proj^M(A)).
\]
see~\cite[\href{https://stacks.math.columbia.edu/tag/02MS}{Tag 02MS}]{stacks-project}.

The following fact is again immediate, by glueing of sections of quasi-coherent sheaves.

\fact
\label{0part}
There is a canonical morphism of $A_0$-modules 
\[
Q_0 \longto \Gamma \bigl( \Proj^M(A),\widetilde{Q} \bigr)
\]
such that for any $S \in \calF_A$ the composition $Q_0 \to  \Gamma ( \Proj^M (A),\widetilde{Q} )\to \Gamma ( D_\dag(S) , \widetilde{Q})=Q_{(S)}$ coincides with the map $ Q_0 \to Q_{(S)}$ given by $x \mapsto \frac{x}{1}$.
\xfact

% \pf
% The morphisms $Q_0 \to Q_{(S)}$ given by $x \mapsto \frac{x}{1}$, where $S$ runs through $\calF_A$, are compatible by restriction, hence glue to the desired morphism.
% \xpf  

%Proposition \ref{proptensorgra} extends \cite[\href{https://stacks.math.columbia.edu/tag/01MK}{Tag 01MK}]{stacks-project}.

The next proposition studies the relation between tensor products of graded modules and of quasi-coherent sheaves.

\prop
\label{proptensorgra} 
%Let $f $ be a relevant family. 
Let $P$, $Q$ be graded $A$-modules. There is a canonical morphism of quasi-coherent $\mathcal{O}_{\Proj^M(A)}$-modules 
  \[
  \widetilde P \otimes_{\mathcal{O}_ X} \widetilde Q \longrightarrow \widetilde{P \otimes_ A Q}
  \]
which induces, for any $S \in \calF_A$, the canonical map 
\[
P_{(S)} \otimes_{A_{(S)}} Q_{(S)} \longrightarrow (P \otimes_A Q)_{(S)}
\]
 on sections over $D_\dag(S)$. Moreover, the diagram
  \[ 
  \xymatrix{
  P_0 \otimes_{A_0} Q_0 \ar[r] \ar[d] &  (P \otimes_A Q)_0 \ar[d] \\  \Gamma (\Proj^M(A), \widetilde P \otimes_{\mathcal{O}_ X} \widetilde Q) \ar[r] &  \Gamma (\Proj^M(A), \widetilde{P \otimes_A Q}) 
  } 
  \]
 commutes, where the upper horizontal arrow is the natural map, and the vertical ones are induced by the morphisms from Fact~\ref{0part}.
\xprop

\pf
 Constructing a morphism as displayed is equivalent to constructing an $\mathcal{O}_{\Proj^M(A)}$-bilinear map 
  $\widetilde P \times \widetilde Q \longrightarrow \widetilde{P \otimes_A Q}$,
  see~\cite[\href{https://stacks.math.columbia.edu/tag/01CA}{Tag 01CA}]{stacks-project}. It suffices to define this map on sections over the opens $(D_\dag(S) : S \in \calF_A)$ compatible with restriction maps. On $D_\dag(S)$, with $S= \{a_i : i \in I \}$, we use the $A_{(S)}$-bilinear map $P_{(S)} \times Q_{(S)} \to (P \otimes_A Q)_{(S)}$ given by $(\frac{x}{a^\nu}, \frac{y}{a^{\nu'}}) \mapsto \frac{x \otimes y}{a^{\nu + \nu'}}$. The commutation of the diagram follows from definitions. 
\xpf

In general, the morphism of Proposition~\ref{proptensorgra} is not an isomorphism, as seen already in the $\N$-graded setting, see~\cite[\href{https://stacks.math.columbia.edu/tag/01ML}{Tag 01ML}]{stacks-project}.

\begin{rmk}
\label{rmk:tilde-equiv}
Consider the setting of Remark~\ref{rmk:action-gp-scheme}, and denote by $\Mod^{M,H}(A)$ the category of $H$-equivariant $M$-graded $A$-modules. Let also $\Mod^{M,H}(A)_{\mathrm{neg}}$ be the Serre subcategory consisting of objects which are negligible as $M$-graded $A$-modules. On the other hand, consider the category $\QCoh^H(\Proj^M(A))$ of $H$-equivariant quasi-coherent sheaves on $\Proj^M(A)$. (For a review of equivariant quasi-coherent sheaves, see e.g.~\cite[Appendix~A]{mr}.) Then very similar considerations to those above allow to construct an exact functor
$\Mod^{M,H}(A) \to \QCoh^H(\Proj^M(A))$
which factors through the quotient $\Mod^{M,H}(A) / \Mod^{M,H}(A)_{\mathrm{neg}}$.
\end{rmk}

%--------------------------------------
\subsection{Twisting sheaves}
\label{ss:twisting-sheaves}
%--------------------------------------

%  Definitions \ref{twist} and \ref{twistmodule} immediately extends \cite[\href{https://stacks.math.columbia.edu/tag/01MN}{Tag 01MN}]{stacks-project}.

We now define the versions in our setting of the twisting sheaves from~\cite[\href{https://stacks.math.columbia.edu/tag/01MN}{Tag 01MN}]{stacks-project}. 
If $Q$ is a graded $A$-module and $\alpha \in M$, we will denote by $Q(\alpha)$ the $M$-graded $A$-module which coincides with $M$ as an $A$-module, but with the $M$-grading defined by $(Q(\alpha))_\beta = Q_{\alpha+\beta}$ for $\beta \in M$.
%Recall the shift-of-degree functor $(-)(\alpha)$ defined in~\S\ref{ss:examples}.

\defi[(Twisting sheaves)]
\label{twist}
Let $\alpha \in M$.
\begin{enumerate}
\item
The quasi-coherent sheaf $\widetilde{A(\alpha)}$ on $\Proj^M(A)$ is 
%called the \emph{$\alpha$-th-twist of the structure sheaf} of $\mathrm{Proj}^M(A)$, and is 
denoted $\mathcal{O}_{\Proj^M(A)} (\alpha)$.
\item
If $\calQ$ is a sheaf of $\calO_{\mathrm{Proj}^M(A)}$-modules, we set $\calQ (\alpha) = \calO_{\mathrm{Proj}^M(A)} (\alpha ) \otimes_{\calO_{\mathrm{Proj}^M(A)}} \calQ$.
\end{enumerate}
\xdefi 

Recall that if $A$ is noetherian as a non-graded ring, then by~\cite[Lemma~2.4]{BS03} the ring $A_0$ is noetherian, and by~\cite[Proposition~2.5]{BS03} the canonical morphism $\Proj^M(A) \to \Spec(A_0)$ is of finite type; as a consequence, $\Proj^M(A)$ is a noetherian scheme.

\begin{lem}
\label{lem:tilde-noetherian-coherent}
Assume that $A$ is a noetherian ring.
\begin{enumerate}
\item
\label{it:twisting-coherent}
For any $\alpha \in M$, the quasi-coherent sheaf $\mathcal{O}_{\mathrm{Proj}^M(A)} (\alpha)$ is coherent.
\item
\label{it:tilde-fg-coherent}
If $Q$ is a finitely generated $M$-graded $A$-module, then $\widetilde{Q}$ is coherent.
\end{enumerate}
\end{lem}

\begin{proof}
\eqref{it:twisting-coherent} We need to prove that for any $S \in \calF_A$ the $A_{(S)}$-module $(A(\alpha))_{(S)}$ is finitely generated. Here $A_{(S)}$, resp.~$(A(\alpha))_{(S)}$, is the degree-$0$, resp.~degree-$\alpha$, component in the localization $A_S$. By~\cite[Lemma~2.4]{BS03}, $A_S$ is noetherian. By a simple argument (see e.g.~\cite[Lemma~2.2]{gy}), this implies that each of its graded components is finitely generated over its degree-$0$ part, which concludes the proof.
%, and is finitely generated as an $A_{(S)}$-algebra. If $(a_i : i \in I)$ is a finite family of homogeneous elements that generate $A_S$ as an $A_{(S)}$-algebra, then 

\eqref{it:tilde-fg-coherent} Any finitely generated $M$-graded $A$-module is a quotient of a finite direct sum of modules $A(\alpha)$ with $\alpha \in M$. The claim therefore follows from~\eqref{it:twisting-coherent} and exactness of the functor $Q \mapsto \widetilde{Q}$, since finite direct sums and quotients of coherent sheaves by quasi-coherent subsheaves are coherent.
\end{proof}

We now drop the assumption that $A$ is noetherian.
Note that $\mathcal{O}_{\mathrm{Proj}^M(A)} (\alpha)$ is not an invertible sheaf in general, even in the $\N$-graded setting. If $\alpha, \alpha' \in M$,
%. Fact \ref{morphismalphaalphaprime} extends \cite[\href{https://stacks.math.columbia.edu/tag/01MO}{Tag 01MO}]{stacks-project}.
%\fact \label{morphismalphaalphaprime} 
since $A(\alpha) \otimes_A A(\alpha ') = A(\alpha + \alpha ')$, Proposition \ref{proptensorgra} implies that there is a canonical map
\begin{equation}
\label{eqn:morphismalphaalphaprime}
\mathcal{O}_ {\mathrm{Proj}^M(A)}(\alpha) \otimes_{\mathcal{O}_ {\mathrm{Proj}^M(A)}} \mathcal{O}_ {\mathrm{Proj}^M(A)}(\alpha') \longrightarrow \mathcal{O}_ {\mathrm{Proj}^M(A)}(\alpha + \alpha').
\end{equation}
These maps define on
\begin{equation}
\label{eqn:sum-twisting-sheaves}
\bigoplus_{\alpha \in M} \mathcal{O}_ {\mathrm{Proj}^M(A)}(\alpha)
\end{equation}
a structure of an $M$-graded sheaf of $\mathcal{O}_ {\mathrm{Proj}^M(A)}$-algebras, and an $M$-graded ring structure on
\begin{equation}
\label{eqn:sum-twisting-sheaves-sections}
\bigoplus_{\alpha \in M} \Gamma (\Proj^M(A) , \calO_{\Proj^M(A)} ( \alpha ) ).
\end{equation}
Note that the morphism~\eqref{eqn:morphismalphaalphaprime} is not an
%These maps are not 
isomorphism in general. (Again, this can already be false in the $\N$-graded setting.)
%\xfact 

More generally, if $Q$ is an $M$-graded $A$-module, and if $\alpha,\alpha' \in M$, we also have $A(\alpha) \otimes_A Q(\alpha ') = Q(\alpha + \alpha ')$; Proposition~\ref{proptensorgra} therefore provides a canonical morphism
\begin{equation}
\label{eqn:morphismalphaalphaprime-Q}
\mathcal{O}_ {\mathrm{Proj}^M(A)}(\alpha) \otimes_{\mathcal{O}_ {\mathrm{Proj}^M(A)}} \widetilde{Q(\alpha')} \to \widetilde{Q(\alpha+\alpha')}.
\end{equation}
These maps define on
\[
\bigoplus_{\alpha \in M} \widetilde{Q(\alpha)}
\]
a structure of an $M$-graded sheaf of modules over~\eqref{eqn:sum-twisting-sheaves}, and on
\[
\bigoplus_{\alpha \in M} \Gamma (\Proj^M(A), \widetilde{Q(\alpha)})
\]
a structure of an $M$-graded module over~\eqref{eqn:sum-twisting-sheaves-sections}.

\fact 
\phantomsection
\label{fact:morp-graded-rings-modules}
\begin{enumerate}
\item
\label{it:morp-graded-rings}
There is a canonical morphism of $M$-graded rings
\[
A \to \bigoplus_{\alpha \in M} \Gamma \Bigl( \Proj^M(A) , \calO_{\Proj^M(A)} ( \alpha ) \Bigr).
\]
\item
\label{it:morp-graded-modules}
For any graded $A$-module $Q$ we have a canonical morphism of $M$-graded abelian groups
%
%\xfact
%\pf Note that $A_{\alpha} = (A (\alpha) )_0$. Now we combine Fact \ref{0part} and Fact \ref{morphismalphaalphaprime}. 
%\xpf 
%
%\fact  \label{factQalphatotilde}
% There is a canonical map 
\[ 
Q \to \bigoplus_{\alpha \in M} \Gamma \Bigl( X,\widetilde{Q(\alpha)} \Bigr).
\]
which is a morphism of $M$-graded $A$-modules with respect to the structure on the right-hand side provided by~\eqref{it:morp-graded-rings}.
\end{enumerate}
\xfact
 
\pf
%Note that $Q_{\alpha} = (Q (\alpha) )_0$. Now we apply 
The morphisms are given by Fact~\ref{0part}, after noticing that $A_{\alpha} = (A (\alpha) )_0$ and $Q_{\alpha} = (Q (\alpha) )_0$. The fact that the morphism in~\eqref{it:morp-graded-rings} is a ring morphism follows from Proposition~\ref{proptensorgra}.
\xpf  

%--------------------------------------
\subsection{Graded modules associated with sheaves}
\label{ss:gr-modules-sheaves}
%--------------------------------------

For any $\calO_{\Proj^M(A)}$-module $\calQ$, the morphism~\eqref{eqn:morphismalphaalphaprime} induces a morphism
%Fact \ref{morphismtensoring} extends \cite[\href{https://stacks.math.columbia.edu/tag/03GJ}{Tag 03GJ}]{stacks-project}.
%\fact \label{morphismtensoring}
%Tensoring with an arbitrary $\calO_X$-module $\calQ$ we get maps
\begin{equation}
\label{eqn:morphismalphaalphaprime-calQ}
\calO_{\Proj^M(A)}(\alpha)\otimes_{\calO_{\Proj^M(A)}} \calQ(\alpha')\longrightarrow\calQ(\alpha + \alpha '),
\end{equation}
%\xfact
and these morphisms define on the direct sum
\[
\bigoplus_{\alpha \in M} \calQ(\alpha)
\]
the structure of an $M$-graded sheaf of modules over~\eqref{eqn:sum-twisting-sheaves}, and on
\[
\bigoplus_{\alpha \in M} \Gamma (\Proj^M(A), \calQ(\alpha))
\]
a structure of an $M$-graded module over the ring~\eqref{eqn:sum-twisting-sheaves-sections}, hence over the ring $A$ by Fact~\ref{fact:morp-graded-rings-modules}\eqref{it:morp-graded-rings}. Denoting by $\Mod(\calO_{\mathrm{Proj}^M(A)})$ the abelian category of sheaves of $\calO_{\mathrm{Proj}^M(A)}$-modules, 
%and by $\Mod^M(A)$ the abelian category of $M$-graded $A$-modules, 
this construction provides a functor
\[
\Gamma_\bullet : \Mod(\calO_{\mathrm{Proj}^M(A)}) \longto \Mod^M(A).
\]

%\fact 
%For an arbitrary $\calO_X$-module $\calQ$ the maps of Fact \ref{morphismtensoring} give a graded module structure
%\[\bigoplus_{\alpha \in M } \Gamma ( X , \calO_X  (\alpha ) ) \times \bigoplus_{\alpha' \in M } \Gamma ( X , \calQ (\alpha ')) \longrightarrow \bigoplus_{\alpha' \in M } \Gamma ( X , \calQ (\alpha ')) \]
%and via Fact \ref{mapofgradedwithsection} also a $A$-module structure. 
%\xfact 

%\fact \label{tildetilde}  Given any graded $A$-module $Q$ we have $Q(\alpha)=Q \otimes_A A( \alpha)$. Hence we get maps \[ \widetilde{Q} (\alpha ) = \widetilde{Q} \otimes_{\calO_X} \calO_X ( \alpha ) \to \widetilde{Q (\alpha )}  .\]
%\xfact
%\fact \label{factQtotilde}
%Fact \ref{factQalphatotilde} defines a map of graded $A$-modules
%$Q\longrightarrow \bigoplus_{\alpha \in M} \Gamma (X,\widetilde{Q(\alpha)}).$
%\xfact 

In the following proposition we use the notation $M[\underline{S}]$ introduced in~\S\ref{ss:relevant-families}.

\prop
\label{prop:twisting-inversible-open} 
Let $S \in \calF_A$, and let $\alpha \in M[\underline{S}]$.
\begin{enumerate}
\item
\label{it:twisting-inversible-open-1}
The sheaf $\calO_{\Proj^M(A)} (\alpha)_{| D_\dag(S)}$ is invertible, and in fact isomorphic to $\calO_{D_\dag(S)}$.
\item
\label{it:twisting-inversible-open-2}
For any $M$-graded $A$-module $Q$, the morphism
\[
\widetilde{Q} (\alpha)_{|D_\dag(S)}  \longto \widetilde{Q(\alpha)}_{|D_\dag(S)} 
\]
obtained by restriction from~\eqref{eqn:morphismalphaalphaprime-Q} is an isomorphism.
\item
\label{it:twisting-inversible-open-3}
For any $\alpha' \in M$, the morphism
%he maps of Fact \ref{morphismalphaalphaprime} restricted to $U_f$ 
\[
\left( \calO_{\Proj^M(A)} (\alpha)_{| D_\dag(S)} \right) \otimes_{\calO_{D_\dag(S)}} \left( \calO_{\Proj^M(A)} (\alpha')_{| D_\dag(S)} \right)\longto \calO_{\Proj^M(A)}( \alpha + \alpha ')_{| D_\dag(S)}
\] 
obtained by restriction from~\eqref{eqn:morphismalphaalphaprime} is an isomorphism.
\item
\label{it:twisting-inversible-open-4}
For any $\calO_{\Proj^M(A)}$-module $\calQ$ and any $\alpha' \in M$, the morphism
\[
 \left( \calO_{\Proj^M(A)}(\alpha)_{|D_\dag(S)} \right) \otimes_{\calO_{D_\dag(S)}} \left( \calQ(\alpha')|_{D_\dag(S)} \right) \longto \calQ (\alpha + \alpha')_{|D_\dag(S)}
\]
obtained by restriction from~\eqref{eqn:morphismalphaalphaprime-calQ} is an isomorphism.
\end{enumerate}
%the maps of Fact \ref{morphismtensoring} restricted to $U_f$,
%\[ \calO_X ( \alpha ) |_{U_f} \otimes_{\calO_{U_f}} \calQ (\alpha ' )|_{U_f} \to \calQ ( \alpha + \alpha ')|_{U_f} ,\] 
%
%and the maps of Fact \ref{tildetilde} restricted to $U_f$,
%\[ \widetilde{Q} ( \alpha ) |_{U_f}  \to \widetilde{ Q( \alpha)}|_{U_f}  \]
%are isomorphisms. 
\xprop 

\pf 
%We can assume that $\underline{S}=S$.
%Write $S= \{ a_i : i \in I \}$. 
Since $\alpha \in M[\underline{S}]$, there exists an invertible element $a \in A_S$ of degree $\alpha$ (see Proposition~\ref{prop:char-relevant}\eqref{it:char-relevant}). Then
%nonzero elements $a,b \in S$ such that $\deg(a)-\deg(b)=\alpha$. Then $\frac{a}{b}$ is a well defined element of $A_S$, and
%$\nu \in \mathbf{Z}_I $ such that $\deg (a^\nu )= \alpha$. 
the map $x \mapsto a \cdot x$ induces an isomorphism of $A_{(S)}$-modules $A_{(S)} \cong (A(\alpha))_{(S)} = (A_S)_\alpha$. We deduce~\eqref{it:twisting-inversible-open-1}.
Similarly, given a graded $A$-module $Q$, the map $q \mapsto a \cdot q$ induces an isomorphism of $A_{(S)}$-modules $Q_{(S)} \cong (Q(\alpha))_{(S)} = (Q_S)_\alpha$, which implies~\eqref{it:twisting-inversible-open-2}. The statement in~\eqref{it:twisting-inversible-open-3} is the special case of~\eqref{it:twisting-inversible-open-2} where $Q=A(\alpha')$. Finally,~\eqref{it:twisting-inversible-open-4} is a direct consequence of~\eqref{it:twisting-inversible-open-3}.
% (cf. \ref{factisoano} and \ref{coroisoano}). 
% This shows that that $\mathcal{O}_X (\alpha) |_{U_f}$ is invertible. Now since $\bigcup_{ f \in \mathcal{F} } U_f = X$, we conclude that $\mathcal{O}_X (\alpha)$ is an invertible $\mathcal{O}_X$-module. The other assertions are similarly proved.
\xpf

\begin{lem}
\label{lem:morph-tilde-Gamma}
For any $\calQ \in \Mod(\calO_{\mathrm{Proj}^M(A)})$, there exists a canonical (in particular, functorial) morphism
$\widetilde{\Gamma_\bullet(\calQ)} \to \calQ$
in $\Mod(\calO_{\mathrm{Proj}^M(A)})$.
\end{lem}

\begin{proof}
Since $\Proj^M(A)$ is covered (by definition) by the open subschemes $D_\dag(S)$ for $S \in \calF_A$ not containing $0$, to prove the lemma it suffices to construct morphisms
\begin{equation}
\label{eqn:morph-tilde-Gamma-1}
\widetilde{\Gamma_\bullet(\calQ)}_{| D_\dag(S)} \to \calQ_{| D_\dag(S)}
\end{equation}
for such $S$, which coincide on intersections of such open subschemes. And given $S$, by basic properties of quasi-coherent sheaves on affine schemes (see~\cite[\href{https://stacks.math.columbia.edu/tag/01I7}{Tag 01I7}]{stacks-project}), to define such a morphism it suffices to define a morphism of $A_{(S)}$-modules
\begin{equation}
\label{eqn:morph-tilde-Gamma-2}
\Gamma_\bullet(\calQ)_{(S)} \to \Gamma(D_\dag(S), \calQ).
\end{equation}
Now an element of $\Gamma_\bullet(\calQ)_{(S)}$ can be represented by a fraction $\frac{m}{s}$ where $s \in S$ and $m \in \Gamma_\bullet(\calQ)_{\deg(s)}$, i.e.~$m$ is a global section of $\calQ(\deg(s))$. Consider
\[
m_{|D_\dag(S)} \in \Gamma(D_\dag(S), \calQ(\deg(s))).
\]
The element $\frac{1}{s} \in A_S$ has degree $-\deg(s)$, hence defines a section of $\calO_{\Proj^M(S)}(-\deg(s))$ on $D_\dag(S)$, which we denote $s^{-1}$. Then the product $m_{|D_\dag(S)} \otimes s^{-1}$ defines a global section of
\[
\calQ(\deg(s))_{|D_\dag(S)} \otimes_{D_\dag(S)} \calO_{\Proj^M(A)}(-\deg(s))_{|D_\dag(S)} \cong \calQ_{|D_\dag(S)},
\]
where the isomorphism is provided by Proposition~\ref{prop:twisting-inversible-open}\eqref{it:twisting-inversible-open-4}.
%We take this section as a definition of the image of $\frac{m}{s}$. 

One easily checks that this section is independent of the representation of the element of $\Gamma_\bullet(\calQ)_{(S)}$ as a fraction, hence that this process defines a map~\eqref{eqn:morph-tilde-Gamma-2}, and then that this map is a morphism of $A_{(S)}$-modules. One also easily sees that the associated morphisms~\eqref{eqn:morph-tilde-Gamma-1} glue on the intersections of open subschemes $D_\dag(S)$, hence provide the desired morphism of sheaves.
\end{proof}

\subsection{Maximally relevant families}
\label{ss:max-relevant}
%--------------------------------------

A family $S \in \calF_A$ will be called \emph{maximally relevant} if $M[\underline{S}]=M$. A homogeneous element $f \in A$ will be called maximally relevant if $f^\N$ is maximally relevant.  We will denote by $\calF_A^{\mathrm{m}} \subset \calF_A$ the subset consisting of maximally relevant families. In this subsection we will explore various consequences of the following condition (which may or may not hold, depending on $A$):
\begin{equation}
\label{eqn:cover-max}
\Proj^M(A) = \bigcup_{S \in \calF_A^{\mathrm{m}}} D_\dag(S).
\end{equation}
The same comments as for~\eqref{eqn:singletonproj} show that this condition is equivalent to the property that $\Proj^M(A)$ is covered by the open subschemes $D_\dag(f)$ with $f$ maximally relevant.

\begin{rmk}
\label{rmk:max-cover}
%Let us note that condition~\eqref{eqn:cover-max} holds at least in some cases of interest.
If in the setting of Corollary~\ref{cor:covering-Proj} one can choose the $a_i$'s to be maximally relevant, then condition~\eqref{eqn:cover-max} holds. This setting covers at least some cases of interest, as follows.
\begin{enumerate}
 \item Assume that $M=\Z$, and that $A$ is generated by $A_1$ as an $A_0$-algebra. 
% (In particular, we are in the $\N$-graded case.) 
 If $(a_i : i \in I)$ is a family of generators of $A_1$ as an $A_0$-module, since any non nilpotent relevant element must have positive degree the condition of Corollary~\ref{cor:covering-Proj} is satisfied.
 %we have $\Proj^{\Z}(A) = \bigcup_{i \in I} D_\dag(a_i)$ by Corollary~\ref{cor:covering-Proj}, and 
 Clearly, each $a_i$ is maximally relevant.
 \item 
 Consider the setting of~\S\ref{ss:flag}. We have explained in Remark~\ref{rmk:flag-covering} how to construct 
 a family $(g_\sigma : \sigma \in E)$ for which Corollary~\ref{cor:covering-Proj} applies. We claim that any $g_\sigma$ in this family is maximally relevant. In fact, fix $\sigma$, and recall the notation introduced in this remark. Since $\bbX$ is generated by $\bbX_+$ as a group, it suffices to justify that $\bbX[\underline{(g_\sigma)^\N}]$ contains $\bbX^+$. Let $\lambda \in \bbX$, and write $\lambda = \sum_\alpha m_\alpha \varpi_\alpha + \lambda_0$ with $m_\alpha \in \N$ and $\lambda_0 \in \widetilde{\bbX}_0$ in the decomposition~\eqref{eqn:X-sc} (for the group $\widetilde{\bG}$). If $a$ is a nonzero vector in the $1$-dimensional vector space $\widetilde{A}_{\lambda_0}$, then the element $a_0 \cdot \prod_\alpha (f^\alpha_{\sigma(\alpha)})^{m_\alpha}$ is of weight $\lambda$ hence belongs to $A$, and divides a power of $g_\sigma$ (in $\widetilde{A}$, hence in $A$ for weight reasons), which justifies that $\lambda \in \bbX[\underline{(g_\sigma)^\N}]$.
% an open covering of $\Proj^{\bbX}(A)$.
 \item
Consider the setting of~\S\ref{ss:vector-bundles}.
 In Remark~\ref{rmk:Springer-covering}, we have constructed a family $(\tilde{g}_\sigma : \sigma \in E)$ for which Corollary~\ref{cor:covering-Proj} applies. It follows from the case treated above using the morphism~\eqref{eqn:morph-A-AV} that each $\tilde{g}_\sigma$ is maximally relevant. Similar comments apply in the setting of~\S\ref{ss:Springer}.
\end{enumerate}
\end{rmk}

First, the following statement is a direct consequence of Proposition~\ref{prop:twisting-inversible-open}.

\coro 
\label{coro:twistings-line-bundles}
Assume that~\eqref{eqn:cover-max} is satisfied.
%\begin{equation}
%\label{eqn:cover-max}
%\Proj^M(A) = \bigcup_{S \in \calF_A^{\mathrm{m}}} D_\dag(S).
%\end{equation}
%there exists a set $\mathcal{F}$ of maximally relevant families of homogeneous elements of $A$ such that 
% $\bigcup_{ f \in \mathcal{F} } U_f = X$.
Then for any $ \alpha \in M$, the quasi-coherent $\calO_{\Proj^M(A)}$-module $\calO_{\Proj^M(A)} (\alpha)$ is an invertible sheaf. Moreover:
\begin{enumerate}
\item
the morphism~\eqref{eqn:morphismalphaalphaprime} is an isomorphism for any $\alpha,\alpha' \in M$;
\item
\label{it:max-rel-2}
the morphism~\eqref{eqn:morphismalphaalphaprime-Q} is an isomorphism for any $\alpha,\alpha' \in M$ and any $M$-graded $A$-module $Q$;
\item
the morphism~\eqref{eqn:morphismalphaalphaprime-calQ} is an isomorphism for any $\alpha,\alpha' \in M$ and any $\calO_{\Proj^M(A)}$-module $\calQ$.
\end{enumerate}
%an invertible $\mathcal{O}_X$-module and the maps (cf. Facts \ref{morphismalphaalphaprime}, \ref{morphismtensoring} and \ref{tildetilde})
%\[ \mathcal{O}_ X(\alpha) \otimes_{\mathcal{O}_ X} \mathcal{O}_ X(\alpha') \longrightarrow \mathcal{O}_ X(\alpha + \alpha'), \]
%
%\[\calO_X(\alpha)\otimes_{\calO_X} \calQ(\alpha')\longrightarrow\calQ(\alpha + \alpha '),\]
%
%\[ \widetilde{Q} (\alpha ) \to \widetilde{Q (\alpha )}  \]
%are isomorphisms.
% Thus the map of Fact \ref{factQalphatotilde} becomes a map
%
%\[ Q_{\alpha} \to \Gamma ( X , \widetilde{Q}(\alpha)) \]
%
%and the map of Fact \ref{factQtotilde} becomes a map
%\[ Q \to \bigoplus_{\alpha \in M } \Gamma ( X , \widetilde{Q}(\alpha))  . \] 
\xcoro

%\pf
%This is a Corollary of Proposition \ref{propinversi}.
%\xpf 

\begin{rmk}
 In the setting of~\S\ref{ss:flag}, for $\lambda \in \bbX$ the invertible sheaf $\calO_{\Proj^{\bbX}(A)}(\lambda)$ corresponds to the invertible sheaf $\calO_{\bG/\bB}(\lambda)$. A similar comment applies in the setting considered in~\S\ref{ss:vector-bundles}.
\end{rmk}

\begin{propo}
\label{prop:max-torsor}
 Assume that~\eqref{eqn:cover-max} is satisfied, and that $M$ is a free abelian group. Then the natural morphism
 \[
  \Spec(A) \smallsetminus V(A_\dag) \longto \Proj^M(A)
 \]
 (see~\eqref{eqn:geometric-quotient})
is a Zariski locally trivial principal bundle for the group scheme $\mathrm{D}_{\Spec(A_0)}(M)$.
\end{propo}

\begin{proof}
 By assumption $\Proj^M(A)$ is covered by the open subschemes $D_\dag(S)$ where $S$ is maximally relevant. Now for such $S$, the ring $A_S$ contains invertible elements in each degree. Choose an isomorphism $M=\Z^n$ and, for any $i \in \{1, \dots, n\}$, choose an invertible element $a_i$ in $A_S$ of degree the $i$-th vector in the canonical basis of $\Z^n$. Then we obtain an isomorphism of $M$-graded rings
 \[
  A_{(S)}[x_i^{\pm 1} : 1 \leq i \leq n] \simto A_S
 \]
sending $\prod_i (x_i)^{n_i}$ to $\prod_i (a_i)^{n_i}$. In other words we have
\[
 \Spec(A_S) \cong \Spec(A_{(S)}) \times_{\Spec(A_0)} \mathrm{D}_{\Spec(A_0)}(M),
\]
so that the restriction of our map to the preimage of $D_\dag(S)$ is a trivial $\mathrm{D}_{\Spec(A_0)}(M)$-bundle.
\end{proof}

\begin{rmk}
\label{rmk:tilde-torsor}
Assume that the conditions in Proposition~\ref{prop:max-torsor} are satisfied. Then this proposition implies that pullback induces an equivalence of categories between $\QCoh(\Proj^M(A))$ and the category of $\mathrm{D}_{\Spec(A_0)}(M)$-equivariant quasi-coherent sheaves on $\Spec(A) \smallsetminus V(A_\dag)$. On the other hand, the category of $M$-graded $A$-modules identifies with the category of $\mathrm{D}_{\Spec(A_0)}(M)$-equivariant quasi-coherent sheaves on $\Spec(A)$. Under these identifications, the functor $Q \mapsto \widetilde{Q}$ corresponds to restriction along the open immersion $\Spec(A) \smallsetminus V(A_\dag) \to \Spec(A)$.
\end{rmk}

For the next lemma, recall the notion of negligible $M$-graded $A$-module from~\S\ref{ss:sheaves-gr-modules}.

\begin{lem}
\label{lem:negligible}
 Assume that~\eqref{eqn:cover-max} is satisfied, and fix a subset $\calF \subset \calF_A^{\mathrm{m}}$ such that
\[
\Proj^M(A) = \bigcup_{S \in \calF} D_\dag(S).
\]
Then if $Q$ is an $M$-graded $A$-module the following conditions are equivalent:
 \begin{enumerate}
  \item $Q$ is negligible;
  \item for any $S \in \calF$ and any $q \in Q$, there exists $s \in S$ such that $s \cdot q=0$.
 \end{enumerate}
\end{lem}

\begin{proof}
 Since $\Proj^M(A)$ is covered by the affine open subschemes $(D_\dag(S) : S \in \calF)$, we have $\widetilde{Q}=0$ if and only if $Q_{(S)}=0$ for any $S \in \calF$. Now since $\calF \subset \calF^{\mathrm{max}}$, for any $S \in \calF$ the ring $A_S$ has invertible elements of all degrees, so that $Q_{(S)}=0$ if and only if $Q_S=0$. Finally, it is clear from the definitions that $Q_S=0$ if and only if for any $q \in Q$ there exists $s \in S$ such that $s \cdot q=0$.
\end{proof}

\begin{rmk}
Let us make the condition in Lemma~\ref{lem:negligible} more explicit in some cases considered in Remark~\ref{rmk:max-cover}.

\begin{enumerate}
 \item 
 First, assume that $A_0$ is noetherian, that $A$ is generated by $A_1$ as an $A_0$-algebra, and moreover that $A_1$ is finite as an $A_0$-module. (This is the setting considered in~\cite[\href{https://stacks.math.columbia.edu/tag/01YR}{Tag 01YR}]{stacks-project}.) For $N \in \Z$ we set $A_{\geq N} = \bigoplus_{m \geq N} A_m$. Then $Q$ is negligible if and only if for any $q \in Q$ there exists $N>0$ such that $A_{\geq N} \cdot q=0$. (Since $A_{\geq N} = (A_{\geq 1})^N$, this condition is equivalent to $Q$ being $A_{\geq 1}$-power torsion in the sense of~\cite[\href{https://stacks.math.columbia.edu/tag/05E6}{Tag 05E6}]{stacks-project}.) In fact, choose $(a_i : i \in I)$ as in Remark~\ref{rmk:max-cover}, with $I$ finite; then we can choose $\calF=\{(a_i)^{\N} : i \in I\}$. If $Q$ satisfies our condition, then it is negligible by Lemma~\ref{lem:negligible}. On the other hand, assume that $Q$ is negligible, and let $q \in Q$. By the lemma, for any $i \in I$ there exists $n_i \in \Z_{>0}$ such that $(a_i)^{n_i} \cdot q=0$. Then if $N = 1 + \sum_{i \in I} (n_i-1)$, we have
 \[
 A_{\geq N} \subset \sum_{i \in I} A \cdot (a_i)^{n_i},
 \]
 %any element in $A_{\geq N}$ is a linear combination (with coefficients in $A_0$) of products of the $a_i$'s, where for some $i$ the element $a_i$ appears at least $a_i$ times. Hence we have $A_{\geq N} \cdot q=0$, 
so that $A_{\geq N} \cdot q=0$.
 \item
 Now, consider the setting of~\S\ref{ss:flag}, and recall the notation of Remarks~\ref{rmk:flag-covering} and~\ref{rmk:max-cover}. For $N \in \Z$ we denote by $\bbX_{\geq N} \subset \bbX$ the submonoid consisting of elements $\lambda$ which satisfy $\langle \lambda, \alpha^\vee \rangle \geq N$ for any simple root $\alpha$. (In particular, we have $\bbX_{\geq 0} = \bbX_+$, and $\bbX_{\geq 1} = \bbX_{++}$.) We also denote by $A_{\geq N}$ the sum of the components in $A$ whose degrees belong to $\bbX_{\geq N}$. We will use similar notation for $\widetilde{\bbX}$ and $\widetilde{A}$.
 Then $Q$ is negligible if and only if for any $q \in Q$ there exists $N>0$ such that $A_{\geq N} \cdot q = 0$. In fact, 
% recall the notation of Remark~\ref{rmk:flag-covering}, and let $\widetilde{\bbX}_0$ be as in~\eqref{eqn:X-sc} for the group $\widetilde{\bG}$. Let also $\widetilde{A}_1$ and $\widetilde{A}_1$ be as in the second proof of Proposition~\ref{prop:flag-var} (for the group $\widetilde{G}$). Then 
 in Lemma~\ref{lem:negligible} we can take $\calF = \{(g_\sigma)^\N : \sigma \in E\}$. Since each $g_\sigma$ belongs to $A_{\geq 1}$, it is clear from this lemma that if $Q$ satisfies our condition, then it is negligible. On the other hand assume that $Q$ is negligible, and fix $q \in Q$. By the lemma, for any $\sigma \in E$ there exists $N_\sigma$ such that $(g_\sigma)^{N_\sigma} \cdot q=0$. 
 %By construction, for any $\sigma$ we have $g_\sigma = (f_\sigma)^{n_\sigma}$ for some element $f_\sigma \in \widetilde{A}$ and some integers $n_\sigma \geq 1$. Set 
 If $N = 1+ \sum_\sigma (N_\sigma n_\sigma - 1)$, we have 
 \begin{equation}
 \label{eqn:inclusion-negligible-flag}
 A_{\geq N} \subset \sum_{\sigma \in E} A \cdot (g_\sigma)^{N_\sigma},
 \end{equation}
 which will imply the claim. In fact, by the surjectivity of the maps~\eqref{eqn:mult-line-bundles-flag} we have
 \[
 \widetilde{A}_{\geq 1} \subset \sum_{\sigma \in E} \widetilde{A} \cdot f_\sigma,
 \]
 and then
 \[
  \widetilde{A}_{\geq N} = ( \widetilde{A}_{\geq 1})^N \subset \sum_{\sigma \in E} \widetilde{A} \cdot (g_\sigma)^{N_\sigma},
 \]
 which implies~\eqref{eqn:inclusion-negligible-flag}.
% 
% The elements $f_\sigma$ generate $\widetilde{A}$ as a $\widetilde{A}_2$-algebra, and each of them is of degree $\sum_\alpha \varpi_\alpha$; hence $A_{\geq N}$ is spanned by products of the form $a_0 \cdot \prod_\sigma (f_\sigma)^{m_\sigma}$ where $a_0 \in \widetilde{A}_2$ and $\sum_\sigma m_\sigma \geq N$. Here we must have $m_\sigma \geq N_\sigma n_\sigma$ for some $\sigma$, so that
% \[
%  A_{\geq N} \subset \sum_\sigma A \cdot (g_\sigma)^{N_\sigma}.
% \]
%This implies that $A_{\geq N} \cdot q=0$ and finishes the proof.
\item
Consider the setting of~\S\ref{ss:vector-bundles}. Defining $\bbX_{\geq N}$ as above, and then $(A_V)_{\geq N}$ as the sum of the graded components of $A_V$ whose label belongs to $\bbX_{\geq N}$, one checks using Lemma~\ref{lem:restriction-vb-surjective} and~\eqref{eqn:inclusion-negligible-flag} that a graded $A_V$-module $Q$ is negligible if and only if for any $q \in Q$ there exists $N>0$ such that $(A_V)_{\geq N} \cdot q = 0$. Similar comments apply in the setting of~\S\ref{ss:Springer}.
\end{enumerate}
\end{rmk}

\begin{propo}
\label{prop:adjunction}
Assume that~\eqref{eqn:cover-max} is satisfied. Then the composition
\[
\Mod^M(A) \xrightarrow{Q \mapsto \widetilde{Q}} \QCoh(\Proj^M(A)) \longto \Mod(\calO_{\Proj^M(A)})
\]
(where the second functor is the obvious forgetful functor)
is left adjoint to the functor $\Gamma_\bullet$.
\end{propo}

\begin{proof}
To prove the proposition we need to define functorial morphisms
$\varepsilon_{\calQ} : \widetilde{\Gamma_\bullet(\calQ)} \to \calQ$
for $\calQ \in \Mod(\calO_{\Proj^M(A)})$ and
$\eta_Q : Q \to \Gamma_\bullet(\widetilde{Q})$
for $Q \in \Mod^M(A)$, which satisfy the usual zigzag relations (see e.g.~\cite[Proposition~A.1.16]{achar}). Here $\varepsilon$ is provided by Lemma~\ref{lem:morph-tilde-Gamma}. (This does not require any assumption.) To define $\eta$, we observe that for $Q \in \Mod^M(A)$ we have
\[
\Gamma_\bullet(\widetilde{Q}) = \bigoplus_{\alpha \in M} \Gamma(\Proj^M(A), \widetilde{Q}(\alpha)) \cong \bigoplus_{\alpha \in M} \Gamma(\Proj^M(A), \widetilde{Q(\alpha)}),
\]
where the isomorphism follows from Corollary~\ref{coro:twistings-line-bundles}\eqref{it:max-rel-2}. The desired morphism is therefore provided by Fact~\ref{fact:morp-graded-rings-modules}\eqref{it:morp-graded-modules}.

We leave it to the reader to check that these morphisms indeed satisfy the zigzag relations.
\end{proof}

Recall that if $\mathcal{L}$ is an invertible sheaf on a scheme $X$ and $s \in \Gamma(X,\mathcal{L})$ is a global section, we have an open subscheme $X_s \subset X$ defined by the nonvanishing of $s$, see~\cite[\href{https://stacks.math.columbia.edu/tag/01CY}{Tag 01CY}]{stacks-project}.

\begin{lem}
\label{lem:open-section}
Assume that~\eqref{eqn:cover-max} is satisfied.
Let $f \in A$ be homogeneous and relevant, and let $\tilde{f}$ be the image of $f$ in $\Gamma(\Proj^M(A), \calO_{\Proj^M(A)}(\deg(f))$ (cf.~Fact~\ref{fact:morp-graded-rings-modules}\eqref{it:morp-graded-rings}). 
%If $\scO_{\Proj^M(A)}(\deg(f))$ is a line bundle (OR OTHER ASSUMPTION?), then w
We have
\[
\Proj^M(A)_{\tilde f} = D_\dag(f).
\]
\end{lem}

\begin{proof}
We can assume that $f$ is nonzero.
Our assumptions imply that $\calO_{\Proj^M(A)}(\deg(f))$ is an invertible sheaf (see Corollary~\ref{coro:twistings-line-bundles}), so that $\Proj^M(A)_{\tilde f}$ is well defined. In view of our assumption, to prove the statement it suffices to prove that for any $S \in \calF^{\mathrm{m}}_A$ we have
\[
D_\dag(S)_{\tilde{f}_{|D_\dag(S)}} = D_\dag(S) \cap D_\dag(f).
\]
%where $\tilde{f}_S$ is the restriction of $\tilde{f}$ to $D_\dag(S)$. 
Fix such an $S$. As in the proof of Proposition~\ref{prop:twisting-inversible-open}\eqref{it:twisting-inversible-open-1}, there exist $a,s \in S$ such that $\deg(a)-\deg(s)=\deg(f)$, and then multiplication by $\frac{a}{s}$ induces an isomorphism of sheaves
\[
\calO_{D_\dag(S)} \simto \calO_{\Proj^M(A)}(\deg(f))_{| D_\dag(S)}.
\]
The inverse image of $\tilde{f}_{|D_\dag(S)}$ under this isomorphism is $\frac{fs}{a}$, which implies that
\[
D_\dag(S)_{\tilde{f}_{|D_\dag(S)}} = \Spec \left( (A_{(S)})_{\frac{fs}{a}} \right).
\]
By Proposition~\ref{prop:magical}\eqref{it:magic-2}, the right-hand side identifies with $\Spec(A_{(S \cdot f^\N)})$, i.e.~with $D_\dag(S) \cap D_\dag(f)$, which finishes the proof.
%COMPLETE! Cf.~\cite[\href{https://stacks.math.columbia.edu/tag/01MV}{Tag 01MV}]{stacks-project}.
\end{proof}

\begin{propo}
\label{prop:tilde-Gamma-isom}
Assume that~\eqref{eqn:cover-max} is satisfied, and moreover that $\Proj^M(A)$ is quasi-compact. Then for any $\calQ$ in $\QCoh(\Proj^M(A))$, the morphism
$\widetilde{\Gamma_\bullet(\calQ)} \to \calQ$
of Lemma~\ref{lem:morph-tilde-Gamma} is an isomorphism.
\end{propo}

\begin{proof}
Since the sheaves under consideration are quasi-coherent, and since~\eqref{eqn:cover-max} holds, to prove the statement it suffices to prove that for any maximally relevant $f \in A$ the morphism
\[
\Gamma(D_\dag(f), \widetilde{\Gamma_\bullet(\calQ)}) \longto \Gamma(D_\dag(f), \calQ)
\]
induced by the morphism of Lemma~\ref{lem:morph-tilde-Gamma} is an isomorphism. Now by definition the left-hand side identifies with $\Gamma_\bullet(\calQ)_{(f)}$. Considering the $\Z$-graded module
\[
Q:= \bigoplus_{n \in \Z} \Gamma \bigl( \Proj^M(A), \calQ(n \cdot \deg(f)) \bigr)
\]
over the $\Z$-graded ring
\[
\bigoplus_{n \in \N} \Gamma \bigl( \Proj^M(A), \calO_{\Proj^M(A)}(n \cdot \deg(f) \bigr),
\]
and denoting by $\tilde{f}$ the image of $f$ in $\Gamma(\Proj^M(A), \calO_{\Proj^M(A)}(\deg(f))$,
then we have
$\Gamma_\bullet(\calQ)_{(f)} = Q_{(\tilde{f})}$.
Since $\Proj^M(A)$ is quasi-compact (by assumption) and quasi-separated (see Lemma~\ref{lem:Proj-qsep}), by~\cite[\href{https://stacks.math.columbia.edu/tag/01PW}{Tag 01PW}]{stacks-project} the right-hand side identifies with
$\Gamma(\Proj^M(A)_{\tilde{f}}, \calQ)$,
i.e.~with $\Gamma(D_\dag(f), \calQ)$ by Lemma~\ref{lem:open-section}, which finishes the proof.
\end{proof}

\begin{rmk}
\label{rmk:tilde-Gamma-isom-torsor}
Under the additional assumption that $M$ is a free abelian group, one can give an alternative proof of Proposition~\ref{prop:tilde-Gamma-isom} as follows. 
%By Proposition~\ref{prop:max-torsor}, 
As explained in the comments following~\eqref{eqn:geometric-quotient}, the natural morphism
$\Spec(A) \smallsetminus V(A_\dag) \to \Proj^M(A)$ is affine, hence quasi-compact,
%and separated (see~\cite[\href{https://stacks.math.columbia.edu/tag/01S7}{Tag 01S7}]{stacks-project}), 
so that the scheme $A \smallsetminus V(A_\dag)$ is quasi-compact. Hence the open immersion $j : \Spec(A) \smallsetminus V(A_\dag) \to \Spec(A)$ is quasi-compact and separated, so that the pushforward functor $j_*$ preserves quasi-coherent sheaves, see~\cite[\href{https://stacks.math.columbia.edu/tag/01LC}{Tag 01LC}]{stacks-project}. Under the identifications considered in Remark~\ref{rmk:tilde-torsor}, the functor $Q \mapsto \widetilde{Q}$ corresponds to $j^*$, while the functor $\Gamma_\bullet$ corresponds to $j_*$, and the statement of the proposition becomes the familiar fact that the adjunction morphism $j^* j_* \to \id$ is an isomorphism.
\end{rmk}

\begin{cor}
\label{cor:tilde-Gamma-equiv-quotient}
Assume that~\eqref{eqn:cover-max} is satisfied, and moreover that $\Proj^M(A)$ is quasi-compact. Then the functor
\[
\mathsf{L} : \Mod^M(A) / \Mod^M(A)_{\mathrm{neg}} \longto \QCoh(\Proj^M(A))
\]
and the composition
\[
\QCoh(\Proj^M(A)) \longto \Mod(\calO_{\Proj^M(A)}) \xrightarrow{\Gamma_\bullet} \Mod^M(A) \longto \Mod^M(A) / \Mod^M(A)_{\mathrm{neg}}
\]
(where the first arrow is the obvious embedding and the third one is the quotient functor) are mutually inverse equivalences of categories.
\end{cor}

This corollary follows immediately from Propositions~\ref{prop:adjunction} and~\ref{prop:tilde-Gamma-isom}, in view of the following general fact. (For the notion of kernel of an exact functor between abelian categories, see~\cite[\href{https://stacks.math.columbia.edu/tag/02MR}{Tag 02MR}]{stacks-project}.)

\begin{lem}
Let $\mathsf{A}$, $\mathsf{B}$ be abelian categories, and let
$L : \mathsf{A} \longto \mathsf{B}$ and $R : \mathsf{B} \longto \mathsf{A}$
be functors. Assume that
\begin{enumerate}
\item
$L$ is left adjoint to $R$;
\item
$L$ is exact;
\item
\label{it:assum-quotient-equiv-3}
the adjunction morphism $LR \to \id$ is an isomorphism.
\end{enumerate}
Then $L$ factors through an equivalence of categories
$\overline{L} : \mathsf{A} / \ker(L) \simto \mathsf{B}$,
whose quasi-inverse is the composition of $R$ with the quotient functor $\mathsf{A} \longto \mathsf{A} / \ker(L)$.
\end{lem}

\begin{proof}
Let us denote by $\pi : \mathsf{A} \to \mathsf{A} / \ker(L)$ the quotient functor.
By~\cite[\href{https://stacks.math.columbia.edu/tag/02MS}{Tag 02MS}]{stacks-project}, $L$ factors through a functor $\overline{L} : \mathsf{A} / \ker(L) \to \mathsf{B}$. Then our assumption~\eqref{it:assum-quotient-equiv-3} shows that $\overline{L} \circ (\pi R) \cong \id$. On the other hand, using adjunction we have a canonical morphism
\[
\pi \longto \pi R L = \pi R \overline{L} \pi.
\]
We claim that this morphism is an isomorphism, which will conclude the proof in view of the fact that if $F,G : \mathsf{A}/\ker(L) \to \mathsf{B}$ are two functor, each morphism of functors $F \pi \to G \pi$ is induced by a unique morphism of functors $F \to G$.

To prove the claim it suffices to prove that for $X \in \mathsf{B}$, the kernel and cokernel of the adjunction morphism $X \to RLX$ belong to $\ker(L)$. Since $L$ is exact, this is equivalent to showing that the image under $L$ of this morphism is an isomorphism, which follows from the zigzag relation and our assumption~\eqref{it:assum-quotient-equiv-3}.
\end{proof}

\begin{rmk}
\label{rmk:tilde-Gamma-equiv-quotient-torsor}
As in Remark~\ref{rmk:tilde-Gamma-isom-torsor}, under the additional assumption that $M$ is a free abelian group, one can give an alternative proof of Corollary~\ref{cor:tilde-Gamma-equiv-quotient} by noticing that the functor $j^*$ induces an equivalence between the category of $\mathrm{D}_{\Spec(A_0)}(M)$-equivariant quasi-coherent sheaves on $\Spec(A) \smallsetminus V(A_\dag)$ and the Serre quotient of the category of $\mathrm{D}_{\Spec(A_0)}(M)$-equivariant quasi-coherent sheaves on $\Spec(A)$ by the Serre subcategory of sheaves supported set-theoretically on $V(A_\dag)$, i.e.~whose restriction to $\Spec(A) \smallsetminus V(A_\dag)$ vanishes.
\end{rmk}

We now consider the case when $A$ is noetherian. Recall that in this case the scheme $\Proj^M(A)$ is noetherian (in particular, quasi-compact), see~\S\ref{ss:twisting-sheaves}. We consider the category $\Mod^M_{\mathrm{fg}}(A)$ of finitely generated $M$-graded $A$-module, its Serre subcategory $\Mod^M_{\mathrm{fg}}(A)_{\mathrm{neg}}$ of objects which are negligible modules, and the category $\Coh(\Proj^M(A))$ of coherent sheaves on $\Proj^M(A)$. Recall (see Lemma~\ref{lem:tilde-noetherian-coherent}) that in this case the functor $Q \mapsto \widetilde{Q}$ restricts to a functor $\Mod^M_{\mathrm{fg}}(A) \to \Coh(\Proj^M(A))$, which must factor through a functor
\begin{equation*}
%\label{eqn:tilde-Coh}
\mathsf{L}_{\Coh} : \Mod^M_{\mathrm{fg}}(A) / \Mod^M_{\mathrm{fg}}(A)_{\mathrm{neg}} \longto \Coh(\Proj^M(A)).
\end{equation*}

\begin{propo}
\label{prop:quotient-equiv-noetherian}
If $A$ is noetherian, the functor $\mathsf{L}_{\Coh}$
%~\eqref{eqn:tilde-Coh} 
is an equivalence of categories.
\end{propo}

\begin{proof}
We have a commutative diagram
\[
\xymatrix{
\Mod^M_{\mathrm{fg}}(A) / \Mod^M_{\mathrm{fg}}(A)_{\mathrm{neg}} \ar[r]^-{\mathsf{L}_{\Coh}} \ar[d] & \Coh(\Proj^M(A)) \ar[d] \\
\Mod^M(A) / \Mod^M(A)_{\mathrm{neg}} \ar[r]^-{\mathsf{L}} & \QCoh(\Proj^M(A))
}
\]
where the lower horizontal arrow is known to be an equivalence (see Corollary~\ref{cor:tilde-Gamma-equiv-quotient}) and the vertical arrows are fully faithful. It follows that $\mathsf{L}_{\Coh}$ is fully faithful.

To prove essential surjectivity, we consider $\mathscr{F} \in \Coh(\Proj^M(A))$, and the $M$-graded $A$-module $Q=\Gamma_\bullet(\mathscr{F})$. Since $A$ is noetherian, $Q$ is the filtered colimit of its finitely generated $M$-graded $A$-submodules; in other words there exists a filtered set $I$ and finitely generated $M$-graded $A$-submodules $Q_i \subset Q$ such that
%\[
%Q_0 \subset Q_1 \subset Q_2 \subset \cdots \subset Q
%\]
%such that each $Q_n$ is finitely generated and 
$Q=\mathrm{colim}_i Q_i$. By exactness of the functor $P \mapsto \widetilde{P}$, and since this functor commutes with colimits (see~\S\ref{ss:sheaves-gr-modules}), each $\widetilde{Q_i}$ is a coherent subsheaf of $\mathscr{F} = \widetilde{Q}$, and we have $\mathscr{F} = \colim_i \, \widetilde{Q_i}$. As in~\cite[\href{https://stacks.math.columbia.edu/tag/01Y8}{Tag 01Y8}]{stacks-project}, this implies that $\mathscr{F} = \widetilde{Q_i}$ for some $i$, hence that $\mathscr{F}$ belongs to the essential image of $\mathsf{L}_{\Coh}$.
%
%we deduce a chain of inclusions
%\[
%\widetilde{Q_0} \subset \widetilde{Q_1} \subset \widetilde{Q_2} \subset \cdots \subset \widetilde{Q}.
%\]
%We claim that the natural morphism
%\[
%\mathrm{colim}_n \, \widetilde{Q_n} \to \widetilde{Q}
%\]
%is an isomorphism. In fact, both sides are quasi-coherent, so that to justify this claim it suffices to show that for any $S \in \calF$ the natural morphism
%\[
%\Gamma(D_\dag(S), \mathrm{colim}_n \, \widetilde{Q_n}) \to \Gamma(D_\dag(S), \widetilde{Q})= Q_{(S)}
%\]
%is an isomorphism. Here the left-hand side identifies with
%\[
%\mathrm{colim}_n \Gamma(D_\dag(S), \widetilde{Q_n}) = \mathrm{colim}_n \, (Q_n)_{(S)},
%\]
%see e.g.~\cite[\href{https://stacks.math.columbia.edu/tag/009F}{Tag 009F}]{stacks-project}, and then the claim follows from the fact that localization commutes with colimits.
%It follows that
%\[
%\mathscr{F} = \widetilde{Q} = \mathrm{colim}_n \widetilde{Q_n}.
%\]
%Since $\mathscr{F}$ is coherent, by~\cite[\href{https://stacks.math.columbia.edu/tag/01Y8}{Tag 01Y8}]{stacks-project} this colimit stabilizes. Hence for $n \gg 0$ we have $\mathscr{F} = \widetilde{Q_n}$, which shows that $\mathscr{F}$ belongs to the essential image of $\mathsf{L}_{\Coh}$.
\end{proof}

\begin{rmk}
\label{rmk:equivalence-equiv}
 In the setting of Remark~\ref{rmk:tilde-equiv} we similarly obtain an equivalence of abelian categories
 \[
 \Mod^{M,H}(A) / \Mod^{M,H}(A)_{\mathrm{neg}} \longto \QCoh^H(\Proj^M(A))
 \]
 and, in case $A$ is noetherian, an equivalence of abelian categories
  \[
 \Mod_{\mathrm{fg}}^{M,H}(A) / \Mod_{\mathrm{fg}}^{M,H}(A)_{\mathrm{neg}} \longto \Coh^H(\Proj^M(A))
 \]
 where we use obvious notation in the left-hand side.
\end{rmk}

%--------------------------------------------------------
\subsection{Derived categories}
%--------------------------------------------------------

We come back to the general setting of~\S\ref{ss:max-relevant}, and consider for $? \in \{+,-,\mathrm{b}\}$ the derived categories $D^? \QCoh(\Proj^M(A))$ and $D^? \Mod^M(A)$. We will denote by $D^?_{\mathrm{neg}} \Mod^M(A)$ the full triangulated subcategory of the latter category consisting of complexes all of whose cohomology objects are negligible. In the following statement we use the Verdier quotient of a triangulated category by a full triangulated subcategory; for this notion, we refer to~\cite[\href{https://stacks.math.columbia.edu/tag/05RA}{Tag 05RA}]{stacks-project}.

\begin{propo}
\label{prop:quotient-cat-derived}
Assume that~\eqref{eqn:cover-max} is satisfied, and moreover that $\Proj^M(A)$ is quasi-compact. Then if $? \in \{+,-,\mathrm{b}\}$ the functor $Q \mapsto \widetilde{Q}$ induces an equivalence of triangulated categories
\[
D^? \Mod^M(A) / D^?_{\mathrm{neg}} \Mod^M(A) \simto D^? \QCoh(\Proj^M(A)).
\]
\end{propo}

\begin{proof}
Since the functor $Q \mapsto \widetilde{Q}$ is exact, it induces a triangulated functor
\[
D^? \Mod^M(A) \to D^? \QCoh(\Proj^M(A))
\]
on derived categories. 
%Now it is easy to see that $\Db_{\mathrm{neg}} \Mod^M(A)$ is generated as a triangulated category by the complexes consisting of a single negligible $M$-graded $A$-module in degree $0$. These objects are clearly sent to $0$, 
This functor sends objects in $D^?_{\mathrm{neg}} \Mod^M(A)$ to complexes all of whose cohomology objects are trivial, i.e.~to the zero object,
so that our functor factors through a triangulated functor
$D^? \Mod^M(A) / D^?_{\mathrm{neg}} \Mod^M(A) \to D^?\QCoh(\Proj^M(A))$.
To check that this functor is an equivalence, we note that by~\cite[Theorem~3.2]{miyachi} there exists a canonical equivalence of triangulated categories
\[
D^? \Mod^M(A) / D^?_{\mathrm{neg}} \Mod^M(A) \simto D^? (\Mod^M(A) / \Mod^M(A)_{\mathrm{neg}}).
\]
This reduces the statement to the abelian case, which was proved in Corollary~\ref{cor:tilde-Gamma-equiv-quotient}.
\end{proof}

The following statement is the version of Proposition~\ref{prop:quotient-cat-derived} for coherent sheaves, in the case of noetherian rings. Here we denote by $D^?_{\mathrm{neg}} \Mod_{\mathrm{fg}}^M(A)$ the full subcategory of $D^? \Mod_{\mathrm{fg}}^M(A)$ consisting of complexes all of whose cohomology objects are negligible. The proof is the same, simply replacing the reference to Corollary~\ref{cor:tilde-Gamma-equiv-quotient} by a reference to Proposition~\ref{prop:quotient-equiv-noetherian}.

\begin{propo}
\label{prop:quotient-cat-derived-noeth}
Assume that $A$ is noetherian, and that~\eqref{eqn:cover-max} is satisfied. Then if $? \in \{+,-,\mathrm{b}\}$ the functor $Q \mapsto \widetilde{Q}$ induces an equivalence of triangulated categories
\[
D^? \Mod_{\mathrm{fg}}^M(A) / D^?_{\mathrm{neg}} \Mod_{\mathrm{fg}}^M(A) \simto D^? \Coh(\Proj^M(A)).
\]
\end{propo}

\begin{rmk}
As in Remark~\ref{rmk:equivalence-equiv}, in the setting of Remark~\ref{rmk:tilde-equiv} one obtains similar equivalences for derived categories of \emph{equivariant} modules and quasi-coherent sheaves.
\end{rmk}

%--------------------------------------------------------
\subsection{On a lemma by Arkhipov--Bezrukavnikov}
%--------------------------------------------------------

Continue with the setting of Proposition~\ref{prop:quotient-cat-derived-noeth}, and denote by $\Mod_{\mathrm{fg},\mathrm{fr}}^M(A)$ the additive category of \emph{free} $M$-graded $A$-modules. Since any finitely generated $M$-graded $A$-module is a quotient of an object of $\Mod_{\mathrm{fg},\mathrm{fr}}^M(A)$, we have a canonical equivalence of triangulated categories
\begin{equation}
\label{eqn:equiv-Kbfree-minus}
K^- \Mod_{\mathrm{fg},\mathrm{fr}}^M(A) \simto D^- \Mod_{\mathrm{fg}}^M(A).
\end{equation}
In view of Proposition~\ref{prop:quotient-cat-derived-noeth}, denoting by $K_{\mathrm{neg}}^- \Mod_{\mathrm{fg},\mathrm{fr}}^M(A)$ the full subcategory of $K^- \Mod_{\mathrm{fg},\mathrm{fr}}^M(A)$ consisting of complexes all of whose cohomology objects are negligible, we deduce an equivalence of categories
\[
K^- \Mod_{\mathrm{fg},\mathrm{fr}}^M(A) / K_{\mathrm{neg}}^- \Mod_{\mathrm{fg},\mathrm{fr}}^M(A) \simto D^- \Coh(\Proj^M(A)).
\]
In~\cite[Sublemma~1]{ab}, the authors state a similar claim for \emph{bounded} homotopy categories.\footnote{Note that for a triangulated category $\mathsf{A}$ and full subcategories $\mathsf{B}$, $\mathsf{C}$, it is not true in general that the natural functor $\mathsf{B} / (\mathsf{B} \cap \mathsf{C}) \to \mathsf{A} / \mathsf{C}$ is fully faithful; the bounded case therefore does not immediately follow from the unbounded case.} The proof is very sketchy, and it is not clear to us if it really applies in the stated generality. Here we provide a complete argument for this claim (under additional assumptions), kindly explained to us by R.~Bezrukavnikov.

\begin{propo}
\label{prop:sublemma-AB}
Assume that $A$ is a finitely generated $\bk$-algebra for some field $\bk$, that $M$ is a free abelian group, and that~\eqref{eqn:cover-max} is satisfied. Consider the bounded homotopy category $K^{\mathrm{b}} \Mod_{\mathrm{fg},\mathrm{fr}}^M(A)$, and the full subcategory $K_{\mathrm{neg}}^{\mathrm{b}} \Mod_{\mathrm{fg},\mathrm{fr}}^M(A)$ of complexes all of whose cohomology objects are negligible. The functor $Q \mapsto \widetilde{Q}$ induces a fully faithful functor
\[
K^{\mathrm{b}} \Mod_{\mathrm{fg},\mathrm{fr}}^M(A) / K_{\mathrm{neg}}^{\mathrm{b}} \Mod_{\mathrm{fg},\mathrm{fr}}^M(A) \longto D^{\mathrm{b}} \Coh(\Proj^M(A)).
\]
\end{propo}

\begin{proof}
As in Remark~\ref{rmk:tilde-torsor}, under the present assumptions, setting $T = \mathrm{D}_{\Spec(\bk)}(M)$, the functor $Q \mapsto \widetilde{Q}$ identifies with the functor
\[
\Mod_{\mathrm{fg}}^M(A) = \QCoh^{T}(\Spec(A)) \xrightarrow{j^*} \Coh^{T}(\Spec(A) \smallsetminus V(A_\dag)) \cong \Coh(\Proj^M(A)),
\]
where $j : \Spec(A) \smallsetminus V(A_\dag) \to \Spec(A)$ is the embedding. What we have to show is therefore that for any bounded complexes $\calF,\calF'$ of objects in $\Mod_{\mathrm{fg},\mathrm{fr}}^M(A)$ the morphism
%functor $j^*$ induces an isomorphism
\[
\phi : \Hom_{K^{\mathrm{b}} \Mod_{\mathrm{fg},\mathrm{fr}}^M(A) / K_{\mathrm{neg}}^{\mathrm{b}} \Mod_{\mathrm{fg},\mathrm{fr}}^M(A)}(\calF, \calF') \to \Hom_{D^{\mathrm{b}} \Coh^{T}(\Spec(A) \smallsetminus V(A_\dag))}(j^* \calF, j^* \calF').
\]
induced by the functor $j^*$ an isomorphism.
%where in the left-hand side we consider morphisms in the quotient category
%\[
%K^{\mathrm{b}} \Mod_{\mathrm{fg},\mathrm{fr}}^M(A) / K_{\mathrm{neg}}^{\mathrm{b}} \Mod_{\mathrm{fg},\mathrm{fr}}^M(A).
%\]
For this we will construct a morphism $\psi$ in the reverse direction, and check that $\phi$ and $\psi$ are inverse to each other. In the course of the proof we will use the obvious fact that the natural functor
\begin{equation}
\label{eqn:functor-Kb-Db}
K^{\mathrm{b}} \Mod_{\mathrm{fg},\mathrm{fr}}^M(A) \to D^{\mathrm{b}} \Mod_{\mathrm{fg}}^M(A)
\end{equation}
is fully faithful.

Fix $\calF,\calF'$ as above, and
consider a morphism $f : j^* \calF \to j^* \calF'$ in $D^{\mathrm{b}} \Coh^{T}(\Spec(A) \smallsetminus V(A_\dag))$. By Proposition~\ref{prop:quotient-cat-derived-noeth}, this morphism can be represented by a diagram
\begin{equation}
\label{eqn:diag-morph-quotient-Kb}
\calF \xrightarrow{g} \calF'' \xleftarrow{h} \calF'
\end{equation}
where $\calF'' \in D^{\mathrm{b}} \Mod^M_{\mathrm{fg}}(A)$ and $h$ is a morphism whose cone $\calG$ has all of its cohomology objects supported set-theoretically on $V(A_\dag)$. Then, for some $n \gg 0$, $\calG$ can be represented by a bounded complex of $M$-graded $A$-modules all of whose components are annihilated by $(A_\dag)^n$ (see e.g.~\cite[Proposition~A.1]{br}). Consider a finite-dimensional graded subspace $E \subset (A_\dag)^n$ that generates this ideal, and the Koszul complex $\mathcal{K}$ of the multiplication morphism $E \otimes_\bk A \to A$, see~\cite[\href{https://stacks.math.columbia.edu/tag/0621}{Tag 0621}]{stacks-project}. By definition, this is a bounded complex of free $M$-graded $A$-modules, concentrated in nonnegative degrees, and whose degree-$0$ component is $A$. Moreover, by~\cite[\href{https://stacks.math.columbia.edu/tag/0663}{Tag 0663}]{stacks-project} the restriction of this complex to $\Spec(A) \smallsetminus V(A_\dag)$ is acyclic. 
%(This follows e.g.~from~\cite[\href{https://stacks.math.columbia.edu/tag/0663}{Tag 0663}]{stacks-project}.) 
Denote by $\mathcal{C}$ the cokernel of the natural embedding of complexes $A \to \mathcal{K}$; then we have a morphism $\mathcal{C}[-1] \to A$ whose cone (namely, $\mathcal{K}$) is supported on $V(A_\dag)$, hence a morphism
$\mathcal{F} \otimes_A \mathcal{C}[-1] \to \calF$
with the same property. We claim that the composition
\[
\mathcal{F} \otimes_A \mathcal{C}[-1] \to \calF \xrightarrow{g} \calF'' \to \calG
\]
vanishes. This will imply that the composition of the first two morphisms factors through a morphism
$\mathcal{F} \otimes_A \mathcal{C}[-1] \to \calF'$
(in $D^{\mathrm{b}} \Mod_{\mathrm{fg}}^M(A)$, or equivalently in $K^{\mathrm{b}} \Mod_{\mathrm{fg},\mathrm{fr}}^M(A)$);
then the diagram
\[
\calF \leftarrow \calF \otimes_A \mathcal{C}[-1] \to \calF'
\]
will define the desired morphism $\psi(f)$ in $K^{\mathrm{b}} \Mod_{\mathrm{fg},\mathrm{fr}}^M(A) / K_{\mathrm{neg}}^{\mathrm{b}} \Mod_{\mathrm{fg},\mathrm{fr}}^M(A)$.

To prove the claim, it suffices to notice that the morphism
\[
\Hom_{D^{\mathrm{b}} \Mod_{\mathrm{fg}}^M(A)}(\calF \otimes_A \mathcal{K}, \calG) \to \Hom_{D^{\mathrm{b}} \Mod_{\mathrm{fg}}^M(A)}(\calF, \calG)
\]
%(where we consider morphisms in $D^{\mathrm{b}} \Mod_{\mathrm{fg}}^M(A)$)
is surjective. This follows from the isomorphism
\[
\Hom_{D^{\mathrm{b}} \Mod_{\mathrm{fg}}^M(A)}(\calF \otimes_A \mathcal{K}, \calG) \cong \Hom_{D^{\mathrm{b}} \Mod_{\mathrm{fg}}^M(A)}(\calF, \calG \otimes_A \mathcal{K}^\vee)
\]
where $\mathcal{K}^\vee$ is the dual complex of $A$-modules, and the fact that $\calG$ is a direct summand in $\calG \otimes_A \mathcal{K}^\vee$ by our choice of $n$.
% this complex can be represented by a complex all of whose component are annihilated by $I^n$.

It is clear from construction that $\phi \circ \psi = \id$. On the other hand, any morphism $f' : \calF \to \calF'$ in the quotient category can be represented by a diagram~\eqref{eqn:diag-morph-quotient-Kb} where now $g$, $h$ are morphisms in $K^{\mathrm{b}} \Mod_{\mathrm{fg},\mathrm{fr}}^M(A)$. For the construction of $\psi(\phi(f'))$ we can take the images of these morphisms in $D^{\mathrm{b}} \Mod_{\mathrm{fg}}^M(A)$; we deduce a commutative diagram
\[
\xymatrix@R=0.1cm{
& \calF \otimes_A \mathcal{C}[-1] \ar[rd]^-{l} \ar[ld]_-{k} & \\
\calF \ar[rd]_-{g} & & \calF' \ar[ld]^-{h} \\
& \calF'' &
}
\]
in $D^{\mathrm{b}} \Mod_{\mathrm{fg}}^M(A)$,
in which $k$ and $h$ have their cones supported on $V(A_\dag)$. Since the functor~\eqref{eqn:functor-Kb-Db}
%\[
%K^{\mathrm{b}} \Mod_{\mathrm{fg},\mathrm{fr}}^M(A) \to D^{\mathrm{b}} \Mod_{\mathrm{fg}}^M(A)
%\]
is fully faithful, we can regard this diagram as a diagram in $K^{\mathrm{b}} \Mod_{\mathrm{fg},\mathrm{fr}}^M(A)$, and its commutativity shows that $\psi(\phi(f'))=f'$.
\end{proof}

\begin{rmk}
\begin{enumerate}
\item
Consider the setting of Remark~\ref{rmk:tilde-equiv}, with $H$ a linearly reductive algebraic group over $\bk$. Then one can consider the full subcategory $\Mod_{\mathrm{fg},\mathrm{fr}}^{M,H}(A)$ of the category $\Mod_{\mathrm{fg}}^{M,H}(A)$ of objects which are sums of objects of the form $V \otimes A(\alpha)$ where $V$ is a finite-dimensional $H$-module and $\alpha \in M$. Let $K_{\mathrm{neg}}^{\mathrm{b}} \Mod_{\mathrm{fg},\mathrm{fr}}^{M,H}(A)$ be the full triangulated subcategory of $K^{\mathrm{b}} \Mod_{\mathrm{fg},\mathrm{fr}}^{M,H}(A)$ whose objects are the complexes all of whose cohomology objects are supported on $V(A_\dag)$. Then the same proof as for
Proposition~\ref{prop:sublemma-AB} shows that the natural functor
\[
K^{\mathrm{b}} \Mod_{\mathrm{fg},\mathrm{fr}}^{M,H}(A) / K_{\mathrm{neg}}^{\mathrm{b}} \Mod_{\mathrm{fg},\mathrm{fr}}^{M,H}(A) \longto D^{\mathrm{b}} \Coh^H(\Proj^M(A))
\]
is fully faithful
\item
Another possible approach to this question, also suggested by R. Bezrukavnikov, would be to use the standard fact that the perfect derived category of an open subscheme $U \subset X$ is the Verdier quotient of the perfect derived category of $X$ by the subcategory of complexes supported on $X \smallsetminus U$, suitably generated to quotient stacks. (For a result of this form, see e.g.~\cite{tt}.) In the presence of a contracting $\mathbb{G}_{\mathrm{m}}$-action, any equivariant vector bundle is free, which relates this statement with Proposition~\ref{prop:sublemma-AB}.
\end{enumerate}
\end{rmk}

%   As in the absolute case, there is an exact functor $\mathcal{F} \to \widetilde{\mathcal{F}}$ from the category of $M$-graded quasi-coherent $\mathcal{A}$-modules to the category of quasi-coherent $\mathcal{O}_X$-modules. In particular, we can define the $\mathcal{O}_X$-modules $\mathcal{O}_X (\alpha):=\widetilde{\mathcal{A} (\alpha)}$ for any $\alpha \in M$.

%\medskip

%\section{Flag variety}
%
%
%\exam[Flag variety] \label{exampleGB} Let $k$ be an algebraically closed field and let $G$ be a connected reductive group scheme over $k$. Let $T$ be a maximal split torus and $B$ be a Borel subgroup such that $B = TN $ where $N$ is the unipotent radical. Then $G/N $ is quasi-affine and the ring $A:=\Gamma (G/N , \mathcal{O} _{G/N })$ is canonically $X^*(T)^+$-graded. We have a schematically dominant open immersion $G/N \to \mathrm{Spec} (A)$ (cf. e.g. \cite[§5]{Wa}). Let $M=(X^*(T)^+)^{gp}$ be the group generated by the cancellative monoid $X^*(T)^+$. Then, for example using \cite[page 7]{BS03} and \cite[page 13]{Wa}, the flag variety $G/B$ identifies with $\mathrm{Proj} (A)$ where $A$ is $M$-graded. Investigations and details around this example will be the topic of this section. 
%\xexam 

%%%%%%%%%%%%%%%%%%%%%%%%%%%%%

%\input{referenc}
\end{document}